\documentclass[preprint,sort&compress,12pt]{elsarticle}
\usepackage{CJK}
\usepackage{bbm}
\usepackage{amstext}
\usepackage{amsmath}
\usepackage{mathrsfs}
\usepackage{stmaryrd}
\usepackage{bm}
\usepackage{pstricks}
\usepackage{pst-coil}
\usepackage{pst-3d}
\usepackage{amsfonts}
\usepackage{amssymb}
\usepackage{amsthm}
\makeatletter
\let\ps@pprintTitle\ps@empty
\makeatother
\allowdisplaybreaks
\newcommand{\mm}{\mathrm}

\newcommand{\mmd}{\mathrm{d}}

\newcommand{\be}{\begin{equation}}
\newcommand{\bea}{\begin{equation}\begin{aligned}}
\newcommand{\beas}{\begin{equation*}\begin{aligned}}
\newcommand{\eeas}{\end{aligned}\end{equation*}}
\newcommand{\eea}{\end{aligned}\end{equation}}
\newcommand{\ee}{\end{equation}}

\renewcommand{\div}{{\rm div }}

\begin{document}
\begin{CJK*}{GBK}{song}
\begin{frontmatter}
\title{Nonlinear Stability of the Rayleigh--Taylor Problem in\\ Quantum Navier--Stokes Equations}
\author[FJ]{Fei Jiang}
\ead{jiangfei0591@163.com}
\author[mmL1]{Yajie Zhang\corref{cor1}}
\ead{zhangyajie315@qq.com}
\author[SJ]{Zhipeng Zhang}\ead{zhangzp@ouc.edu.cn}
\cortext[cor1]{Corresponding author. }
\author[FJ]{Youyi Zhao}\ead{pde\_fzu@163.com}
\address[FJ]{School of Mathematics and Statistics, Fuzhou University, Fuzhou, 350108, China.}
\address[mmL1]{Institute of Applied Physics and Computational Mathematics,
Beijing, 100088, China.}
\address[SJ]{School of Mathematical Sciences, Ocean University of China, Qingdao, 266100, China.}
\begin{abstract}	
It is well-known that the Rayleigh--Taylor (abbr. RT) instability can be completely inhibited by the quantum effect stabilization in proper circumstances leading to a cutoff wavelength in the \emph{linear} motion equations. Motivated by the linear theory, we further investigate the {stability} for the \emph{nonlinear} RT problem of quantum Navier--Stokes equations in a slab with Navier boundary condition, and rigorously prove the inhibition of RT instability by the quantum effect under a proper setting. More precisely, if the RT density profile $\bar\rho$ satisfies an additional stabilizing condition, then there is a threshold ${\varepsilon_{\mm{c}}}$ of the scaled Planck constant, such that if the scaled Planck constant is bigger than ${\varepsilon_{\mm{c}}}$, the small perturbation solutions around an RT equilibrium state are algebraically stable in time. The mathematical proof is realized by a complicated multi-layer energy method with  anisotropic norms of spacial derivatives.
\end{abstract}
\end{frontmatter}
\newtheorem{thm}{Theorem}[section]
\newtheorem{lem}{Lemma}[section]
\newtheorem{pro}{Proposition}[section]
\newtheorem{concl}{Conclusion}[section]
\newtheorem{cor}{Corollary}[section]
\newproof{pf}{Proof}
\newdefinition{rem}{Remark}[section]
\newtheorem{definition}{Definition}[section]
	
\section{Introduction}\label{introud}
	
The equilibrium of a heavier fluid on the top of a lighter one, subject to gravity, is
unstable. In fact, small disturbances acting on the equilibrium will grow and lead to the release of
potential energy, as the heavier fluid moves down under gravity, and the lighter one is displaced
upwards. This phenomenon was first studied by Rayleigh \cite{RLISa} and then Taylor \cite{TGTP}, and is called therefore the Rayleigh--Taylor (abbr. RT) instability. In the last decades, the RT instability has been extensively investigated from physical, numerical, and mathematical aspects, see \cite{CSHHSCPO,WJH,desjardins2006nonlinear,hateau2005numerical,GBKJSAN}, and the references cited therein. It has been also widely analyzed on how the physical factors, such
as elasticity \cite{JFJWGCOSdd}, rotation \cite{CSHHSCPO,BKASMMHRJA}, surface tension \cite{GYTI2,WYJTIKCT}, magnetic fields \cite{JFJSWWWOA,JFJSJMFMOSERT,WYJ2019ARMA,JFJSARMA2019}, capillarity force \cite{LFCZZP}, quantum  {effect} \cite{bychkov2008rayleigh,dolai2017rayleigh} and so on,  {influence} the dynamics of RT instability. In this paper, we are interested in the phenomenon of quantum effect inhibiting RT instability.

Bychkov--Marklund--Modestov might have been the first to find that quantum effect always stabilizes the RT instability by the \emph{linear} motion equations of quantum fluids \cite{bychkov2008rayleigh}. In particular, they obtained the following maximal RT instability growth rate $\sigma$ corresponding to a special perturbation mode:
\begin{align}
\sigma=\sqrt{g\gamma-(\varepsilon\gamma \kappa)^2},
\label{2025072421531}
\end{align}
where $g>0$, $\gamma:=\bar{\rho}'/\bar{\rho}>0$, $\varepsilon>0$ and $\kappa$ are the gravity constant, the (local) inverse length scale  of the density profile $\bar{\rho}$, the scaled Planck
constant and the perturbation wave number. As can be seen from the above relation, the smaller the growth rate, the bigger the scaled Planck
constant. Since Bychkov et.al.'s pioneering work, it has been further widely analyzed on how the quantum {effect} {influence} the RT instability in inhomogeneous rotating plasmas \cite{hoshoudy2009quantum},  quantum
magnetized plasmas \cite{wang2012stabilization}, viscoelastic plasmas models through porous mediums \cite{hoshoudy2011quantum} and so on.

Recently, Chang--Duan investigated the RT instability
in the \emph{nonlinear} motion equations of viscous quantum fluids (i.e. quantum Navier--Stokes equations) defined on the whole space $\mathbb{R}^3$, and derived a threshold $\varepsilon_{\mm{c},\mathbb{R}^3}$ defined by a general density {profile} and $g$ for the scaled Planck constant \cite{MR4625960}. They found a class of density profiles, which always  {result in} the RT instability for any $\varepsilon$ due to $\varepsilon_{\mm{c},\mathbb{R}^3}=\infty$. However they {could not answer} whether $\varepsilon_{\mm{c},\mathbb{R}^3}$ is finite for some density profiles, which {lead to} the RT instability for too small $\varepsilon$. This is also one of the reasons why they did not further investigate whether the quantum effect can inhibit the RT instability.

Recently, the inhibition of RT instability by magnetic fields has been mathematically verified in the non-resistive magnetohydrodynamic (MHD) equations in slabs, see \cite{WYJ2019ARMA,JFJSZYYO} for examples. It is well-known that the RT instability can be completely quenched by the quantum effect stabilization, as well as the magnetic field stabilization, in proper circumstances leading to a cutoff wavelength \cite{wang2012stabilization}.
This motivates us to also consider the RT instability of quantum Navier--Stokes (QNS) equations in a slab. Under such case, we obtain a finite threshold $\varepsilon_{\mm{c}}$ due to the finite height of the slab. By developing a new frame of energy method, which is different {from} the ones for the MHD equations in \cite{WYJ2019ARMA,JFJSZYYO}, we realize the mathematical proof for the inhibition of RT instability by the quantum effect under the case  $\varepsilon>\varepsilon_{\mm{c}}$. Before presenting our detailed result, we shall introduce the so-called quantum RT problem.

\subsection{A quantum RT problem}\label{subsec:01}
\numberwithin{equation}{section}
To conveniently investigate the inhibition of RT instability by the quantum effect, we consider the following system of inhomogeneous incompressible QNS
equations in the presence of a gravitational field as in \cite{MR4625960}:
\begin{align}\label{a0101}
\begin{cases}
\rho_t+ v\cdot \nabla  \rho = 0, \\
\rho v_t+\rho v\cdot \nabla v - \mu \Delta v  + \nabla  {P}  =  2\varepsilon^2\rho\nabla\left(\rho^{-1/2}\Delta\sqrt{\rho}\right)-g \rho \mathbf{e}^3,\\
\mm{div}v=0,
\end{cases}
\end{align}
where the unknowns ${\rho (x,t)} \in \mathbb{R}^{+}$, $v(x,t) \in {\mathbb{R}^3}$ and $P(x,t)$ denote the  particle density, particle velocity and kinetic pressure of the fluid {resp.} at the spacial position $x \in {\mathbb{R}^3}$ for time $t\in\mathbb{R}^+_0:=[0, {+\infty} )$. The letters $\varepsilon>0$,  $\mu>0$ and $g$ {represent} the scaled Planck constant, the viscosity constant and the gravitational constant  {{resp.}}. The expression $- g\rho \mathbf{e}^3$ represents the gravity, where $\mathbf{e}^3 $ is the unit vector with the third component being $1$. In addition, the expression $\rho^{-1/2}\Delta\sqrt{\rho}$ can be interpreted as a quantum potential, i.e. the so-called Bohm potential. In particular, the quantum effect term $2\varepsilon^2\rho\nabla\left(\rho^{-1/2}\Delta\sqrt{\rho}\right)$ can be rewritten as follows
\begin{align} 2\varepsilon^2\rho\nabla\left(\rho^{-1/2}\Delta\sqrt{\rho}\right)=\varepsilon^2\mm{div}\left(
\Delta\rho\mathbb{I}- {\rho }^{-1}\nabla \rho  \otimes \nabla \rho \right) .\label{1pa}\end{align}
We mention that the well-posedness problem of QNS equations has been widely investigated, see \cite{MR2964097,MR1458149,MR2175779,MR2776529,MR2334851,MR2823997,MR3961293,MR2481754,MR2644915} and the references cited therein. In this paper, we consider that the fluid {domain} $\Omega$ is a slab, i.e.
$$\Omega:=\mathbb{R}^2\times (0,h).$$

To introduce the quantum RT problem, we should choose an RT equilibrium state $(\bar{\rho},0)$ to \eqref{a0101}, where the density profile $\bar{\rho}$ only depends on the third spacial variable $x_3$ and satisfies
\begin{gather}
\label{0102}
\inf_{ x_3\in (0,h)}\{\bar{\rho}(x_3)\}>0,\\[1mm]
\label{0102n}\bar{\rho}'|_{x_3=s}>0 \ \mbox{ for some} \ s\in (0,h).
\end{gather}
Here and in what follows, $\bar{\rho}':=\mm{d}\bar{\rho}/\mm{d}x_3$. Then the corresponding pressure profile $\bar{P}$ can be determined by the hydrostatic relation
\begin{equation}
\bar{P}' = {\varepsilon ^2}\left( {\bar \rho '' - {{\left( {{{{{\left| {\bar \rho '} \right|}^2}}}/{{\bar \rho }}} \right)}}} \right)^\prime  - g\bar \rho . \label{equcomre}
\end{equation}	
We remark that the condition  \eqref{0102}
prevents us from treating vacuum, while the condition in \eqref{0102n} is called the RT condition,
which assures that there is at least a region in which the density is larger with increasing height $x_3$, thus leading to the classical RT instability, see \cite[Theorem 1.2]{JFJSO2014}. However we will see that such instability can inhibited by the quantum effect.
	
Denoting the perturbations around the RT equilibrium  state by
$$\varrho=\rho -\bar{\rho}\mbox{ and } v= v- {0}, $$
and then recalling the relations \eqref{1pa} and \eqref{equcomre},
we obtain the perturbation system from \eqref{a0101}:
\begin{equation}\label{1a}
\begin{cases}
\varrho _t+\bar{\rho}'v_3+v\cdot\nabla \varrho=0,\\
(\bar{\rho}+{\varrho}) v_{t}+   (\bar{\rho}+{\varrho}) v\cdot  \nabla v+\nabla  \beta =\mu\Delta v-g{\varrho}\mathbf{e}^3- {\varepsilon ^2}\mathbf{Q}, \\
\mm{div} v=0,
\end{cases}
\end{equation}
where $\beta:=P - {\varepsilon ^2}\Delta \varrho -\bar P$ and
$$\mathbf{Q}:=\div \left( {\frac{\nabla \left( {\varrho  + \bar \rho } \right) \otimes \nabla \left( {\varrho  + \bar \rho } \right) }{{\varrho  + \bar \rho }}- \frac{\nabla \bar \rho  \otimes \nabla \bar \rho }{{\bar \rho }}} \right).$$

Now we further impose the initial-boundary value conditions
\begin{align}\label{1ab}
&(\varrho,v)|_{t=0}=(\varrho^0,v^0),\\
& \nonumber
v|_{\partial\Omega}\cdot{\mathbf{n}}=0\mbox{ and } ((\mathbb{D}v|_{\partial\Omega}){\mathbf{n}})_{\mm{tan}}=0,
\end{align}
where $\partial\Omega$ is the boundary of  $\Omega$,  $\mathbf{n}$ denotes the outward unit normal vector to
$\partial\Omega$, and the subscript ``$\mm{tan}$" means the tangential component of a vector, i.e.   $v_\mm{tan}=v-(v\cdot {\mathbf{n}}){\mathbf{n}}$.
 {The differential operator $\mathbb{D}$ is defined by $\mathbb{D}v = \nabla v+\nabla v^{\top}$, where the superscript $\top$ denotes the transposition.} Since $\Omega$ is a slab, the Navier boundary condition is equivalent to the boundary condition
\begin{align}
(v_3,\partial_3 v_1,\partial_3 v_2)|_{\partial\Omega}=0,
\label{n1}
\end{align}
where we have defined $\partial_3:=\partial_{x_3}$.
We mention that the above boundary condition contributes to the mathematical verification of the inhibition phenomenon in this paper.
For (the sake of) simplicity, we call the initial-boundary value problem \eqref{1a}--\eqref{n1} {\it the  quantum RT (abbr.QRT) problem}. The phenomenon of the quantum effect inhibiting the RT instability mathematically reduces to the stability in time of small perturbation solutions to the above QRT problem with non-trivial initial data.

To analyze the inhibiting effect of the quantum effect term ${\varepsilon ^2}\mathbf{Q}$, we should rewrite it as a sum of both linear and nonlinear terms, i.e.
\begin{align}
\mathbf{Q}=&\div \left( {\frac{\nabla \left( {\varrho  + \bar \rho } \right) \otimes \nabla \left( {\varrho  + \bar \rho } \right) }{{\varrho  + \bar \rho }}- \frac{\nabla \bar \rho  \otimes \nabla \bar \rho }{{\bar \rho }} } \right)\nonumber\\ =&\div \left( \frac{\nabla \left( {\varrho  + \bar \rho } \right)}{{\bar \rho }} \otimes \nabla \varrho  + \frac{\nabla \varrho }{{\bar \rho }}  \otimes \nabla \bar \rho -\frac{{  \varrho }}{{\bar \rho \left( {\varrho  + \bar \rho } \right)}}\nabla \left( {\varrho  + \bar \rho } \right) \otimes \nabla \left( {\varrho  + \bar \rho } \right)  \right)
\nonumber\\ =&\div\left(\frac{{\nabla \bar \rho }}{{\bar \rho }} \otimes \nabla \varrho + \frac{\nabla \varrho }{{\bar \rho }}  \otimes {\nabla \bar \rho }  -\frac{\varrho }{{{{\bar \rho }^2}}}\nabla \bar \rho  \otimes \nabla \bar \rho   \right)\nonumber\\& +\div\bigg( \frac{{\nabla  \varrho }}{{\bar \rho }} \otimes \nabla \varrho + \frac{\varrho^2 }{{{{\bar \rho }^2(\varrho+\bar \rho)}}}\nabla \bar \rho  \otimes \nabla \bar \rho -\frac{{  \varrho }}{{\bar \rho \left( {\varrho  + \bar \rho } \right)}}\nabla \varrho  \otimes \nabla \varrho \nonumber\\&   -\frac{  \varrho }{\bar \rho \left( {\varrho  + \bar \rho } \right)}\left( \nabla \bar \rho \otimes \nabla \varrho+\nabla \varrho \otimes\nabla \bar \rho\right) \bigg)=\mathbf{Q}^{\mm L}+{\mathbf{Q}^{\mm N}},\label{07291549}
\end{align}
{where} we have defined that
\begin{align}
&\mathbf{Q}^{\mm L}:={\partial _3}\left( {  {\frac{{\bar \rho '}}{{\bar \rho }}}  \nabla \varrho} \right)+\left(\left(  \frac{\bar \rho '}{\bar \rho } \right) '{\partial_3}\varrho - \partial _3\left( \left|\frac{\bar \rho '}{\bar \rho }\right|^2\varrho  \right) + \frac{\bar \rho '}{\bar \rho }\Delta \varrho \right)\mathbf{e}^3,\nonumber \\
&\mathbf{Q}^{\mm N}:=\partial_{3}\left(\frac{ \bar \rho' \varrho }{\bar \rho  ( \varrho  + \bar \rho) }     \left(\frac{  \bar{\rho}'\varrho \mathbf{e}^3}{\bar{\rho} }  - {  \nabla \varrho} \right)  \right)
+ \text{div}\left( \frac{ \nabla \varrho \otimes \nabla \varrho }{ \varrho + \bar{\rho} } -\frac{\bar \rho'  \varrho\nabla \varrho \otimes\mathbf{e}^3}{\bar \rho \left( {\varrho  + \bar \rho } \right)} \right)
\label{2508171749}.\end{align}

Inserting the equality \eqref{07291549} into \eqref{1a} yields
\begin{equation}\label{1abx}
\begin{cases}
\varrho _t+\bar{\rho}'v_3=-v\cdot\nabla \varrho,\\
\bar{\rho}  v_{t}+\nabla  \beta-\mu\Delta v+g{\varrho}\mathbf{e}^3+ {\varepsilon ^2}\mathbf{Q}^{\mm L} = - {\varrho}  v_{t}-   (\bar{\rho}+{\varrho}) v\cdot  \nabla v-{\varepsilon ^2}\mathbf{Q}^{\mm N}, \\
\mm{div} v=0.
\end{cases}
\end{equation}

By the standard linear analysis in \cite{MR4625960,LFCZZP},  it is easy to derive a threshold ${\varepsilon_{\mm{c}}}$ and verify that the linearized problem of \eqref{1abx}  is stable under the sharp stability condition
\begin{equation}
\varepsilon>{\varepsilon_{\mm{c}}}: = \mathop {\sup}\limits_{w \in {H_\sigma^* }} \sqrt{\frac{{g\int {\bar \rho 'w_3^2\mm{d}x} }}{\int \left| {{\underline{\bar{\rho}}}}{\nabla {w_3}} \right|^2\mm{d}x}},\label{2saf01504}
\end{equation}
where $\bar{\rho}$ must satisfies the stabilizing condition
\begin{align}
\label{2022205071434}
\inf_{ {x_3}\in (0,h)} \{|\bar{\rho}' {(x_3)}|\}>0
\end{align} and we have used the notations
\begin{align}
&\int:=   \int_{ \mathbb{R}^2\times (0,h)},\ {\underline{\bar{\rho}}}:=\frac{\bar{\rho}'}{\sqrt{\bar{\rho}}} ,
\  H_{\sigma}^*:=\{w\in H_{\sigma}~|~{w_3}\neq 0\}, \label{2508162239} \\
&\mbox{and }H_{\sigma}:=\{w\in  {W^{1,2}(\Omega)}~|~\mm{div}w=0,\ w_3|_{\partial\Omega}=0\}.
\label{20224013022226}
\end{align}
We mention that if $w \in H_{\sigma}$ additionally satisfies $\int \left| {\underline{\bar{\rho}}} {\nabla {w_3}} \right|^2\mm{d}x= 0$,  then ${w_3}=0$ due to \eqref{0102}, \eqref{2022205071434} and \eqref{2202402040948}.

The aim of this paper is to show that \eqref{2saf01504} is also the sharp stability condition of  the QRT problem \eqref{1a}--\eqref{n1} for the small perturbation solutions.

\subsection{Notations}\label{202402081729}
Before stating our main result,  we shall introduce simplified notations which will be used throughout this paper.

(1) Basic notations: $\langle t\rangle:=1+t$. $A:B:=a_{ij}b_{ij}$, where $A:=(a_{ij})_{n\times n}$ and $B:=(b_{ij})_{n\times n}$ are $n\times n$ matrices, and we have used the Einstein convention of summation over repeated indices. $a\lesssim b$ means that $a\leqslant cb$, where the generic positive constant $c$ may vary from line to line, and may depend on the domain $\Omega$, the parameter $\theta$ in Theorem \ref{thm2} and the other known physical functions/parameters, such as $\bar{\rho}$, $\mu$, $g$,  $\varepsilon$ in the QRT problem. $I_a:=(0,a)$ denotes a time interval, in particular, $I_\infty=\mathbb{R}^+:=(0,\infty)$. $\overline{S}$ denotes the closure of a set $S\subset \mathbb{R}^n$ with $n\geqslant 1$, in particular, $\overline{I_T} =[0,T]$ and $\overline{I_\infty} = \mathbb{R}^+_0$.

$\partial_i:=\partial_{x_i}$, where $i=1$, $2$, $3$. $\nabla^0\phi:=\phi$. $\nabla_{\mm{h}}:=(\partial_{1},\partial_{2})^{\top}$,  $\nabla_{\mm{h}}^\bot:=(-\partial_{2},\partial_{1})^{\top}$,  $\Delta_\mm{h}:=\partial^2_{1}+\partial^2_{2}$.  Let $f:=(f_1,f_2,f_3)^{\top}$ and $\psi:=(\psi_1,\psi_2,\psi_3)^{\top}$ be vector functions defined in a 3D domain, we then define that $f_{\mm{h}}:=(f_1,f_2)^{\top}$,, $ \mm{div}_{\mm{h}} f_{\mm{h}}:=\partial_{1} f_1 +\partial_{2} f_2$,
$\mm{curl}{f}:=(\partial_{2}f_3-\partial_{3}f_2,
\partial_{3}f_1-\partial_{1}f_3,\partial_{1}f_2-\partial_{2}f_1)^{\top}$ and $f\land \psi:=(f_2\psi_3 - f_3\psi_2,f_3\psi_1- f_1\psi_3, f_1\psi_2-f_2\psi_1)^\top$.

 For a given function $f\in L^2(\Omega)$, the horizontal Fourier transform of $f$ can be given as follows
\begin{equation}\label{hftxh}
\hat{f}(\xi_{\mm{h}},x_3) =\frac{1}{2\pi} \int_{\mathbb{R}^2} f(x_{\mm{h}},x_3) e^{-\mm{i}x_{\mm{h}}\cdot \xi_{\mm{h}}}\mm{d}x_{\mm{h}},
\end{equation}
where $\xi_{\mm{h}}:=(\xi_1,\xi_2)^\top$. Thus we further define the operator $\Lambda^{s}_{\mm{h}}$ with $s\in\mathbb{R}$ by
$$\Lambda^{s}_{\mm{h}}f(x_\mm{h},x_3)=\int_{\mathbb{R}^2}|\xi_\mm{h}|^{s}\hat{f}
(\xi_{\mm{h}},x_3)e^{ \mm{i}x_{\mm{h}}\cdot \xi_{\mm{h}}}\mm{d}\xi_{\mm{h}}.$$

(2) Simplified Banach spaces, norms and semi-norms:
\begin{align}
&L^p:=L^p (\Omega),\ {H}^i:=W^{i,2}(\Omega ),\  H^1_0:=\{\phi\in H^1~|~\phi|_{\partial\Omega}=0\},  \nonumber \\
&    {H}^j_{\diamond}:=\{w\in H^j~|~ \partial_3w_1 =\partial_3w_2=w_3 =0\mbox{ on }\partial\Omega \},\  {H}^{\sigma,j}_{\diamond}:=H_\sigma\cap  {H}^j_{\diamond},\nonumber\\
&{ {H}}^6_{\bar{\rho}}:=\{\phi\in H^6~|~\phi|_{\partial\Omega}=0,\ \partial^l_{3}(\phi+\bar{\rho})|_{\partial\Omega}=0\mbox{ for }1\leqslant l\leqslant 4\},\nonumber\\
&{H}^{k}_{-s,i}: =\{w\in H^k~|~\Lambda^{-s}_{\mm{h}}w\in H^{i} \},\  {H}^{\diamond,\sigma, 5}_{-s,i}: =  {H}^{\sigma,5}_\diamond\cap {H}^{5}_{-s,i},\ {H}^{\bar{\rho},6}_{-s,i}:= { {H}}^6_{\bar{\rho}}\cap {H}^{6}_{-s,i},
\nonumber  \\
& \|\cdot \|_i :=\|\cdot \|_{H^i},\   \|\cdot\|_{{-s},k}:= \|\Lambda^{-s}_{\mm{h}}\cdot\|_k, ~\|\cdot\|_{{i-s},k}:= \sum_{\alpha_1+\alpha_2=i}\|\Lambda^{-s}_{\mm{h}}\partial_{1}^{\alpha_1}\partial_{2}^{\alpha_2}\cdot\|_k,\nonumber
\end{align}
where $s\in\mathbb{R}$, $\bar{\rho}$ is the density profile provided in Subsection \ref{subsec:01}, $1\leqslant p\leqslant \infty$ is a real number,  $i$, $k\geqslant 0$, $j\geqslant 2$ are integers. It should be noted that the space ${H}^{k}_{-s,i}$ is endowed with the norm {$\|\cdot\|_{k}+\|\cdot\|_{-s,i}$}.
For simplicity, we denote $\sqrt{\sum_{i=1}^n\|w_i\|_X^2}$ by $\|(w_1,\cdots,w_n)\|_X$, where $\|\cdot\|_X$ represents a norm or a semi-norm, and $w_i$ are scalar functions or vector functions for $1\leqslant i\leqslant n$.

(3) Simplified spaces of functions with values in a Banach space:
\begin{align}
& {{\mathfrak{P}} _{T}:=\{\varrho\in
C^0(\overline{I_T} ,   {H}^{\bar{\rho},6}_{-s,3} )~|~ \varrho_t\in C^0(\overline{I_T} , {H}^5_{-s,2})\cap  L^2(I_T,{H}^5_{-s,2})\}}, \nonumber\\
& {\mathfrak{W}}_{ T}: =   \{v\in C^0(\overline{I_T},{H}^{\diamond,\sigma, 5}_{-s,2} )\cap L^2(I_T,{H}^6_{-s,3}) ~|~
v_t\in C^0(\overline{I_T},{H}^3_{-s,0})\cap  L^2(I_T,{H}^4_{-s,1})\}, \nonumber
\end{align}
where $0\leqslant T\leqslant \infty$.
	
(4) Functionals of linearized potential energy: for $r\in H^1_0$,
\begin{align}&E_{\mm{L}}(r):=\varepsilon ^2 \| {\underline{\bar{\rho}}}\nabla r\|^2_0- g\|\sqrt{\bar{\rho}'} r\|^2_0\label{edeE},\\
\label{eeE}
&E(r):=\varepsilon ^2\|\nabla r/\sqrt {\bar \rho }\|^2_{0}+\int\frac{\varepsilon ^2 (\bar{\rho}''/\bar \rho)'-g}{\bar{\rho}'}r^2\mm{d}x ,
\end{align}
where the definition of  $\underline{\bar{\rho}}$ can be found in \eqref{2508162239}.
	
(5) Functionals of tangential (generalized) energy/dissipation:
\begin{align}
&	{\underline{\mathcal{E}}}:= \|\varrho\|_{{1},1}^2+  \|v\|_{{1},0}^2,\ \underline{\mathcal{D}} := \|\nabla v\|_{1,0}^2,\label{2022402200520957}\\
&\underline{\mathfrak{E}}:= \|\varrho\|^2_{-s,3}+\|(\varrho_{t},v)\|^2_{-s,2}+\|v_t\|_{-s,0}^2, \nonumber\\& \underline{\mathfrak{D}}:= \|\varrho\|^2_{1-s,2}+\|(\varrho_t,\nabla v)\|^2_{-s,2}+\| v_t\|^2_{-s,1}.
\label{202240220052095s7}
\end{align}
	
(6) Functionals of energy/dissipation:
\begin{align}
&\mathfrak{E}:= \|\varrho\|^2_3+\|(\varrho_{t},v)\|^2_2+\|v_t\|_0^2,\ \mathfrak{D}:=\|\varrho\|^2_{1,2}		
+\|(\varrho_t,\nabla v)\|^2_2+\| v_t\|^2_1,\label{diD}	\\		
&\mathcal{E}:= \|\varrho\|^2_6+\|(\varrho_{t},v)\|^2_5+\|v_t\|_3^2,\  \mathcal{D}:=\|\varrho\|^2_{1,5}+\|(\varrho_t,\nabla v)\|^2_5+\| v_t\|^2_4.\label{di}\end{align}

\subsection{A stability result}

Now  we state the stability result for the QRT problem, which presents that the quantum effect can inhibit the RT instability for properly large scaled Planck constant.
\begin{thm}\label{thm2}
Let $\mu$, $\varepsilon $ be positive constants, 
\begin{align}
\label{2510031559}
\theta\in(0,({19-\sqrt{349}})/{6}),\ s=1-\theta,
\end{align} and $\bar{\rho}\in {C^8}[0,h]$ satisfy the stabilizing condition \eqref{2022205071434} and the boundary condition
\begin{align}
\label{0102n1}	\bar{\rho}''|_{\partial\Omega}=	\bar{\rho}''''|_{\partial\Omega}=0.
\end{align}
Under the sharp stability condition \eqref{2saf01504},
there is a sufficiently small constant $\delta\in (0,1)$, such that for any $( \varrho^0,v^0)\in  { {H}}_{-s,3}^{\bar{\rho},6}\times  {H}_{-s,2}^{\diamond,5}$  satisfying necessary compatibility conditions and a smallness condition
\begin{align*}
\|(\nabla \varrho^0,v^0)\|_5+ \|(\nabla \varrho^0,v^0)\|_{-s,2}\leqslant\delta,
\end{align*}
the QRT problem  \eqref{1a}--\eqref{n1} admits a unique global classical  solution $(\varrho,v,\nabla\beta)\in {\mathfrak{P}} _{\infty}\times {\mathfrak{W}_{\infty} }\times C^0(\mathbb{R}_0^+,{H}^5)$. Moreover, the solution enjoys the following stability estimate with algebraic decay-in-time	
\begin{align}
& \mathcal{E}_{\mathrm{total}}(t) +\int_0^t\mathcal{D}_{\mathrm{total}}(\tau)    \mm{d}\tau
\lesssim \|(\nabla \varrho^0,v^0)\|_{-s,2}^2+\|(\nabla \varrho^0,v^0)\|_5^2 \label{eh3}
\end{align}	
for any $t>0$, where we have defined that
\begin{align}
&
\mathcal{E}_{\mathrm{total}}(t) := \underline{\mathfrak{E}}(t)+ \mathcal{E}(t) + \langle t \rangle^{a} \mathfrak{E}(t) + \langle t \rangle^{3a/2} \underline{\mathcal{E}}(t), \\
&  \mathcal{D}_{\mathrm{total}}(t):= \underline{\mathfrak{D}}(t) + \mathcal{D}(t) + \langle t \rangle^{a} \mathfrak{D}(t)+ \langle t \rangle^{3a/2} \underline{\mathcal{D}}(t) ,\\
&
a:=\theta+2/3.
 \label{disfa}\end{align}
\end{thm}
\begin{rem}It is well-known that
$$ {\int^{h}_{0} \phi^2 {\mm{d}s}}\leqslant \phi^{-2}h^2 {\|\phi'\|^2_{L^2(0,h)}}\mbox{ for any }
\phi \in H^1_0(0,h).$$
In view of the above fact and the relation of $\varepsilon_{\mm{c}}=a_3$ in \eqref{2022309032042}, it is easy to check that
$$ {\varepsilon_{\mm{c}}} \leqslant h\pi^{-1}\sqrt{g\|\bar \rho '\|_{L^\infty(0,h)}\| |\bar \rho '|^{-2}
\bar\rho\|_{L^\infty(0,h)}},
$$
which means that the stronger the inhibition effect of the quantum effect, the smaller the height of the slab. In particular, if $\bar{\rho}'\equiv \alpha>0$, then the above inequality reduces to
\begin{align} {\varepsilon_{\mm{c}}} \leqslant{h}{\pi}^{-1} \sqrt{{g{\alpha}^{-1}\|
\bar\rho\|_{L^\infty(0,h)}}} . \label{rem1}
\end{align}
It is well-known that if $\alpha$ is bigger, the RT instability is more likely to occur.
However, we can clearly see  from the above estimate \eqref{rem1} with given $g$, $h$ and $\|\bar \rho \|_{L^\infty(0,h)}$ that the inhibition effect from the linearized term $-{\varepsilon ^2}\mathbf{Q}^{\mm L}$ of the quantum effect will become stronger, as the value of $\alpha$ increases.
\end{rem}
\begin{rem} When $\varepsilon< \varepsilon_{\mm{c}}$ or $\varepsilon_{\mm{c}} = \infty$, we can construct RT instability solutions by further exploiting the bootstrap instability method in \cite{jiang2023rayleigh}.
It should be noted that $\varepsilon_{\mm{c}} = +\infty$, if $\bar{\rho}$  satisfies $\bar{\rho}' \geqslant 0$ and, for some given positive constants $\tilde{c}_1$, $\tilde{c}_2$, $\tau$, $\eta$ and an interval $(x_3^0-\eta,x_3^0+\eta) \subset (0,h)$,
\begin{equation}
\tilde{c}_1 \leqslant  {\bar{\rho}'}{|x_3 - x_3^0|^{-({2+\tau})}} \leqslant \tilde{c}_2 \text{ for any } x_3 \in (x_3^0-\eta,x_3^0+\eta),
\end{equation}
please refer to \cite[Proposition 2.1]{zhang2022rayleigh}. This is the reason why we shall need the stabilizing condition \eqref{2022205071434} to exclude the case $\varepsilon_{\mm{c}} = \infty$.
\end{rem}
\begin{rem}	
We mention that the boundary conditions in \eqref{0102n1} will be used to estimate  the {higher-order} spacial derivatives of $(\varrho,v)$ in Lemma \ref{lem1a}, which finally contribute to establishing the dissipative estimate of $\|\varrho\|_{1,5}$ in  \eqref{r4}.
\end{rem}

\subsection{Sketch the proof for the stability result}
Recently, the inhibition of RT instability by the magnetic tension in non-resistive MHD equations was mathematically verified under the Lagrangian coordinates by exploiting the dissipative estimates of
$\eta$ under small perturbation, where $\eta$ denotes the
departure function of fluid particles \cite{JFJSJMFMOSERT,WYJ2019ARMA}. However, such method seems to be extremely difficulty to further investigate our problem, since we can not capture the dissipative
estimates of $(\eta_1, \eta_2)$. In this paper, we will develop another proof frame for the inhibition of RT instability by the quantum effect stated in Theorem \ref{thm2}. The key step is to derive the \emph{a priori} stability energy estimate \eqref{eh3}.	Next  we briefly sketch the proof of Theorem \ref{thm2}, the details of which will be presented in Sections \ref{sec:global} and \ref{subsec:08a}.

Let $(\varrho, v)$ be a solution {of} the QRT problem \eqref{1a}--\eqref{n1} defined on $\Omega\times (0,T)$ and satisfy the
\emph{a priori} assumption
\begin{align}
{B(T)} \leqslant \delta^* \in (0,1) ,\label{imvarrho}\end{align}
where we have defined that
\begin{align}
B(T):=\sup_{0\leqslant t\leqslant T}\{ \mathcal{E}_{\mm{\mm{total}}}(t)\}+{\int_0^T\mathcal{D}_{\mathrm{total}}(t)    \mm{d}t}.
\label{202508172246}
\end{align}

  It is easy to infer that $(\varrho, v)$ enjoys the following energy identity {in} a differential form
\begin{align}
{\frac{1}{2}}\frac{\mathrm{d}}{\mathrm{d}t} \left( E(\varrho) + \| \sqrt{\rho} v \|^2_{0} \right) +{\mu \|\nabla v \|^2_{0} }=\tilde{\mathcal{I}}, \label{lemiqasfd11xz}
\end{align}
where $\tilde{\mathcal{I}}$ is an integral over $\Omega$ involving nonlinear terms, see \eqref{2508111549} for {its} expression. In particular, integrating the above identity over $(0,T)$ yields
\begin{align}
 E(\varrho)+\|\sqrt{\rho}v\|^2_{0}  + {2\mu\int^T_0\|\nabla v \|^2_{0}\mathrm{d}t} = E(\varrho^0)+\| \sqrt{\rho^0} v^0 \|^2_{0}+2\int_0^T\tilde{\mathcal{I}}\mm{d}t.
 \label{2508311429}
\end{align}

 Under the sharp stability condition \eqref{2saf01504}, one has (referring to Lemma \ref{lem:20201224215})
  $$\|\varrho \|_1^2\lesssim E(\varrho).$$
Thanks to the above stabilizing estimate and the lower/upper bounds of density (see \eqref{im1}), we further deduce from \eqref{2508311429} that
  \begin{align}
\|\varrho \|_1^2+\| v\|^2_{0}  + \mu \int_0^T\|\nabla v \|^2_{0}\mm{d}t \lesssim  \|\varrho^0\|_1^2+\|v^0\|_0^2+\int_0^T\tilde{\mathcal{I}}\mm{d}t.
\label{2508131322}
\end{align}
In addition, we can capture the following dissipative estimate of the tangential derivatives of $\varrho$ from the vortex equations of \eqref{1a}$_2$ and the sharp stability condition (referring to \eqref{Pfsda2} with $l=0$): \begin{align}
&c^{-1}\|\varrho\|_{1,1}^2\leqslant E( \partial_2\varrho)+E( \partial_1\varrho)
={\mathcal{V}}^{\mm{L}}+ {\mathcal{V}}^{\mm{N}}\lesssim \|\nabla v\|^2_{1}+\|v_t\|^2_0 + {\mathcal{V}}^{\mm{N}} ,\label{LEM2P1asd}
\end{align}
where the definitions of ${\mathcal{V}}^{\mm{L}}$ resp. ${\mathcal{V}}^{\mm{N}}$ can be found in \eqref{2025081313071} resp. \eqref{202508131307}, and  ${\mathcal{V}}^{\mm{N}}$ is an integral over $\Omega$ involving nonlinear terms.

Obviously, we can not expect to close energy estimates from the above basic estimates {\eqref{2508131322} and \eqref{LEM2P1asd}}. Therefore we shall further extend the above basic estimates to the case of {higher-order} derivatives for closing the integrals involving nonlinear terms. By performing careful energy estimates, we will establish
the following high-order energy inequality for small perturbation solutions:
\begin{align}
\label{omessetsimQa}
&\sup_{0 \leqslant t \leqslant T}\{\mathcal{E}(t)\} + \int_0^T \mathcal{D}(t) \mathrm{d}t \lesssim \mathcal{E}(0) + \sup_{0 \leqslant t \leqslant T}\{\|\varrho\|_1\|\varrho\|_6^2\} + \int_0^T \Big( \sqrt{\mathcal{E}} \mathcal{D} +\mathfrak{N}(t) \Big) \mathrm{d}t,
\end{align}
where we have defined that
\begin{align}
\mathfrak{N}(t):=\|\nabla v\|_{1,0}\mathcal{E} +\sqrt{\mathcal{E}}
\left( \|\varrho\|_{3}\left(\|\nabla v\|_2 \sqrt{\mathcal{E}}+ \|\varrho\|_{3}\sqrt{\mathcal{D}}\right)+(\|\nabla v\|_{1,0}\|\nabla v\|_{1})^{\frac{1}{2}}(\mathcal{E} \mathcal{D})^{\frac{1}{4}}\right).
\label{msaft}
\end{align}	
 
To control the time integral of $\mathfrak{N}$, we require the decay estimates of $\|\nabla v\|_{1,0}$, $\|\varrho\|_{3}$ and $\|\nabla v\|_2$.
To this purpose, we deduce the tangential energy inequality with decay-in-time
\begin{align}
&\sup_{0 \leqslant t \leqslant T}\left\lbrace \langle t \rangle^{3a/2} \underline{\mathcal{E}}(t)\right\rbrace  + \int_0^T \langle t \rangle^{3a/2} \underline{\mathcal{D}}(t) \mathrm{d}t\lesssim \underline{\mathcal{E}}(0) +B^{3/2}(T) +\int_0^T \langle t \rangle^{a/2} \mathfrak{D}(t) \mathrm{d}t \label{qecc1aa}
\end{align}
	and  the low-order energy inequality  with decay-in-time
\begin{align}
\label{qecc1a}
&\sup_{0 \leqslant t \leqslant T}\left\lbrace  \langle t \rangle^{a} \mathfrak{E}(t) \right\rbrace  + \int_0^T \langle t \rangle^{a} \mathfrak{D}(t) \mathrm{d}t \lesssim \mathfrak{E}(0) +B^{3/2}(T)+ \int_0^T \langle t \rangle^{a-1} \left( \|\varrho\|^2_{3} + \|v\|^2_{2} \right) \mathrm{d}t.
\end{align}

To control $\int_0^T \langle t \rangle^{a-1} \left( \|\varrho\|^2_{3} + \|v\|^2_{2} \right) \mathrm{d}t$ in
\eqref{qecc1a}, we will use the interpolation inequality in \eqref{v1} and Young's inequality to infer that, for any given constant $\chi$ (see \eqref{2508311447} for the detailed derivation),
\begin{align}\label{lem2dx2aac}
&\int_0^T\langle t \rangle^{a-1}(\|\varrho\|^2_{3} + \|v\|^2_{2}) \mathrm{d}t
\nonumber\\&\lesssim  \chi^{1+s}\left(  \sup_{0 \leqslant t \leqslant T} \left\lbrace \underline{\mathfrak{E}}(t)\right\rbrace\int_0^T \langle t \rangle^{(a-1)(1+s)-{2as}/{3}}\mathrm{d}t \right)  + \int_0^T \chi^{-\frac{1+s}{s}} \left( \langle t \rangle^{a} \mathfrak{D}(t) + \|\varrho\|^{2}_{1,5} \right)\mathrm{d}t.
\end{align}
Thus we shall estimate for the negative tangential derivatives of $\varrho$ and $v$ in $\underline{\mathfrak{E}}(t)$. Fortunately, we can build the following energy inequality of negative tangential derivatives
\begin{align}
\label{lem2dx2aa}
&\sup_{0 \leqslant t \leqslant T} \left\lbrace  \underline{\mathfrak{E}}(t)\right\rbrace  + \int_0^T  \underline{\mathfrak{D}}(t) \mathrm{d}t  \lesssim  \underline{\mathfrak{E}}(0) +  B^{3/2}(T).
\end{align}
In addition, we need the following condition for the parameters $a$ and $s$
	 \begin{align}\label{ab1}
		 (a-1)(1+s)<{2as}/{3} -1 
		\end{align}to obtain the uniform estimate with respect to $T$:
\begin{align}
\label{2510031506}
\int_0^T \langle t \rangle^{(a-1)(1+s)-{2as}/{3}}\mathrm{d}t\lesssim 1.
\end{align}
Thanks to \eqref{lem2dx2aac} and \eqref{2510031506}, one easily deduces from   \eqref{qecc1aa}, \eqref{qecc1a} and \eqref{lem2dx2aa} that
\begin{align}
\label{lem2dx3a}
&\sup_{0 \leqslant t \leqslant T} \left\lbrace  \underline{\mathfrak{E}}(t) + \langle t \rangle^{a} \mathfrak{E}(t) +\langle t \rangle^{3a/2} \underline{\mathcal{E}}(t) \right\rbrace   + \int_0^T \Big( \langle t \rangle^{a} \mathfrak{D}(t) + \langle t \rangle^{3a/2} \underline{\mathcal{D}}(t)  + \underline{\mathfrak{D}}(t) \Big) \mathrm{d}t \nonumber \\
&\lesssim \mathcal{E}_{\mathrm{total}}(0)+ \int_0^T \|\varrho\|^{2}_{1,5} \mathrm{d}t +B^{3/2}(T).
\end{align}

Exploiting the above estimate with time-decay, we can control $\mathfrak{N}$. More precisely, we have
	\begin{align*}
		\int_0^T\|\nabla v\|_{1,0}{\mathcal{E}}\mm{d}t \lesssim B(T)\left( \int_0^T \langle t \rangle^{-3a/2} \mathrm{d}t\right)^{1/2}\left( \int_0^T \langle t \rangle^{3a/2} \underline{\mathcal{D}}(t) \mathrm{d}t\right)^{1/2}\lesssim B^{3/2}(T),
	\end{align*} 
if $a$ additionally satisfies
\begin{align}
&{3a}>{2}. \label{ab21} 
\end{align}
Similarly, the integrations of the rest terms in $\mathfrak{N}$ over $(0,T)$ can be bounded by $B^{3/2}(T)$ (see the detailed derivation in \eqref{2510031549}). Therefore we arrive at
	\begin{align*}
		\int_0^T\mathfrak{N}\mm{d}t \lesssim  B^{3/2}(T).
	\end{align*} 
Finally, putting \eqref{omessetsimQa} and \eqref{lem2dx3a} together, we obtain the desired \emph{a priori} stability estimate \eqref{eh3}   under the \emph{a priori} assumption \eqref{imvarrho} with sufficiently small $\delta^*$ (see Proposition \ref{2508191327} for the details). This fact, together with the unique local(-in-time) solvability of the QRT problem in Proposition \ref{202102182115},  immediately implies the unique global solvability of small perturbation solutions in Theorem \ref{thm2}. 
We mention that the detailed {derivations} of the above key four energy inequalities \eqref{omessetsimQa}, \eqref{qecc1aa}, \eqref{qecc1a} and \eqref{lem2dx2aa} will be presented in Proposition \ref{lem3s}.
 
We end this section by explaining why we pose the condition \eqref{2510031559} for $s$ and $\theta$.
For simplicity, we choose $\theta=a -{2}/{3}\in (0,1)$  and $s=1-\theta\in (0,1)$. Substituting this two relations into 
\eqref{ab1} yields $3\theta^2 -{19}\theta + 1> 0$. Solving this quadratic inequality for $\theta\in (0,1)$, we find that it is satisfied for $\theta \in \left(0, ({19 - \sqrt{349}})/{6} \right)$. In particular, the condition $s\in (0,1)$ makes sure 
that the product estimates in $\dot{H}^{-s}$ in Lemma \ref{2024v3} can be applied to estimate the nonlinear terms.

\section{\emph{A priori} estimates}\label{sec:global}
This section is devoted to establishing \emph{a priori} stability estimate \eqref{eh3} for the small perturbation solutions. To this purpose, let $(\varrho, u,\beta)$ be the sufficiently smooth solution of the QRT problem  \eqref{1a}--\eqref{n1} defined on $\Omega\times I_T$.
Moreover, $(\varrho, u)\in {\mathfrak{P}} _{T}\times {\mathfrak{W}_{T} }$ satisfies \eqref{imvarrho}. In addition, $\varrho^0$ satisfies both the upper and lower bounds of the initial density:
\begin{equation}
0< {2}^{-1}{\inf\limits_{x_3\in(0,h)}\big\{\bar{\rho}(x_3) \big\}} \leqslant \inf\limits_{x\in\Omega}\big\{{\rho}^{0}(x)\big\} \leqslant \sup\limits_{x\in\Omega}\big\{{\rho}^{0}(x)\big\}  \leqslant 2\sup\limits_{x_3\in(0,h)}\big\{\bar{\rho}(x_3)\big\},\label{im1a}
\end{equation}
where $\rho^0:=\rho|_{t=0}$ and $\rho:= \bar{\rho}+\varrho$. It follows from the mass equation \eqref{a0101}$_{1}$ that
\begin{equation}
0<\inf\limits_{x\in\Omega}\big\{{\rho}^{0}(x)\big\}\leqslant {\rho}(x,t)\leqslant \sup\limits_{x\in\Omega}\big\{{\rho}^{0}(x)\big\}\mbox{ for any }(x,t)\in \Omega\times I_T.\label{im1}
\end{equation}

\subsection{Preliminaries}
To begin with, we derive more
boundary conditions for $(\varrho, v , \omega)$   and the elliptic   estimates for $\varrho$, $v$ and $v_t$, where $\omega: = \mm{curl} v$.
\begin{lem}
\label{2018051410720}
For any given $t\in [0,T]$,
 the solution $(\varrho,v)$ satisfies the following boundary conditions
\begin{align}\label{varrho}
&( \partial^{2[i/2]+1}_3\rho, \partial_3^{2i}\varrho,\partial_3^{2i}\tilde{\varrho}) |_{\partial\Omega}=0,\\
&(\partial^{2i+1}_3v_{\mm{h}},\partial^{2i}_3v_{3}) |_{\partial\Omega}=0,\label{vbian} \\
&\label{omega} (\partial^{2i}_{3}\omega_{\mm{h}},\partial_{3}^{2i+1}\omega_3)|_{\partial\Omega}=0,
\end{align}	
where $0\leqslant i\leqslant 2$ and $[i/2]$ means {the largest integer, which is not greater than} $i/2$, and we have defined
\begin{equation}
\label{hftxha} \tilde{\varrho}:=\varrho/\bar{\rho}'.
\end{equation}
\end{lem}
\begin{pf}
Since $\varrho\in {\mathfrak{P}} _{T}$, it holds obviously {that}
\begin{align}
&( \partial^{2[i/2]+1}_3\rho, \partial_3^{2i}\varrho ) |_{\partial\Omega}=0.
\label{2508171135}
\end{align}
Using the boundary conditions of $\bar{\rho}$ in \eqref{0102n1} and $\partial_3^{2i}\varrho|_{\partial\Omega}=0$ in \eqref{2508171135}, we obtain
\begin{equation}
\partial_3^{2i} \left(\frac{\varrho}{\bar{\rho}'}\right)=\sum_{k=0}^{2i}{C_{2i}^k}\partial_3^{k}\left(\frac{1}{\bar{\rho}'}\right)
\partial_3^{2i-k} \varrho=0\  \text{on}\ \partial \Omega\label{nfdsa22aa},
\end{equation}
where {$C_{2i}^k=(2i)!/k!(2i-k)!$}. Putting \eqref{2508171135} and \eqref{nfdsa22aa} together yields \eqref{varrho}.

Thanks to the boundary condition of $v$ in \eqref{n1}, we directly get	
\begin{align}
&\omega_{\mm{h}}=\partial_{3}\omega_{3}=\partial_{3}(\partial_{1}v_2-\partial_{2}v_1) =0 \text{ on }\partial \Omega, \label{omegsfdaa1s}
\end{align}
which implies that \eqref{vbian} and \eqref{omega} hold for $i=0$ by \eqref{n1} and \eqref{omegsfdaa1s}.
Next we derive \eqref{vbian} and \eqref{omega} with $i=1$ and $2$.

Applying curl to the momentum equation \eqref{1a}$_2$,  we obtain the following vortex equation
\begin{align}
\label{l41a}
& {\rho}  \omega_t  =\mu\Delta \omega+ {\mathbf{L}}- {\mathbf{M}}-   {\rho} v\cdot  \nabla \omega - \mathbf{N} ,
\end{align}		
where we have defined that
\begin{align}
\label{curla}
\begin{cases}
{\mathbf{M}}  :=(-\bar{\rho}'\partial_{t}v_2,\bar{\rho}'\partial_{t}v_1,0)^{\top},\\
 {\mathbf{L}}: =
\left(g-
{{\varepsilon ^2}{{\left( { {{\bar \rho ''}}/{{\bar \rho }}} \right)}^\prime } } \right){\left( {- {\partial _2},{\partial _1},0} \right)^{\top}}\varrho \\
\qquad\ +\varepsilon ^2 \left( {\bar \rho '}\left(-\partial _2,\partial _1,0\right)^\top\Delta \varrho /{\bar \rho }-  { {{{{\left| {\bar \rho '}/ {\bar \rho } \right|}^{2}}}}{{\left({ - {\partial _2},{\partial _1},0} \right)}^{\top}}{\partial _3}\varrho  } \right),\\\mathbf{N}:=  \mathbf{N}^{\mm{m}} +\mathbf{N}^{\mm{c}} + {\varepsilon ^2}\mm{curl}\mathbf{Q}^{\mm N},\\
\mathbf{N}^{\mm{m}}:=  (\partial_{2} {\varrho}\partial_{t}v_3-\partial_{3} {\varrho}\partial_{t}v_2,\partial_{3} {\varrho}\partial_{t}v_1-\partial_{1} {\varrho}\partial_{t}v_3,\partial_{1} {\varrho}\partial_{t}v_2-\partial_{2} {\varrho}\partial_{t}v_1)^{\top},\\
\mathbf{N}^{\mm{c}}:=\big(\partial_{2}( {\rho}v)\cdot  \nabla v_3-\partial_{3}( {\rho}v)\cdot  \nabla v_2,\partial_{3}( {\rho}v)\cdot  \nabla v_1-\partial_{1}( {\rho}v)\cdot  \nabla v_3,\\
\qquad\quad  \partial_{1}( {\rho}v)\cdot  \nabla v_2-\partial_{2}( {\rho}v)\cdot  \nabla v_1\big)^{\top}.
\end{cases}	\end{align}	
In particular, we can compute out that
\begin{align}
\mm{curl}\mathbf{Q}^{\mm N}= &\partial_3 \left(\nabla \frac{(\bar{\rho}')^2 \varrho^2}{\bar{\rho}^2 (\varrho + \bar{\rho})} \land  \mathbf{e}^3 - \nabla \frac{\bar{\rho}' \varrho}{\bar{\rho} (\varrho + \bar{\rho})} \land  \nabla \varrho\right) \nonumber \\
&+ \mm{curl}\left( \mm{div} \left( \frac{\nabla \varrho \otimes \nabla \varrho}{\varrho + \bar{\rho}} \right) \right) - \nabla \mm{div} \left( \frac{\bar{\rho}' \varrho}{\bar{\rho} (\varrho + \bar{\rho})}\nabla \varrho \right)\land  \mathbf{e}^3\nonumber \\
=&\Bigg(\partial_3\left(\frac{\bar{\rho}^{\prime\prime}\varrho\partial_2\varrho}{\varrho+\bar{\rho}}\right)-\partial_2\left(\frac{\bar{\rho}^{\prime}\varrho}{\bar{\rho}(\varrho+\bar{\rho})}\Delta\varrho+\nabla\left(\frac{\bar{\rho}^{\prime}\varrho}{\bar{\rho}(\varrho+\bar{\rho})}\right)\cdot\nabla\varrho\right)\nonumber\\&+\sum_{j=1}^3 \left( \partial_2 \left( \frac{ \partial_j^2 \varrho \partial_3 \varrho + \partial_j \varrho \partial_j \partial_3 \varrho }{\varrho + \bar{\rho}} - \frac{ \partial_j \varrho \partial_3 \varrho \partial_j (\varrho + \bar{\rho}) }{(\varrho + \bar{\rho})^2} \right) \right. \nonumber \\
& \left. - \partial_3 \left(\frac{ \partial_j^2 \varrho \partial_2 \varrho + \partial_j \varrho \partial_j \partial_2 \varrho }{\varrho + \bar{\rho}} - \frac{ \partial_j \varrho \partial_2 \varrho \partial_j (\varrho + \bar{\rho}) }{(\varrho + \bar{\rho})^2} \right) \right),\nonumber \\
 &-\partial_3\left(\frac{\bar{\rho}^{\prime\prime}\varrho\partial_1\varrho}{\varrho+\bar{\rho}}\right)+\partial_1\left(\frac{\bar{\rho}^{\prime}\varrho}{\bar{\rho}(\varrho+\bar{\rho})}\Delta\varrho+\nabla\left(\frac{\bar{\rho}^{\prime}\varrho}{\bar{\rho}(\varrho+\bar{\rho})}\right)\cdot\nabla\varrho\right)\nonumber\\&+\sum_{j=1}^3 \left( \partial_3 \left( \frac{ \partial_j^2 \varrho \partial_1 \varrho + \partial_j \varrho \partial_j \partial_1 \varrho }{\varrho + \bar{\rho}} - \frac{ \partial_j \varrho \partial_1 \varrho \partial_j (\varrho + \bar{\rho}) }{(\varrho + \bar{\rho})^2} \right) \right. \nonumber \\
& \left. - \partial_1 \left( \frac{ \partial_j^2 \varrho \partial_3 \varrho + \partial_j \varrho \partial_j \partial_3 \varrho }{\varrho + \bar{\rho}} - \frac{ \partial_j \varrho \partial_3 \varrho \partial_j (\varrho + \bar{\rho}) }{(\varrho + \bar{\rho})^2} \right)\right),\nonumber \\
 & \sum_{j=1}^3 \left( \partial_1 \left(\frac{ \partial_j^2 \varrho \partial_2 \varrho + \partial_j \varrho \partial_j \partial_2 \varrho }{\varrho + \bar{\rho}} - \frac{ \partial_j \varrho \partial_2 \varrho \partial_j (\varrho + \bar{\rho}) }{(\varrho + \bar{\rho})^2} \right) \right. \nonumber \\
& \left. - \partial_2 \left( \frac{ \partial_j^2 \varrho \partial_1 \varrho + \partial_j \varrho \partial_j \partial_1 \varrho }{\varrho + \bar{\rho}} - \frac{ \partial_j \varrho \partial_1 \varrho \partial_j (\varrho + \bar{\rho}) }{(\varrho + \bar{\rho})^2} \right) \right)\Bigg)^{\top}.\label{qn3a}
\end{align}

Exploiting the boundary conditions of $(\varrho,\partial_3\rho, v_3,\partial_3 v_{\mm{h}},\omega_{\mm{h}})$ in \eqref{varrho}--\eqref{omega}, we deduce from \eqref{l41a} that
\begin{align} &
\mu \partial^2_{3}\omega_{\mm{h}}=\mu\Delta \omega_{\mm{h}}={\mathbf{M}}_{\mm{h}}-{\mathbf{L}}_{\mm{h}}+ {\rho}  \partial_{t}\omega_{\mm{h}}+ {\rho} v\cdot  \nabla \omega_{\mm{h}} + \mathbf{N}_{\mm{h}}  =0 \text{ on }\partial \Omega.
\label{2505171219}
\end{align}		
By virtue of the boundary conditions of $(v_{\mm{h}}, \partial_3^2\omega_{\mm{h}})$ in  \eqref{omega} and  \eqref{2505171219}, and  the incompressible condition  \eqref{1a}$_3$, we  further obtain \begin{align}(\partial^3_3v_{\mm{h}},\partial^2_3v_{3}) |_{\partial\Omega}=0 ,\label{202223090s9315a4}
\end{align}		which implies
\begin{align}\partial^3_{3}\omega_{3}=\partial^3_{3}(\partial_{1}v_2-\partial_{2}v_1)=0 \text{ on } \partial \Omega.\label{202223090s93154}
\end{align}		
Thanks to \eqref{2505171219}--\eqref{202223090s93154}, we  arrive at \eqref{vbian} and \eqref{omega} for $i=1$.

Similarly to \eqref{2505171219}, we use \eqref{varrho} and \eqref{omega} with $i=0$, $1$ to deduce that
\begin{align} &
\mu \partial^4_{3}\omega_{\mm{h}}=\mu {\partial^2_3}\Delta \omega_{\mm{h}}=\partial^2_3({\rho}  \partial_{t}\omega_{\mm{h}}-{\mathbf{L}}_{\mm{h}}+ {\rho} v\cdot  \nabla \omega_{\mm{h}} + {\mathbf{M}}_{\mm{h}} + \mathbf{N}_{\mm{h}} ) =0 \text{ on }\partial \Omega.
\label{2505171saf219}
\end{align}In addition, exploiting {the boundary condition of $\partial_3^3v_{\mm{h}}$} in \eqref{vbian} and \eqref{2505171saf219}, we have \begin{align}(\partial^{5}_3v_{\mm{h}},\partial^{4}_3v_{3},\partial^5_{3}\omega_{3}) |_{\partial\Omega}=0.\label{202223090s9315a4x}
\end{align}		 	
Thanks to \eqref{2505171saf219} and \eqref{202223090s9315a4x}, we  arrive at \eqref{vbian} and \eqref{omega} for $i=2$.  This completes the proof.
\hfill $\Box$
\end{pf}
\begin{lem}	\label{201805141072}For any given $t\in [0,T]$, we have
\begin{align}
&\| \nabla v \|_{\ell,j}\lesssim
\begin{cases} \|(\nabla  v , \Delta\omega)\|_{\ell,0} & \hbox{for }\ell=0\mbox{ or } -s,\mbox{ and } j=2; \\
\|(\nabla v , \nabla\Delta^2\omega)\|_0 & \hbox{for }\ell=0,\ j=5, \end{cases}
\label{nav1}
\\
&\|  v \|_{\ell,j}\lesssim\begin{cases}\|( v , \nabla\omega)\|_{\ell,0} & \hbox{for }\ell=0\mbox{ or }-s,\
\mbox{ and } j=2; \\\|( v , \Delta^2\omega)\|_0 & \hbox{for }\ell=0,\ j=5, \end{cases} \label{nav111}
\\
&\|\nabla v_t \|_{3}\lesssim
\|(  \nabla v_t , \Delta^2v_t)\|_0, \label{nav111x} \\
&	\| \varrho\|_{\ell,j}\lesssim \begin{cases}
\|\nabla\Delta \tilde\varrho\|_{\ell,0} & \hbox{for }\ell=0\mbox{ or }-s, \mbox{ and } j=3; \\
\|\Delta^3\tilde\varrho\|_{0} & \hbox{for }\ell=0,\ j=6,
\label{nav1var}
\end{cases}
\end{align}
where the definition of $\tilde\varrho$ has been provided by \eqref{hftxha}. \end{lem}
\begin{pf}
Recalling the boundary condition of $ \omega $ in  \eqref{omega} with $i=0$ and $1$, we use both the elliptic estimates in Lemmas \ref{xfsddfs2212} and \ref{xfs05072212} to deduce that
$$\|\nabla \Delta^i\partial_k^i\Lambda^{\ell}_{\mm{h}}\omega\|_1\lesssim\|\nabla\Delta^i \partial_k^i \Lambda^{\ell}_{\mm{h}} \omega\|_0+ \| \Delta^{1+i} \partial_k^i\Lambda^{\ell}_{\mm{h}} \omega\|_0,$$
where $i=[j/2]-1$  and $k=1$, $2$.
Similarly, we also have
$$\| \nabla \Delta^i\partial_3^i\Lambda^{\ell}_{\mm{h}}\omega\|_1\lesssim\| \nabla\Delta^i\partial_3^i \Lambda^{\ell}_{\mm{h}}\omega\|_0+ \|  \Delta^{1+i}\partial_3^i\Lambda^{\ell}_{\mm{h}} \omega\|_0.$$
Putting the above  two estimates together yields
\begin{align}\label{20171124150s2}\| \nabla^{1+i} \Delta^i\Lambda^{\ell}_{\mm{h}}\omega\|_1\lesssim\| \nabla^{1+i}\Delta^i\Lambda^{\ell}_{\mm{h}}\omega\|_0+ \| \nabla^i \Delta^{1+i} \Lambda^{\ell}_{\mm{h}} \omega\|_0.\end{align}
Thus it follows from \eqref{20171124150s2} that
\begin{align}\|  \Lambda^{\ell}_{\mm{h}}\omega\|_2\lesssim\|   \nabla\Lambda^{\ell}_{\mm{h}} v\|_1+ \|  \Delta   \Lambda^{\ell}_{\mm{h}}\omega\|_0
\label{250831}
\end{align}
and
\begin{align}\| \Delta \omega\|_3\lesssim\|  \Delta  v\|_3+ \|  \nabla\Delta^2  \omega\|_0.\label{2508311}\end{align}

Making use of the incompressible condition  \eqref{1a}$_3$, the boundary conditions of $(v_3,\partial_3^2v_3)$ in \eqref{vbian}, the Hodge-type elliptic estimate \eqref{202005021302} and  the interpolation inequality \eqref{201807291850}, we can further deduce from
\eqref{250831}, resp. \eqref{2508311} that
\begin{align}\|  \Lambda^{\ell}_{\mm{h}}v\|_2\lesssim\|   \nabla\Lambda^{\ell}_{\mm{h}} v\|_0+ \|  \Delta   \Lambda^{\ell}_{\mm{h}}\omega\|_0
\label{2508171552}
\end{align}
resp.
\begin{align}\| \Delta  v\|_4\lesssim\|   \Delta v\|_0+ \|  \nabla \Delta^2 \omega\|_0.
\label{2508171601}
\end{align}
In particular, the estimate \eqref{2508171552} yields \eqref{nav1} with $\ell=0$ or $-s$, and $j=2$.

In addition,  in view of the elliptic estimates in Lemmas \ref{xfsddfs2212} and \ref{xfs05072212}, it holds that
\begin{align}\|  \nabla  v\|_5 \lesssim\| \nabla v\|_4+ \| \Delta  v\|_4,
\label{2508171602}
\end{align}
which, together with \eqref{2508171601}, yields \eqref{nav1} with
{$(\ell,j)=(0,5)$}.

Utilizing the Hodge-type elliptic estimate {\eqref{202005021302}}, the incompressible condition {\eqref{1a}$_3$},  and the interpolation inequality \eqref{201807291850},  we easily get
{\eqref{nav111}} with $\ell=0$ or $-s$, and $j=2$.

Similarly to \eqref{nav1} with $(\ell,j)=(0,5)$, it is easy to deduce that
\begin{align}
\label{08171621}
\| \nabla  v \|_4\lesssim \|(\nabla v, \Delta^2\omega)\|_0,
\end{align}
and
\begin{align} \label{081716212}
\|\nabla v_t \|_{3}\lesssim
\|( \nabla v_t , \Delta^2v_t)\|_0.
\end{align}
Hence \eqref{nav111} with $(\ell,j)=(0,5)$ holds by \eqref{08171621} and the interpolation inequality \eqref{201807291850}. In addition, using the Hodge-type elliptic estimate \eqref{202005021302} and the incompressible condition {\eqref{1a}$_3$ again}, we infer \eqref{nav111x} from \eqref{081716212}.
 		
Finally, recalling the definition of $\tilde{\varrho}$  in \eqref{hftxha} and the boundary conditions of $\varrho$ in \eqref{varrho},  we  easily follow the  argument of
\eqref{nav1} to obtain \eqref{nav1var}.
\hfill $\Box$
\end{pf}	
\subsection{A basic energy estimate and a stabilizing estimate}
Now we derive a basic energy estimate for the QRT problem.
\begin{lem}\label{varrho1}For any given $t\in [0,T]$, we have
\begin{align}
&\frac{\mm{d}}{\mm{d}t}\left( E\left( \varrho\right)+\| \sqrt{\rho}  v\|^2_{0} \right) +2\mu \|\nabla v \|^2_{ 0}\lesssim\sqrt{\mathfrak{E}}\mathfrak{D}.
\label{lemiq11x}
\end{align}
\end{lem}
\begin{pf}
Taking the inner product  of \eqref{1a}$_{2}$ and $v$ in $L^2$, and then making use of
the mass equation \eqref{1a}$_1$, the incompressible condition \eqref{1a}$_3$,
the boundary condition of $v_3$ in \eqref{n1}, the relation \eqref{07291549} and the
integration by parts, we  obtain		
\begin{align}	
&\frac{1}{2}\frac{\mm{d}}{\mm{d}t}\int {\rho} | v|^2\mm{d}x +\mu\int |\nabla v|^2\mm{d}x
\nonumber\\&=\int \left(\frac{{\varrho}_t v^2  }{2}- {\rho} v\cdot  \nabla v-g{\varrho}\mathbf{e}^3- {\varepsilon ^2}\mathbf{Q}^{\mm L}-{\varepsilon ^2}\mathbf{Q}^{\mm N}\right)\cdot v \mm{d}x
\nonumber\\&=\int \left(\frac{{\varrho}_t v^2  }{2}- {\rho} v\cdot  \nabla v-{\varepsilon ^2}\mathbf{Q}^{\mm N}\right)\cdot v \mm{d}x-\int \left( g{\varrho}v_3+{\varepsilon ^2}\mathbf{Q}^{\mm L}v\right)  \mm{d}x
\nonumber \\
&=-\int\left(g{\varrho}v_3+{\varepsilon ^2}\mathbf{Q}^{\mm L}\cdot v \right)    \mm{d}x-\int{\varepsilon ^2}\mathbf{Q}^{\mm N}\cdot v  \mm{d}x=: {\mathcal{I}}-\int{\varepsilon ^2}\mathbf{Q}^{\mm N}\cdot v   \mm{d}x.
\label{2022401sdfa211727}
\end{align}

Making use of {the mass equation \eqref{1a}$_1$}, the incompressible condition \eqref{1a}$_3$, the boundary condition of $v_3$ in \eqref{n1}, the integration by parts, and the relation
\begin{align}
\label{2509010855}
\mathbf{Q}^{\mm L}-\nabla \partial_3{\left(\bar \rho '/\bar \rho \right)} = \left(\frac{\bar \rho '}{\bar \rho }\Delta \varrho - \partial _3\left( \left|\frac{\bar \rho '}{\bar \rho }\right|^2\varrho  \right) - \left( \frac{\bar \rho '}{\bar \rho } \right) ''\varrho\right)\mathbf{e}^3=:\tilde{\mathbf{Q}}^{\mm L} ,
\end{align}
the integral $ {\mathcal{I}}$ can be rewritten as
\begin{align*}
{\mathcal{I}}	&= \int \left({\varepsilon ^2}\left(
{\partial _3}\left( {{\left| {\frac{{\bar \rho '}}{{\bar \rho }}} \right|^2}\varrho } \right)- \frac{{\bar \rho '}}{{\bar \rho }}\Delta \varrho - {\left( {\frac{{\bar \rho '}}{{\bar \rho }}} \right)}^\prime{\partial_3}\varrho  \right)v_3- g{\varrho}v_3-{\varepsilon ^2}\partial _3\left( {\left( {\frac{{\bar \rho '}}{{\bar \rho }}} \right)\nabla \varrho} \right)\cdot v \right)  \mm{d}x
\\&=\int \left({\varepsilon ^2}\left(\left(\frac{\bar \rho '}{\bar \rho } \right)''\varrho- \frac{{\bar \rho '}}{{\bar \rho }}\Delta \varrho  + {\partial _3}\left( {{\left| {\frac{{\bar \rho '}}{{\bar \rho }}} \right|^2}\varrho }\right)  \right)- g{\varrho}   \right)v_3  \mm{d}x
\\	&=\int \frac{1}{\bar{\rho}'}  \left(g{\varrho}-{\varepsilon ^2}\left( \left( \frac{\bar \rho ''}{\bar \rho } \right)^\prime{\varrho}+\left| \frac{\bar \rho '}{\bar \rho }\right|^2\partial _3\varrho -  \frac{{\bar \rho '}}{{\bar \rho }}\Delta \varrho \right)\right)\left(\varrho _t+v\cdot\nabla \varrho \right)  \mm{d}x
\\	&=-\frac{1}{2}\frac{\mm{d}}{\mm{d}t}E\left( \varrho\right)+\int \left(\frac{g-\varepsilon ^2 (\bar{\rho}''/\bar \rho)'}{\bar{\rho}'} {\varrho} + \left( \frac{{\varepsilon ^2}}{\bar \rho }\right)'\partial _3\varrho +\frac{\varepsilon ^2}{\bar \rho }\Delta\varrho   \right) v\cdot\nabla\varrho  \mm{d}x.
\end{align*}
Inserting the above identity into \eqref{2022401sdfa211727} yields			
\begin{align}
&{\frac{1}{2}}\frac{\mm{d}}{\mm{d}t}\left( E\left( \varrho\right)+\| \sqrt{\rho}  v\|^2_{0} \right) +{\mu}\|\nabla v \|^2_{ 0}=\tilde{\mathcal{I}}, \label{e011}
\end{align}
where we have defined that
\begin{align}\tilde{\mathcal{I}}:=\int \left( \left(\frac{g-\varepsilon ^2 (\bar{\rho}''/\bar \rho)'}{\bar{\rho}'} {\varrho} +\left( \frac{{\varepsilon ^2} }{\bar \rho }\right)'\partial _3\varrho  +\frac{ {\varepsilon ^2}}{\bar \rho }\Delta\varrho  \right) v\cdot\nabla\varrho -{\varepsilon ^2}\mathbf{Q}^{\mm N}\cdot v\right) \mm{d}x.
\label{2508111549}
\end{align}

In view of \eqref{imvarrho}, \eqref{im1a}, \eqref{esmmdforinfty} and the regularity of $\bar{\rho}$, it is easy to prove that
\begin{align}
& \left\|\left(\frac{\bar{\rho}'}{\bar \rho}\right)^j (\varrho + \bar{\rho})^{-1}f\right\| _{i}\lesssim \|f\|_i \mbox{ for }1\leqslant i\leqslant 5\mbox{ and }0\leqslant j\leqslant  2.\label{im3a}	\end{align}
{In addition, exploiting the boundary condition of $v_3$ and the incompressible condition  \eqref{1a}$_3$, we  get
\begin{align}
 \|  v_3\|_{i}\lesssim	\begin{cases}\|  v\|_{1,0}&\mbox{for }i=0;\\
\| \nabla v\|_{1,0}&\mbox{for }i=2, \end{cases}   \label{im3}
\end{align}
where we have used the  Poincar\'e-type inequality  \eqref{2202402040948}   in the case $i=1$ and the elliptic estimate  \eqref{xfsdsaf41252} ($ \|  v_3\|_{2}\lesssim  \|  \Delta v_3\|_{0}$) in the case $i=2$.
Then}, making use of  \eqref{im3a}, \eqref{im3}, the interpolation inequality \eqref{leml3}, H\"older's inequality and the regularity of $\bar{\rho}$,  we obtain		
\begin{align}
\tilde{\mathcal{I}}\lesssim&  \|\varrho\|_{W^{2,4}(\Omega)}  ( \| \varrho\|_{W^{1,4}(\Omega)} \|v_3\|_{0}+\|\varrho\|_{1,1}\|v_{\mm{h}}\|_{L^4}) \nonumber \\
\lesssim&\|\varrho\|_{1,2} \big(\|\varrho\|_{3}\left(\|\varrho\|_{3}\|\nabla v\|_{0} +\|\varrho\|_{1,1}\| v\|_{1}\right)\|\nabla v\|_{0}\big)^{\frac{1}{2}}
\lesssim\sqrt{\mathfrak{E}}\mathfrak{D}. \label{202221401177321aaa}
\end{align}
Putting the above estimate into \eqref{e011} yields  \eqref{lemiq11x}.	 This completes the proof. \hfill $\Box$
\end{pf}
	
It should be noted that $E(\varrho)$ in \eqref{lemiq11x} is positive. More precisely, we have the following
stabilizing estimate.
\begin{lem}\label{lem:20201224215}
 {Under the sharp stability condition \eqref{2saf01504}},  it holds that
\begin{align}
\|  r\|^2_1\lesssim E(r)\mbox{ for any }r\in H_0^1.\label{F}
\end{align}
\end{lem}
\begin{pf}
Thanks to the relation ${\varepsilon_{\mm{c}}}^2=ga_1$ in  \eqref{2022309032042}, we get
\begin{align}
g\int \bar{\rho}'\varpi^2\mm{d}{x}\leqslant {\varepsilon_{\mm{c}}}^2  \|\underline{ {\rho}}\nabla  \varpi  \|^2_0 \mbox{ for any } \varpi\in H_0^1.
\label{2022309032055}
\end{align}
Exploiting \eqref{2022309032055}, the sharp stability condition  \eqref{2saf01504}, the stabilizing condition  \eqref{2022205071434} and Poincar\'e-type inequality \eqref{2202402040948}, we obtain that
\begin{align}
\|\varpi\|_1^2\lesssim 	&	\| \nabla  \varpi \|_{0}^2\lesssim  (\varepsilon^2-{\varepsilon_{\mm{c}}}^2){\|\bar{\rho}'\nabla  \varpi \|_0^2}\nonumber \\
\lesssim &  \| \varepsilon\underline{\bar{\rho}}\nabla  \varpi \|_0^2 -g\|\sqrt{\bar{\rho}'}\varpi\|_0^2=c E_{\mm{L}}(\varpi) . \label{2022401311434}
\end{align}		
		
Let $\Upsilon =r/\bar{\rho}'$.
Exploiting the boundary condition of $r$ and the integration by parts, it is easy to compute out that 	
\begin{align}
&\int\left( \frac{|\nabla r|^2}{\bar \rho}+\frac{ (\bar{\rho}''/\bar \rho)'}{\bar{\rho}'}r^2\right) \mmd x\nonumber \\&=\int\Big{(}{\bar \rho}^{-1}\left( |  \bar \rho ' \nabla_{\mm{h}}\Upsilon   |^2+(\bar{\rho}''\Upsilon +\bar{\rho}'\partial_{3}\Upsilon )^2\right) + (\bar{\rho}''/\bar \rho)'\bar{\rho}'\Upsilon ^2\Big{)} \mmd x\nonumber \\&=\int {\bar \rho}^{-1}\left( |\bar \rho '\nabla_{\mm{h}}\Upsilon   |^2   + (\bar{\rho}''\Upsilon +\bar{\rho}'\partial_{3}\Upsilon )^2  - {|\bar \rho ''|^2} \Upsilon ^2- {2\bar \rho '\bar \rho ''} \Upsilon\partial_{3}\Upsilon  \right) \mmd x \nonumber \\&= \int \left|\underline{\bar{\rho}}\nabla  \Upsilon\right| ^2 {\mm{d}x}.  \label{2022401311434a}
\end{align}
Thus, we derive from  \eqref{2022401311434} and \eqref{2022401311434a}   with $\varpi= \Upsilon$ that
\begin{align*}
&  \|   r/\bar{\rho}' \|^2_1  =\|  \Upsilon\|^2_1  \lesssim E_{\mm{L}}( \Upsilon)=  E( r) ,
\end{align*}
which, together with  Poincar\'e-type inequality \eqref{2202402040948}, the sharp stability condition,  the stabilizing condition and Young's inequality, yields \eqref{F}. This completes the proof. \hfill $\Box$
\end{pf}

\subsection{Estimates of high-order spacial derivatives}
This section is devoted to establishing {the} estimates of the {high-order} spacial derivatives  of $(\varrho,v)$ and {the} dissipative estimate of $v_t$ in the following two lemmas resp..
\begin{lem}\label{lem1a}
For any given $t\in\overline{I_T}$, we have
\begin{align}
& \frac{\mm{d}}{\mm{d}t}\bigg(\varepsilon ^2\| {{\underline{\bar{\rho}}}}\nabla^{k}\Delta^{1+j}\tilde{\varrho}\|^2_{0} +\|\sqrt{{\rho}}\nabla^{k}\Delta^{j}\omega\|^2_{0} \nonumber\\
&-(2+2j+k)\varepsilon ^2\int \frac{{\underline{\bar{\rho}}^2}}{\bar{\rho}'}
\partial_{3}\varrho(\partial^{2+2j+k}_3 {\tilde{\varrho}})^2\mmd x
\bigg)+c\|\nabla^{k}\cdot\nabla \Delta^{j}\omega\|^2_{0}\nonumber \\
& \lesssim 	\begin{cases}
(\|\varrho\|_{1,1}+\|v_t\|_{1})\|\nabla v\|_2+\sqrt{\mathfrak{E}}  \mathfrak{D}+ \mathfrak{T}(t) &\mbox{for }\,k=1,\ j=0
;\\	
(\|\varrho\|_{1,4}+\|v_t\|_{4})\|\nabla v\|_5+\sqrt{\mathcal{E}}  \mathcal{D}+\mathfrak{N}(t)
&\mbox{for }\,k=0,\ j=2,
\end{cases}\label{202221401177321a}
\end{align}
where the definitions of ${\underline{\bar{\rho}}}$, $\mathfrak{N}(t)$  and  $\tilde{\varrho}$
can be found in \eqref{20224013022226}, \eqref{msaft} and \eqref{hftxha}   resp., and
we have defined that
\begin{equation}
\label{202508091522}
\mathfrak{T}(t):= \sqrt{\mathfrak{E}}
\left( \|\nabla v\|_{1,0} \sqrt{\mathfrak{E}}+\|\nabla v\|_2{\mathfrak{E}}+(\|\nabla v\|_{1,0}\|\nabla v\|_{1})^{\frac{1}{2}}(\mathfrak{E} \mathfrak{D})^{\frac{1}{4}}\right).
\end{equation}
\end{lem}
\begin{pf}
Multiplying the mass equation \eqref{1a}$_1$ by $1/\bar{\rho}'$, and then applying ${\underline{\bar{\rho}}^2}\nabla^{k}\Delta^{1+j}$ to the resulting identity, we have
\begin{align}
{\underline{\bar{\rho}}^2} \nabla^{k}\Delta^{1+j} (\tilde{\varrho}_t+v_3-\underline{\varrho}v_3+v\cdot\nabla \tilde{\varrho})=0.
\label{l4asf1a}
\end{align}	
Here and in what follows, we have defined that  $\underline{\varrho}:=(1/\bar{\rho}')'\varrho$. Obviously
\begin{align}
{\|\underline{\varrho}\|_{i,l}\lesssim}\|\varrho\|_{i,l}\lesssim \|\tilde{\varrho}\|_{i,l}\lesssim \|\varrho\|_{i,l}\mbox{ for }i+l\leqslant 6.
\label{l4assaff1a}
\end{align}

Taking the inner product of \eqref{l4asf1a} and $\varepsilon ^2  \nabla^{k}\Delta^{1+j}
\tilde{\varrho}$  in $L^2$ yields
\begin{align}	
&\frac{1}{2}\frac{\mm{d}}{\mm{d}t} \| \varepsilon \underline{\bar{\rho}} \nabla^{k}\Delta^{1+j} \tilde{\varrho}\|^2_{0}
=\varepsilon ^2\int{{\underline{\bar{\rho}}^2}}\big( \nabla^{k}\Delta^{1+j}( \underline{\varrho}v_3- v_3-v\cdot\nabla \tilde{\varrho}) \nabla^{k}\Delta^{1+j}\tilde{\varrho}\big)\mm{d}x. 	\label{e21}
\end{align}	
Applying  $\Delta ^{j}$ to  the vortex equation \eqref{l41a} yields
\begin{align}
&\Delta ^{j} \left( {\rho}  \omega_t \right)  =\Delta ^{j} \left(\mu\Delta \omega +{\mathbf{L}} - {\mathbf{M}}-   {\rho} v\cdot  \nabla \omega - \mathbf{N}\right).	\label{l41acc1}	\end{align}	
Taking the inner product  of the above identity
and $\Delta^{k+j}\omega $ in $L^2$, and then exploiting the integration by parts   and the  boundary condition of $(\omega_{\mm{h}},\partial_3^2\omega_{\mm{h}},\partial_3^{2i+1}\omega_3) $ in \eqref{omega}, we  arrive at	\begin{align}	
&\frac{1}{2}\frac{\mm{d}}{\mm{d}t}\| \sqrt{{\rho}} \nabla^{k}\Delta^{j}\omega  \|^2_{0} +\mu\|\nabla^{k}\cdot\nabla \Delta^{j}\omega \|^2_{0}
\nonumber\\&=
\int  \bigg(  \left(\Delta^{j}( {\mathbf{L}}-{\mathbf{M}}_{\mm{h}}-{\rho}
v\cdot  \nabla \omega+\mathbf{N} )\right)\cdot\Delta^{k+j}\omega
+k(
\partial_i\rho \omega_t\cdot \partial_i\omega + 2^{-1}\varrho_t|\nabla  \omega|^2)\nonumber\\&\quad  +\frac{j}{2}\left(\rho\Delta^2\omega_t-\Delta^2(\rho\omega_t) +\frac{\varrho_t \Delta^2\omega  }{2} \right)\cdot\Delta^2\omega  \bigg)   \mm{d}x.	\label{e22}
\end{align}
Adding \eqref{e21} and \eqref{e22} together,  we get
\begin{align}
&\frac{1}{2}\frac{\mm{d}}{\mm{d}t}\left( \| \varepsilon\underline{\bar{\rho}}\nabla^{k}\Delta^{1+j}\tilde{\varrho}\|^2_{0} +\|\sqrt{{\rho}}\nabla^{k}\Delta^{j}\omega\|^2_{0}
\right)+\mu\|\nabla^{k}\cdot\nabla \Delta^{j}\omega\|^2_{0}=\sum_{m=1}^3 	\mathcal{J}^{k,j}_{m},\label{1e21}
\end{align}		
where we have defined that
\begin{align*}
\mathcal{J}^{k,j}_{1}:=&(-1)^k\int\bigg( \Delta^{j}( {\mathbf{L}}-{\mathbf{M}} ) \cdot\Delta^{k+j}\omega-\varepsilon ^2{\underline{\bar{\rho}}^2} \nabla^{k}\Delta^{1+j} v_3\cdot\nabla^{k}\Delta^{1+j}\tilde{\varrho}
\\&+k \bar{\rho}'\omega_t\cdot\partial_{3}\omega+\frac{j}{2}
\left(\bar{\rho}\Delta^2\omega_t -\Delta^2(\bar{\rho}\omega_t)\right)\cdot\Delta^2\omega\bigg)  \mm{d}x,\\
\mathcal{J}^{k,j}_{2}:=	& 	\int  \bigg( k\left( \frac{\varrho_t|\nabla  \omega|^2 }{2}-\partial_i {\varrho}\partial_t \omega \cdot \partial_i\omega \right) -   \Delta ^{j}  ({\rho} v\cdot  \nabla \omega+\mathbf{N} ) \cdot\Delta^{k+j}\omega \\&\quad+\frac{j}{2}\left(\varrho\Delta^2\omega_t
-\Delta^2(\varrho\omega_t) +\frac{\varrho_t \Delta^2\omega }{2} \right)\cdot\Delta^2\omega  \bigg)   \mm{d}x,	 \\
\mathcal{J}^{k,j}_{3}:=& \varepsilon ^2\int  {\underline{\bar{\rho}}^2} \left( \nabla^{k}\Delta^{1+j}\left( \tilde\varrho v_3 -v\cdot\nabla \tilde{\varrho}  \right)\right) \nabla^{k}\Delta^{1+j}\tilde{\varrho}  \mm{d}x. \end{align*}
Next we estimate for $	\mathcal{J}^{k,j}_{1}$, $	\mathcal{J}^{k,j}_{2}$ and $	\mathcal{J}^{k,j}_{3}$ in sequence.
		
(1) Exploiting	 the boundary  conditions  of $(\partial_3^{2i}\tilde{\varrho},\partial_3^{2i}\omega_{\mm{h}})$ in \eqref{omega},  the integration by parts and the relation $$
\Delta v_3=\partial_{2}\omega_1-\partial_{1}\omega_2,$$  the integral $\mathcal{J}^{k,j}_{1}$ can be estimated as follows
\begin{align}
\mathcal{J}^{k,j}_{1}=& \int\bigg( k\left((( {\mathbf{L}}- \varepsilon ^2{\underline{\bar{\rho}}^2} \left(-\partial _2,\partial _1,0\right)^\top\Delta  \tilde\varrho )-{\mathbf{M}}) \cdot\Delta\omega+ \bar{\rho}'\omega_t\cdot\partial_{3}\omega\right)\nonumber
\\&-\frac{j}{2}  \bigg(\nabla\Delta( {\mathbf{L}}- \varepsilon ^2{\underline{\bar{\rho}}^2} \left(-\partial _2,\partial _1,0\right)^\top\Delta  \tilde\varrho -{\mathbf{M}}):\nabla\Delta^{2}\omega\nonumber\\&+
\left(\bar{\rho}\Delta^2\omega_t -\Delta^2(\bar{\rho}\omega_t)\right)\cdot\Delta^2\omega\bigg)  \bigg)  \mm{d}x\nonumber
\\\lesssim&	\begin{cases}
(\|\varrho\|_{1,1}+\|v_t\|_{1})\|\nabla v\|_2&\mbox{for }k=1,\ j=0;
\\(\|\varrho\|_{1,4}+\|v_t\|_{4})\|\nabla v\|_5 &\mbox{for }k=0,\ j=2.
 \end{cases}\label{fkj1}
\end{align}		
		
(2)  Recalling the expressions of $(\mm{curl}\mathbf{Q}^{\mm N})_3$  in \eqref{qn3a} and $\mathbf{N}_3^{\mm{m}}$ in \eqref{curla}, we derive from  the boundary conditions of  $\partial_3^{2i}(\varrho,v)$ in \eqref{varrho} and \eqref{vbian} that \begin{align} \partial_{3}\Delta(\mathbf{N}_3^{\mm{m}}+(\mm{curl}\mathbf{Q}^{\mm N})_3)|_{\partial\Omega}=0.
\label{2022240123120405zc}
\end{align}
 Thanks to the expressions of $(\mm{curl}\mathbf{Q}^{\mm N})_i$ with $i=1,2$ and $3$  in \eqref{qn3a}, we can compute out that
\begin{align}
\|\mm{curl}\mathbf{Q}^{\mm N}\|_{i}\lesssim &	\begin{cases}
\|\varrho\|_{3}\|\varrho\|_{1,2} &\mbox{for }i=0;\\\|\varrho\|_{6}\|\varrho\|_{1,5}	&\mbox{for }i=3.\label{qn3}
\end{cases}
\end{align}

Making use of \eqref{im3a}, \eqref{2022240123120405zc}, the boundary conditions of $\partial_3^{2i}\omega_{\mm{h}}$ in \eqref{omega},  the product estimates in \eqref{fgestims} and the integration by parts,  the integral $\mathcal{J}^{k,j}_{2}$ can be  controlled as follows
\begin{align}
\mathcal{J}^{k,j}_{2}\lesssim & \big( (1+\|\varrho\|_{2+2j+k})\|v\|_{1+2j+k}\|\omega\|_{1+2j+k}+\|\varrho_t\|_{1+2j+k}\|v\|_{1+2j+k}\nonumber \\&+\|\varrho\|_{2+2j+k}\|v_t\|_{2j+k}\big)\|\nabla v\|_{1+2j+k}
+k \left| \int   (\mathbf{N}^{\mm{m}}+\mm{curl}\mathbf{Q}^{\mm N})   \cdot \Delta\omega    \mm{d}x\right|\nonumber\\&+ j \left| \int   \nabla\Delta(\mathbf{N}^{\mm{m}}+ \mm{curl}\mathbf{Q}^{\mm N})   :  \nabla\Delta^2\omega    \mm{d}x\right|\nonumber \\
\lesssim &\big( (1+\|\varrho\|_{2+2j+k})\|v\|_{1+2j+k}\|\omega\|_{1+2j+k}+\|\varrho_t\|_{1+2j+k}\|v\|_{1+2j+k}\nonumber
\\&+\|\varrho\|_{2+2j+k}(\|\varrho\|_{1,1+2j+k}+\|v_t\|_{2j+k})\big)\|\nabla v\|_{1+2j+k}\nonumber \\\lesssim &	\begin{cases}
\sqrt{\mathcal{E}}\mathcal{D} &\mbox{for }k=0,\ j=2;\\\sqrt{\mathfrak{E}}\mathfrak{D}	&\mbox{for }k=1,\ j=0,
\end{cases}  \label{fkj2}
\end{align}
where we have used \eqref{qn3} in the second inequality, and \eqref{imvarrho} in the last inequality.

(3) To estimate for the integral $	\mathcal{J}^{k,j}_{3}$, we first use the integration by parts to split it into three terms:
\begin{align*}
\mathcal{J}^{k,j}_{3}
&= \varepsilon ^2 \sum_{1\leqslant m\leqslant 3}	\mathcal{J}^{k,j}_{3,m},
\end{align*}
where we have defined that
\begin{align}
\mathcal{J}^{k,j}_{3,1} := &   -\left( 2+2j+k\right)\int  {\underline{\bar{\rho}}^2}  \partial_3  v_3(\partial^{2+2j+k}_3 {\tilde{\varrho}} )^2\mm{d}x, \label{fkj31a}\\
\mathcal{J}^{k,j}_{3,2} :=&\begin{cases}
-\int  {\underline{\bar{\rho}}^2}((\partial_{l}  v \cdot\nabla\Delta {\tilde{\varrho}}	+  2(\partial_{i}  v: \nabla\partial_{l}\partial_{i}   {\tilde{\varrho}}+\partial_{l} \partial_{i}  v\cdot \nabla\partial_{i}   {\tilde{\varrho}} )+\Delta  v \cdot\nabla\partial_{l} {\tilde{\varrho}}\\+ \partial_{l} \Delta v \cdot\nabla {\tilde{\varrho}}
)\partial_{l} \Delta_{\mm{h}} \tilde{\varrho}+ (\nabla_{\mm{h}}  v \cdot\nabla\Delta {\tilde{\varrho}}	+  2\nabla_{\mm{h}} (\partial_{i}  v\cdot \nabla\partial_{i}   {\tilde{\varrho}})
\\+\nabla_{\mm{h}}(\Delta v \cdot\nabla {\tilde{\varrho}}))\nabla_{\mm{h}} \partial^2_{3}    \tilde{\varrho}+(
  \partial_{3}  v_{\mm{h}} \cdot\nabla_{\mm{h}}\Delta {\tilde{\varrho}}+\partial_{3}  v_{3} \partial_{3}\Delta_{\mm{h}} {\tilde{\varrho}}\\
+  2\partial_{3}(\partial_{i}  v_{\mm{h}}\cdot \nabla_{\mm{h}}\partial_{i}   {\tilde{\varrho}}+\nabla_{\mm{h}} v_{3}\cdot \nabla_{\mm{h}}\partial_{3}   {\tilde{\varrho}})+\partial_{3}(\Delta v_{\mm{h}} \cdot\nabla_{\mm{h}} {\tilde{\varrho}})) \partial^3_{3} \tilde{\varrho}) \mm{d}x&\mbox{for }k=1,\ j=0;\\
\int {\underline{\bar{\rho}}^2} (( v_{\mm{h}}\cdot \nabla_{\mm{h}}\Delta^3\tilde{\varrho}
 - \Delta^3(v_{\mm{h}}\cdot \nabla_{\mm{h}}\tilde{\varrho}))\Delta^3\tilde{\varrho}+
\Delta^3(\underline{\varrho}v_3&
\\-\partial_{3}\tilde{\varrho}v_3)(\partial^2_{3}\Delta_{\mm{h}}\Delta\tilde{\varrho}	+\partial^4_{3}\Delta_{\mm{h}}\tilde{\varrho}+\Delta_{\mm{h}}\Delta^2\tilde{\varrho})
+v_3\partial_{3}\Delta^3\tilde{\varrho}(\Delta_{\mm{h}}\Delta^2\tilde{\varrho}&\\+\partial^2_{3}\Delta_{\mm{h}}\Delta\tilde{\varrho}	+\partial^4_{3}\Delta_{\mm{h}}\tilde{\varrho})) \mm{d}x&\mbox{for }k=0,\ j=2,
\end{cases}\nonumber \\  	\mathcal{J}^{k,j}_{3,3}:= &   \int(({\underline{\bar{\rho}}})' v_3\nabla^{k}\Delta^{1+{j}}\tilde{\varrho} \cdot \nabla^{k}\Delta^{1+{j}}\tilde{\varrho}+{\underline{\bar{\rho}}^2} ( v_3\partial^{1+k}_3\Delta^{1+j} {\tilde{\varrho}}+  \partial^k_3\Delta^{1+j}(({\underline{\varrho}}-\partial_3 {\tilde{\varrho}} )v_3)\nonumber \\
& \quad +(2+2j+k) \partial_3  v_3\partial^{2+2j+k}_3 {\tilde{\varrho}} )  \partial^{2+2j+k}_3 {\tilde{\varrho}} )\mm{d}x.\label{2.59}
\end{align}
It should be noted that we have used the Einstein convention of summation over repeated indices, where $1\leqslant i$, $l\leqslant 3$.
Next we estimate for $\mathcal{J}^{k,j}_{3,1}$, $\mathcal{J}^{k,j}_{3,2}$ and $	\mathcal{J}^{k,j}_{3,3}$ in sequence.
		
(i) Applying $\partial_{3}$ to \eqref{1a}$_{1}$  yields
\begin{align*}
\bar{\rho}'\partial_{3}v_3+\bar{\rho}''v_3+
\partial_{3}(\varrho _t+v\cdot\nabla \varrho)=0.
\end{align*}
Thanks to the above identity, we can compute out that
\begin{align}
\mathcal{J}^{k,j}_{3,1}&=(2+2j+k)\int\frac{{\underline{\bar{\rho}}^2} }{\bar{\rho}'} \left( \bar{\rho}''v_3+
\partial_{3}(\varrho _t+v\cdot\nabla \varrho)\right)(\partial^{2+2j+k}_3 {\tilde{\varrho}} )^2\mmd x\nonumber\\
&=\left( 2+2j+k\right)\sum_{1\leqslant i\leqslant 3} {\mathcal{J}}^{k,j}_{3,1,i}, \label{fkj12}
\end{align}
where we have defined that
\begin{align}
&{\mathcal{J}}^{k,j}_{3,1,1}:=\frac{\mmd}{\mmd t}\int \frac{{\underline{\bar{\rho}}^2}}{\bar{\rho}'}
\partial_{3}\varrho(\partial^{2+2j+k}_3 {\tilde{\varrho}})^2\mmd x,\label{fkj311}\\
&{\mathcal{J}}^{k,j}_{3,1,2}:=-\int\frac{{\underline{\bar{\rho}}^2} }{\bar{\rho}'}\left( 2
\partial_{3}\varrho\partial^{2+2j+k}_3 {\tilde{\varrho}}\partial^{2+2j+k}_3 {\tilde{\varrho}}_t-v\cdot\partial_3\nabla \varrho(\partial^{2+2j+k}_3\tilde{\varrho})^2\right) \mmd x,\label{fkj312}\\
&{\mathcal{J}}^{k,j}_{3,1,3}:= \int \frac{{\underline{\bar{\rho}}^2}}{\bar{\rho}'}\left(  {\bar{\rho}''} v_3+\partial_3v\cdot\nabla \varrho  \right)(\partial^{2+2j+k}_3 {\tilde{\varrho}})^2\mmd x .\label{fkj313}
\end{align}	
Next we shall estimate for ${\mathcal{J}}^{k,j}_{3,1,2}$ and ${\mathcal{J}}^{k,j}_{3,1,3}$.		
		
Making use of \eqref{l4assaff1a}, the product estimates in \eqref{fgestims}, the incompressible condition {\eqref{1a}$_3$}, the boundary condition of $v_3$ and the integration by parts, ${\mathcal{J}}^{k,j}_{3,1,2}$ can be further rewritten as follows
\begin{align*}
{\mathcal{J}}^{k,j}_{3,1,2}&=\int\frac{{\underline{\bar{\rho}}^2}}{\bar{\rho}'}\left( 2\left(
\partial_{3}\varrho\partial^{2+2j+k}_3\tilde{\varrho}\partial^{2+2j+k}_3(v_3-\underline{\varrho}v_3+v\cdot\nabla \tilde{\varrho})
\right) +v\cdot\partial_3\nabla \varrho(\partial^{2+2j+k}_3\tilde{\varrho})^2\right) \mmd x\nonumber
\\
&={\mathcal{J}}^{k,j}_{3,1,2,1}+{\mathcal{J}}^{k,j}_{3,1,2,2},
\end{align*}
where we have defined that
\begin{align}
{\mathcal{J}}^{k,j}_{3,1,2,1}&:=2\int\frac{{\underline{\bar{\rho}}^2}}{\bar{\rho}'}
\partial_{3}\varrho\partial^{2+2j+k}_3\tilde{\varrho}\partial^{2+2j+k}_3v_3\mmd x,
\nonumber\\
{\mathcal{J}}^{k,j}_{3,1,2,2}&:=2\int\bigg(\frac{{\underline{\bar{\rho}}^2}}{\bar{\rho}'}
\partial_{3}\varrho\partial^{2+2j+k}_3\tilde{\varrho}\left(\partial^{2+2j+k}_3(v\cdot\nabla \tilde{\varrho}) -v\cdot\nabla \partial^{2+2j+k}_3\tilde{\varrho}\right.\nonumber\\
&\quad\left.-\partial^{2+2j+k}_3(\underline{\varrho} v_3) \right)-\partial_3{\varrho} (\partial^{2+2j+k}_3\tilde{\varrho})^2\mm{div}\left( \frac{{\underline{\bar{\rho}}^2}v}{\bar{\rho}'}\right)\bigg)\mmd x.\nonumber
\end{align}
Using  the incompressible condition, the boundary conditions of $(v,\partial^{2i}_3\tilde{\varrho})$ in \eqref{n1} and \eqref{vbian}, the product estimates in  \eqref{fgestims} and the integration by parts, we obtain
\begin{align}
{\mathcal{J}}^{k,j}_{3,1,2,1}&=-2\int\frac{{\underline{\bar{\rho}}^2}}{\bar{\rho}'}
\partial_{3}\varrho\partial^{2+2j+k}_3\tilde{\varrho}\partial^{1+2j+k}_3 \mm{div}_{\mm{h}} v_{\mm{h}} \mmd x
\nonumber\\&=2\int\frac{{\underline{\bar{\rho}}^2}}{\bar{\rho}'}\left(
\partial^{2+2j+k}_3\tilde{\varrho} \partial^{1+2j+k}_3v_{\mm{h}}\cdot\nabla _{\mm{h}}\partial_{3}\varrho+\partial_{3}\varrho\partial^{1+2j+k}_3v_{\mm{h}}\cdot \nabla _{\mm{h}}\partial^{2+2j+k}_3\tilde{\varrho}\right) \mmd x
\nonumber	\\&=2\int\left(\frac{{\underline{\bar{\rho}}^2}}{\bar{\rho}'}
\partial^{1+2j+k}_3v_{\mm{h}}\cdot \partial^{2+2j+k}_3\tilde{\varrho}\nabla _{\mm{h}}\partial_{3}\varrho-\partial_3\left(  \frac{{\underline{\bar{\rho}}^2}}{\bar{\rho}'}\partial_{3}\varrho\partial^{1+2j+k}_3v_{\mm{h}}\right) \cdot\nabla _{\mm{h}}\partial^{1+2j+k}_3\tilde{\varrho}\right) \mmd x
\nonumber	\\& 	\lesssim\|\varrho\|_{1,1+2j+k}\|\varrho\|_{2+2j+k}\|\nabla v\|_{1+2j+k}.
\label{20224012292110vA}
\end{align}	
In addition, utilizing  \eqref{imvarrho}, \eqref{im3}, the interpolation inequality \eqref{leml3} and the product estimates in \eqref{fgestims}, we obtain
\begin{align}
{\mathcal{J}}^{k,j}_{3,1,2,2}&\lesssim\|\varrho\|^2_{2+2j+k}\left(\|v_3\|_{2}+\|\varrho\|_{3}\|\nabla v\|_2 \right)+\|\varrho\|_{1,1+2j+k}\|\varrho\|_{2+2j+k}\|\nabla v\|_{1+2j+k}\nonumber\\&\quad +\|\varrho\|_{3}\|\varrho\|_{2+2j+k}\left( \|\partial_{3}\tilde{\varrho}\partial_{3}v_3\|_{1+2j+k}+\|\underline{\varrho}v_3\|_{2+2j+k}\right)  \nonumber	\\&\lesssim\|\varrho\|_{2+2j+k}\big{(}\|\varrho\|_{2+2j+k}\| v\|_{1,1}+\|\varrho\|_{3} ( \|\varrho\|_{2+2j+k}\|\nabla v\|_2\nonumber\\&\quad  +\|\varrho\|_{3}\|\nabla v\|_{1+2j+k}) +  \|\varrho\|_{1,1+2j+k}\|\nabla v\|_{1+2j+k}\big{)}.\label{20224012292110vAC}
\end{align}	
			
Similarly to \eqref{20224012292110vAC}, we can also have
\begin{align}		{\mathcal{J}}^{k,j}_{3,1,3}&\lesssim\|\varrho\|^2_{2+2j+k}\left(\|v\|_{1,1}+\|\varrho\|_{3} \|\nabla v\|_2\right).\label{20224012292110vACc}
\end{align}			
Therefore,  we immediately deduce from \eqref{fkj12}-- \eqref{20224012292110vACc} that
\begin{align}
&\mathcal{J}^{k,j}_{3,1}-(2+2j+k)\frac{\mmd}{\mmd t}\int \frac{{\underline{\bar{\rho}}^2}}{\bar{\rho}'}
\partial_{3}\varrho(\partial^{2+2j+k}_3 {\tilde{\varrho}})^2\mmd x\nonumber
\\&\lesssim\|\varrho\|_{2+2j+k}\Big(\big(\|\varrho\|_{2+2j+k}\|\nabla v\|_{1,0}+\|\varrho\|_{3}( \|\varrho\|_{2+2j+k}\|\nabla v\|_2\nonumber\\&\quad+\|\varrho\|_{3}\|\nabla v\|_{1+2j+k}) \big) + \|\varrho\|_{1,1+2j+k}\|\nabla v\|_{1+2j+k}\Big).
\label{fkj31}
\end{align}			
	
(ii) Arguing  analogously to  \eqref{20224012292110vAC}, we  have		
\begin{align}
&\mathcal{J}^{k,j}_{3,2}\lesssim\|\varrho\|_{1,1+2j+k}\|\varrho\|_{2+2j+k}\|\nabla v\|_{1+2j+k}.
\label{fkj32}
\end{align}			
		
(iii) Now we  rewrite	$\mathcal{J}^{k,j}_{3,3}$ as follows
\begin{align}
\mathcal{J}^{k,j}_{3,3}&=	\mathcal{J}^{k,j}_{3,3,1}+	\mathcal{J}^{k,j}_{3,3,2},\nonumber
\end{align}
where	we have defined that			
\begin{align*}
\mathcal{J}^{k,j}_{3,3,1}:=&\begin{cases}
-\int  {\underline{\bar{\rho}}^2} \partial^3_3 {\tilde{\varrho}} (3\partial^2_{3} v_3(\partial_{3}\underline{\varrho}+\partial^2_{3}\tilde{\varrho})+
\partial_{3}\Delta v_3(\underline{\varrho}+\partial_{3}\tilde{\varrho})) \mm{d}x&\mbox{for }k=1,\ j=0;\\-\int  {\underline{\bar{\rho}}^2} \partial^6_3 {\tilde{\varrho}} ((\partial^4_{3}\underline{\varrho}+\partial^5_{3}\tilde{\varrho})(12\partial^2_{3}v_3+3\Delta v_3) \\
+(\partial^3_{3}\underline{\varrho}+\partial^4_{3}\tilde{\varrho})(8\partial^3_{3} v_3+12\partial_{3}\Delta v_3) \nonumber\\+(\partial^2_{3}\underline{\varrho}+\partial^3_{3}\tilde{\varrho})(12\partial^2_{3}\Delta v_3+2\Delta^2 v_3) \nonumber\\+6(\partial_{3}\underline{\varrho}+\partial^2_{3}\tilde{\varrho})\partial_{3}\Delta^2 v_3+(\underline{\varrho}+\partial_{3}\tilde{\varrho})
\Delta^3 v_3)\mm{d}x&\mbox{for }k=0,\ j=2
\end{cases}
\end{align*}
and
\begin{align*}
\mathcal{J}^{k,j}_{3,3,2}:=&\begin{cases}
-\int  {\underline{\bar{\rho}}^2} \partial^3_3 {\tilde{\varrho}} (3\partial^2_{3}\underline{\varrho}\partial_{3} v_3+\partial^3_{3}\underline{\varrho} v_3)\mm{d}x&\mbox{for }k=1,\ j=0;\\
-\int  {\underline{\bar{\rho}}^2} \partial^6_3 {\tilde{\varrho}} (6\partial^5_{3}\underline{\varrho}\partial_{3} v_3+ \partial^6_{3}\underline{\varrho}v_3) \mm{d}x&\mbox{for }k=0,\ j=2.
\end{cases}
\end{align*}
	
In view of the boundary conditions of $
(\partial_3^{2i}\tilde{\varrho},\partial_3v_{\mm{h}},v_3)$ in \eqref{varrho} and \eqref{vbian}, we can follow the argument of the derivation of \eqref{20224012292110vA} to establish
\begin{align}
\mathcal{J}^{k,j}_{3,3,1}&	\lesssim\|\varrho\|_{1,1+2j+k}\|\varrho\|_{2+2j+k}   \|\nabla v\|_{1+2j+k}. \label{fkj331}
\end{align}	
Following the argument of  \eqref{202221401177321aaa}, we can derive  that
\begin{align}
\mathcal{J}^{k,j}_{3,3,2}&\lesssim\|\varrho\|_{2+2j+k}\left( \|\varrho\|_{2+2j+k}\|v_3\|_{2} +\|\partial^{1+2j+k}_3\underline{\varrho}\|_{L^4}\|\partial_3v\|_{L^4}\right) \nonumber
\\&\lesssim\|\varrho\|_{2+2j+k}\left(\|\varrho\|_{2+2j+k} \|\nabla v\|_{1,0} +(\|\varrho\|_{1,1+2j+k}\|\varrho\|_{2+2j+k}\|\nabla v\|_{1,0}\|\nabla v\|_{1})^{\frac{1}{2}}\right).   \label{fkj332}
\end{align}	
Thus  it follows  from  \eqref{fkj331} and \eqref{fkj332}
that
\begin{align}
\mathcal{J}^{k,j}_{3,3}	&\lesssim\|\varrho\|_{2+2j+k}\Big(  \|\varrho\|_{1,1+2j+k}\|\nabla v\|_{1+2j+k}+\|\varrho\|_{2+2j+k} \|\nabla v\|_{1,0}\nonumber
\\&\quad+(\|\varrho\| _{1,1+2j+k}\|\varrho\|_{2+2j+k}\|\nabla v\|_{1,0}\|\nabla v\| _{1})^{\frac{1}{2}}\Big).
\label{fkj33}
\end{align}
Consequently, we  deduce from \eqref{fkj31}, \eqref{fkj32} and \eqref{fkj33} that
\begin{align}
&\mathcal{J}^{k,j}_{3}-\varepsilon ^{2}(2+2j+k)\frac{\mmd}{\mmd t}\int \frac{{\underline{\bar{\rho}}^2}}{\bar{\rho}'}
\partial_{3}\varrho(\partial^{2+2j+k}_3 {\tilde{\varrho}})^2\mmd x\nonumber
\\&\lesssim	\begin{cases}
\sqrt{\mathfrak{E}}\left( \mathfrak{D}+\|\nabla v\|_{1,0}\sqrt{\mathfrak{E}}+\|\nabla v\|_2 {\mathfrak{E}} +(\|\nabla v\|_{1,0}\|\nabla v\|_{1})^{\frac{1}{2}}(\mathfrak{E}\mathfrak{D})^{\frac{1}{4}}\right) 	&\mbox{for }k=1,\ j=0;\\	\sqrt{\mathcal{E}}\bigg( \mathcal{D}+\|\nabla v\|_{1,0}\sqrt{\mathcal{E}}+\|\varrho\|_{3}\left( \|\nabla v\|_2 \sqrt{\mathcal{E}}+\|\varrho\|_{3} \sqrt{\mathcal{D}}\right) &\\+(\|\nabla v\|_{1}\|\nabla v\|_{1,0})^{\frac{1}{2}}(\mathcal{E}\mathcal{D})^{\frac{1}{4}}\bigg) &\mbox{for }k=0,\ j=2.
\end{cases}
\label{fkj3}
\end{align}	
	
Finally, inserting \eqref{fkj1}, \eqref{fkj2} and \eqref{fkj3} into \eqref{1e21}, we  arrive at \eqref{202221401177321a}. This completes the proof.
\hfill $\Box$	
\end{pf}
\begin{lem}\label{lemt}
For any given $t\in\overline{I_T}$, we have
\begin{align}
& \frac{\mm{d}}{\mm{d}t}\left(
2\varepsilon ^2\int\Delta^2\left(  \frac{\bar{\rho}'}{\bar{\rho}}\Delta{\varrho}-\left( \frac{\bar{\rho}'}{\bar{\rho}}\right)'\partial_{3}{\varrho}\right) \Delta^2  v_3\mm{d}x +\mu\|\nabla \Delta^2 v\|_0^2\right)+c\| \Delta^2  v_t\|_{0}^2\nonumber \\
&\lesssim \|\varrho\|_{1,4}^2+\|  v_3\|_5^2+\|  v_t\|_0^2+\| \nabla v_t\|_2^2+ \sqrt{\mathcal{E}} \mathcal{D} \label{2022sfa2140117732}. \end{align}
\end{lem}	
\begin{pf}
Applying $\Delta^2 $ to \eqref{1a}$_2$ yields
\begin{align}
{\rho}\Delta^2 v_{t}=&\Delta^2 (\mu\Delta v-g{\varrho}\mathbf{e}^3-\nabla  \beta- {\varepsilon ^2}\mathbf{Q}-   \rho v\cdot  \nabla v)
- (\Delta^2 ({\rho} v_{t}) -{\rho}\Delta^2 v_{t})  . 	
\label{20220020913}
\end{align}
Utilizing  the boundary conditions of $(v,\partial^3_{3}v_{\mm{h}})$ in   \eqref{vbian} and  the incompressible condition \eqref{1a}$_{3}$, we obtain
\begin{align}
\partial_3 \Delta^2    v_{1}=		\partial_3 \Delta^2    v_{2} = \Delta^2    v_3=0 \mbox{ on } 	\partial\Omega.\label{vtt2}
\end{align}
	
Taking the inner product of \eqref{20220020913} and $\Delta^2 v_{t} $
in $L^2$, and then exploiting  the mass equation \eqref{1a}$_{1}$, \eqref{vtt2}, the incompressible condition \eqref{1a}$_{3}$ and the integration by parts, we get
\begin{align}
&\frac{\mu  }{2} \frac{\mm{d}}{\mm{d}t}  \|\nabla \Delta^2  v\|_0^2 +\|\sqrt{\rho}\Delta^2  v_t\|_0^2= \sum_{i=1}^3\mathcal{K}_i,
\label{2024012200334}
\end{align}
where we have defined that
\begin{align*}
&\mathcal{K}_1:=-\int \Delta^2( g{\varrho} \mathbf{e}^3+ \varepsilon ^2\tilde{\mathbf{Q}})\cdot \Delta^2 v_t\mm{d}x,\\
&\mathcal{K}_2:=\varepsilon ^2\int \Delta^2\left(\left(\left( \frac{\bar{\rho}'}{\bar{\rho}}\right)'\partial_{3}{\varrho}- \frac{\bar{\rho}'}{\bar{\rho}}\Delta{\varrho}\right) \mathbf{e}^3-\nabla\left( \frac{\bar{\rho}'}{\bar{\rho}}\partial_{3}{\varrho}\right) \right)\cdot  \Delta^2  v_t\mm{d}x,\\
&\mathcal{K}_3:=\int \big( {\rho}\Delta^2 v_{t}-  \Delta^2(\rho v\cdot  \nabla v
+ {\rho} v_{t})  \big) \cdot \Delta^2 v_t\mm{d}x,\\
&\tilde{\mathbf{Q}}:={\mathbf{Q}}-{\bar{\rho}'}(\Delta{\varrho}\mathbf{e}^3+ \nabla\partial_{3}{\varrho})/{\bar{\rho}}.
\end{align*}
We next  estimate for the above three integrals $\mathcal{K}_1$--$\mathcal{K}_3$ in sequence.
	
Exploiting  the product estimates in \eqref{fgestims},  the boundary condition of $\Delta^2 v_3$ in \eqref{vtt2}, H\"older's inequality, the incompressible condition and the integration by parts, the integral $\mathcal{K}_1$  can be controlled by
\begin{align}
\mathcal{K}_1 =&  \int  \left( \varepsilon ^2  \nabla_{\mm{h}}\Delta\partial_3 (   g{\varrho} +\varepsilon ^2\tilde{\mathbf{Q}}_3 ) \cdot\partial_t \partial^2_3 \Delta v_{\mm{h}}-\Delta_{\mm{h}}\Delta( g{\varrho} \mathbf{e}^3 \right. \nonumber
\\
&\quad \left.+\varepsilon ^2\tilde{\mathbf{Q}} )\cdot \Delta^2 v_t- \varepsilon ^2\partial_3^2   \Delta\tilde{\mathbf{Q}}_{\mm{h}} \cdot \partial_t \Delta^2 v_{\mm{h}}  \right) \mm{d}x \nonumber \\
\lesssim & (\|\varrho\|_{1,4}+ \|\varrho\|_{1,5}\|\varrho\|_6)\|\Delta^2v_t\|_{0}.\label{20220020a913}
\end{align}		
Using the mass equation \eqref{1abx}$_1$, the boundary condition of $\Delta^2 v_3$ in \eqref{vtt2} and the incompressible condition, we have
\begin{align*}
\mathcal{K}_2 =& \varepsilon ^2\int  \nabla\Delta \left( \frac{\bar{\rho}'}{\bar{\rho}} \Delta(\bar{\rho}'v_3+v\cdot\nabla \varrho )\right) \cdot \nabla \Delta^2 v_3\mm{d}x  \\
&-\varepsilon ^2\int  \nabla\Delta \left( \left( \frac{\bar{\rho}'}{\bar{\rho}}\right)'\partial_{3}(\bar{\rho}'v_3+v\cdot\nabla \varrho )\right) \cdot \nabla \Delta^2 v_3\mm{d}x\nonumber \\
&	-\varepsilon ^2\frac{\mm{d}}{\mm{d}t}  \int \Delta^2\left( \frac{\bar{\rho}'}{\bar{\rho}}\Delta{\varrho}-\left( \frac{\bar{\rho}'}{\bar{\rho}}\right)'\partial_{3}{\varrho}\right) \Delta^2  v_3\mm{d}x\nonumber \\
\leqslant& c(\|v_3\|_5+\|\varrho\|_{1,5}\|v\|_5+\|\varrho\|_{6}\|v_3\|_5)\|v_3\|_5\nonumber \\
&-\varepsilon ^2\frac{\mm{d}}{\mm{d}t}  \int \Delta^2\left( \frac{\bar{\rho}'}{\bar{\rho}}\Delta{\varrho}-\left( \frac{\bar{\rho}'}{\bar{\rho}}\right)'\partial_{3}{\varrho}\right) \Delta^2  v_3\mm{d}x .\end{align*}
	Similarly to \eqref{20220020a913}, we obtain
\begin{align*}
\mathcal{K}_3\lesssim  \big( (1+\| \varrho \|_6)\|v\|_5\|\nabla v\|_5 +(1+\|\varrho\|_6)(\| v_{t}\|_0+\|\nabla v_{t}\|_2)\big) \|\Delta^2v_t\|_0 .\end{align*}

Finally, putting the above three estimates into \eqref{2024012200334},  and then using \eqref{im1}, \eqref{im3} and Young's inequality, we can obtain \eqref{2022sfa2140117732}. This completes the proof.		
\hfill $\square$ \end{pf}

\subsection{Tangential estimates}
This section is devoted to building some  estimates for the tangential derivatives of $(\varrho,v)$.
\begin{lem}
\label{lem2dxfa}Let $a$ and $s$ be the constants appearing in Theorem \ref{thm2}. For any given $t\in\overline{I_T}$, we have
\begin{align}
&\sup_{0\leqslant  t\leqslant  T}\left\lbrace  \langle t\rangle^{\ell_1}(\|\varrho\|^2_{{\ell_2},1}+\|v\|^2_{{\ell_2},0})\right\rbrace\nonumber \\
&  +\int_0^{T}\langle t\rangle^{\ell_1}\left( \| \nabla v\|^2_{{\ell_2},0}
-c \langle t\rangle^{-1}(\| \varrho\|^{ {2}}_{\ell_2,1}+\| v\|^2_{{\ell_2},0})
\right) \mm{d}t\nonumber\\&\lesssim \|\varrho^0\|^2_{{\ell_2},1}+\|v^0\|^2_{{\ell_2},0}+ 	B^{3/2}(T) \hbox{ for } \begin{cases}
{\ell_1}=0,\ {\ell_2}=-s;\\	
{\ell_1}=3a/2,\ {\ell_2}=1,
\end{cases}	
\label{dctv1}  \\
& \sup_{0\leqslant  t\leqslant  T} \big{\{}\langle t\rangle^{\ell_1}(  \|\nabla v\|^2_{{\ell_2},0}-c\|\varrho\|_{\ell_2,0}\|v\|_{\ell_2,1}  ) \big{\}}\nonumber\\&+\int_{0}^{T}\langle t\rangle^{\ell_1}\big(\|v_t\|^2_{{\ell_2},0}- c \langle t\rangle^{-1} \|\nabla v\|^2_{{\ell_2},0}\big)\mm{d}t \nonumber \\
&\lesssim\int_0^{T} \langle t\rangle^{\ell_1}\left( \|\varrho_{t}\|^2_{{\ell_2},1}+ \langle t\rangle^{-1}\| v\|_{{\ell_2},1}\| \varrho\|_{{\ell_2},1}\right) \mm{d}t\nonumber\\&\quad +\|(\varrho^0,v^0)\|_{\ell_2,1}^2+	 B^{3/2}(T)\label{dctv3} \hbox{ for }\begin{cases}
{\ell_1}=0,\ {\ell_2}=-s;\\
{\ell_1}=a\hbox{ or }0,\ {\ell_2}=0,
\end{cases}\\
& \sup_{0\leqslant  t\leqslant  T}\left\lbrace  \langle t\rangle^{\ell_1}\|(\nabla \Delta \tilde{\varrho}, \nabla\omega  )\|^2_{{\ell_2},0}\right\rbrace \nonumber\\& +\int_0^{T}\langle t\rangle^{\ell_1}\left(
 \|\Delta \omega\|^2_{{\ell_2},0}-c\langle t\rangle^{-1}(\|  \nabla \Delta \tilde{\varrho} \|^2_{\ell_2,0}+\| \nabla\omega  \|^2_{\ell_2,0}) \right)\mm{d}t\nonumber \\&	
\lesssim 	\int_0^{T}\langle t\rangle^{\ell_1}\big((\|\varrho\|_{1+{\ell_2},1}+\| v_t\|_{{\ell_2},1})\|\nabla v\|_{{\ell_2},2}\big)\mm{d}t
\nonumber \\&\quad 	+
\|(\nabla {\varrho}^0,v^0  )\|^2_{{\ell_2},2}+	 B^{3/2}(T)	 \hbox{ for }
{\ell_1}=0,\ {\ell_2}=-s  \label{dctv2}
\end{align}	
and	\begin{align}
& \sup_{0\leqslant  t\leqslant  T}\left\lbrace \langle t\rangle^{\ell_1} \left( \| \nabla \omega\|^2_{\ell_2,0}- c\| \Delta v_3\|_{{\ell_2},0}\|\Delta \tilde{\varrho}\|_{{\ell_2},0}\right)\right\rbrace  + \int_0^{T}\langle t\rangle^{\ell_1}\| \omega_t\|^2_{{\ell_2},0}\mm{d}t\nonumber
\\&\lesssim\int_0^{T}\big( \langle t\rangle^{\ell_1}(\|\varrho\|^2_{1+{\ell_2},2}+\|\varrho_{t}\|^2_{{\ell_2},2}+\|v_t\|^2_{{\ell_2},0}
 )+ \langle t\rangle^{ -1}(\| \varrho\|_{1+{\ell_2},2}\| \nabla v\|_{{\ell_2},0}\nonumber
\\&\quad +\|\nabla v\|^2_{{\ell_2},1})\big) \mm{d}t+\|(\varrho^0,v^0)\|_{\ell_2,2}^2 +	 B^{3/2}(T)\label{dctv4}\hbox{ for }\begin{cases}
{\ell_1}=0,\ {\ell_2}=-s;\\
{\ell_1}=a\hbox{ or }0,\ {\ell_2}=0.
\end{cases}
\end{align}		
\end{lem}
\begin{pf}
Applying $\Lambda^{\ell_2}_{\mm{h}}$ to the mass equation \eqref{1a}$_1$ and the momentum equations in \eqref{1a}$_2$, one has
\begin{equation}\label{h1f}
\begin{cases}
\Lambda^{\ell_2}_{\mm{h}} (\varrho _t+\bar{\rho}'v_3 +v\cdot\nabla \varrho)=0,\\	
\Lambda^{\ell_2}_{\mm{h}} (\bar{\rho} v_{t} +\nabla  \beta) =
\Lambda^{\ell_2}_{\mm{h}} (\mu\Delta v-g{\varrho}\mathbf{e}^3- {\varepsilon ^2}\mathbf{Q}-{\varrho} v_{t}-   \rho v\cdot  \nabla v). 		\end{cases}
\end{equation}
Taking the inner product of \eqref{h1f}$_1$ and
$$-\langle t\rangle^{\ell_1}\Lambda^{\ell_2}_{\mm{h}}\left(\frac{g-\varepsilon ^2 (\bar{\rho}''/\bar \rho)'}{\bar{\rho}'} {\varrho} +  \left( \frac{\varepsilon ^2}{\bar \rho }\right)'\partial _3\varrho+ \frac{\varepsilon ^2}{\bar \rho }\Delta\varrho \right)$$ in $L^2$ yields
\begin{align}	
&\frac{1}{2}\left( \frac{\mm{d}}{\mm{d}t}	\left( \langle t\rangle^{\ell_1}E(\Lambda^{\ell_2}_{\mm{h}}\varrho)  \right)  -{\ell_1}\langle t\rangle^{{\ell_1}-1}E(\Lambda^{\ell_2}_{\mm{h}}\varrho)\right) \nonumber\\&
= \int \langle t\rangle^{\ell_1}\Lambda^{\ell_2}_{\mm{h}}\left(\frac{g-\varepsilon ^2 (\bar{\rho}''/\bar \rho)'}{\bar{\rho}'} {\varrho} + \left( \frac{{\varepsilon ^2}}{\bar \rho }\right)'\partial _3\varrho+ \frac{{\varepsilon ^2}}{\bar \rho }\Delta\varrho   \right)\Lambda^{\ell_2}_{\mm{h}}(\bar{\rho}'v_3+ v\cdot \nabla \varrho)\mm{d}x. \label{h1f1} \end{align}
	
Taking the inner product of  \eqref{h1f}$_2$ and  $	\langle t\rangle^{\ell_1}\Lambda^{\ell_2}_{\mm{h}}v$ in $L^2$, and then using  the mass equation \eqref{1a}$_1$,  the incompressible condition \eqref{1a}$_3$, the boundary condition of $v_3$ in \eqref{n1}, \eqref{2509010855} and the integration by parts, we can obtain	
\begin{align}	
&\frac{1}{2}\left( \frac{\mm{d}}{\mm{d}t}\int	\langle t\rangle^{\ell_1} {\bar \rho}|\Lambda^{\ell_2}_{\mm{h}}v|^2\mm{d}x-{\ell_1}\langle t\rangle^{{\ell_1}-1}\|\sqrt{\bar\rho} \Lambda^{\ell_2}_{\mm{h}}v\|^2_{0}\right) +\mu\langle t\rangle^{\ell_1}\|\nabla v\|^2_{{\ell_2},0}\nonumber\\& = -\int\langle t\rangle^{\ell_1} \Lambda^{\ell_2}_{\mm{h}} \left({\varrho} v_{t}+g{\varrho}\mathbf{e}^3+  \rho v\cdot  \nabla v+{\varepsilon ^2}\tilde{\mathbf{Q}}^L  \right)\cdot\Lambda^{\ell_2}_{\mm{h}} v \mm{d}x
.\label{h1f2}
\end{align}
Thus, it follows from \eqref{h1f1} and \eqref{h1f2} that
\begin{align}	
&\frac{1}{2}\left( \frac{\mm{d}}{\mm{d}t}\left(\langle t\rangle^{\ell_1} \left( E(\Lambda^{\ell_2}_{\mm{h}}\varrho)+\|\sqrt{\bar\rho} \Lambda^{\ell_2}_{\mm{h}}v\|^2_{0}\right)\right)-{\ell_1}\langle t\rangle^{{\ell_1}-1} \left( E(\Lambda^{\ell_2}_{\mm{h}}\varrho)+\|\sqrt{\bar\rho} \Lambda^{\ell_2}_{\mm{h}}v\|^2_{0}\right)\right) \nonumber\\&+\mu\langle t\rangle^{\ell_1}\|\nabla v\|^2_{{\ell_2},0} ={\mathcal{L}}_1 +{\mathcal{L}}_2,
\label{1ee1}
\end{align}
where we have defined that
\begin{align}
{\mathcal{L}}_1:=&\int\langle t\rangle^{\ell_1} \Lambda^{\ell_2}_{\mm{h}}\left(\frac{g-\varepsilon ^2 (\bar{\rho}''/\bar \rho)'}{\bar{\rho}'} {\varrho}+  \left( \frac{\varepsilon ^2}{\bar \rho }\right)'\partial _3\varrho +  \frac{\varepsilon ^2}{\bar \rho }\Delta\varrho \right)\Lambda^{\ell_2}_{\mm{h}}( v\cdot \nabla \varrho)\mm{d}x,
\nonumber\\
{\mathcal{L}}_2:=	& -\int\langle t\rangle^{\ell_1}\Big( {\varepsilon ^2}\left( \Lambda^{\ell_2}_{\mm{h}} \mathbf{Q}^{\mm N}_{\mm{h}}\cdot \Lambda^{\ell_2}_{\mm{h}}v_{\mm{h}}+\Lambda^{\ell_2}_{\mm{h}} \mathbf{Q}^{\mm N}_{3} \Lambda^{\ell_2}_{\mm{h}}v_{3}\right)+\Lambda^{\ell_2}_{\mm{h}}(\varrho v_{t}+  \rho v\cdot  \nabla v)\cdot \Lambda^{\ell_2}_{\mm{h}}v \Big) \mm{d}x.\nonumber
\end{align}
		
   	Taking inner products of
\eqref{1a}$_{1}$ resp. \eqref{1a}$_{2}$ and
$$\langle t\rangle^{\ell_1}\Lambda^{\ell_2}_{\mm{h}}\partial_t\bigg(\frac{g-\varepsilon ^2 (\bar{\rho}''/\bar \rho)'}{\bar{\rho}'} {\varrho} + \left( \frac{{\varepsilon ^2}}{\bar \rho }\right)'\partial _3\varrho + \frac{{\varepsilon ^2}}{\bar \rho }\Delta\varrho  \bigg)$$ resp. $\langle t\rangle^{\ell_1}\Lambda^{\ell_2}_{\mm{h}}v_t$ in $L^2$, then summing up the resulting identities, and finally using \eqref{2509010855}, the   boundary condition of $v_3$, the incompressible condition and the integration by parts, we arrive at
\begin{align}	
& \frac{\mm{d}}{\mm{d}t} \left(\frac{\mu\langle t\rangle^{\ell_1}}{2}\|\nabla v\|^2_{{\ell_2},0}+\int\langle t\rangle^{\ell_1} \left( g\Lambda^{\ell_2}_{\mm{h}}{\varrho}\Lambda^{\ell_2}_{\mm{h}}v_3+{\varepsilon ^2}\Lambda^{\ell_2}_{\mm{h}}\mathbf{Q}^{\mm L}\cdot \Lambda^{\ell_2}_{\mm{h}}v\right)  \mm{d}x\right) +\langle t\rangle^{\ell_1}\| \sqrt{{\bar\rho}}\Lambda^{\ell_2}_{\mm{h}}v_t\|^2_{0}\nonumber
\\&	=\langle t\rangle^{\ell_1}E(\Lambda^{\ell_2}_{\mm{h}}\varrho_{t})+\frac{\mu {\ell_1}\langle t\rangle^{{\ell_1}-1}}{2}\|\nabla v\|^2_{{\ell_2},0}\nonumber
\\&\quad+{\ell_1}\int\langle t\rangle^{{\ell_1}-1} \left( g\Lambda^{\ell_2}_{\mm{h}}{\varrho}\Lambda^{\ell_2}_{\mm{h}}v_3+{\varepsilon ^2}\Lambda^{\ell_2}_{\mm{h}}\tilde{\mathbf{Q}}^{\mm L}\cdot \Lambda^{\ell_2}_{\mm{h}}v\right)  \mm{d}x+{\mathcal{L}}_{3}+{\mathcal{L}}_{4},\label{1ee3}
\end{align}			
where we have defined that		
\begin{align}
{\mathcal{L}}_{3}:=&  \int\langle t\rangle^{\ell_1}\Lambda^{\ell_2}_{\mm{h}}\partial_t \left(\frac{\varepsilon ^2 (\bar{\rho}''/\bar \rho)'-g}{\bar{\rho}'} {\varrho } - \left( \frac{ {\varepsilon ^2}}{\bar \rho }\right)'\partial _3\varrho  -\frac{ {\varepsilon ^2}}{\bar \rho }\Delta\varrho \right) \Lambda^{\ell_2}_{\mm{h}}\left(v\cdot\nabla \varrho \right)\nonumber  \mm{d}x,\\	{\mathcal{L}}_{4}:=&-\int\langle t\rangle^{\ell_1} \Lambda^{\ell_2}_{\mm{h}}\left(\varrho v_t+{\rho} v\cdot  \nabla v+ {\varepsilon ^2}\mathbf{Q}^{\mm N}\right) \cdot \Lambda^{\ell_2}_{\mm{h}}v_t \mm{d}x.\nonumber
\end{align}

{Applying ${\underline{\bar{\rho}}^2}\Lambda^{\ell_2}_{\mm{h}}\nabla\Delta$ to the mass equation \eqref{1a}$_1$  and $\Lambda^{\ell_2}_{\mm{h}}$ to the vortex equations \eqref{l41a} yield}
\begin{equation}\label{1aaaB}
\begin{cases}
{\underline{\bar{\rho}}^2}\Lambda^{\ell_2}_{\mm{h}}\nabla\Delta (\tilde{\varrho}_t+v_3+v\cdot\nabla \tilde{\varrho}-\underline{\varrho}v_3)=0, \\
\Lambda^{\ell_2}_{\mm{h}} (\bar{\rho} \omega_t+\varrho \omega_t+   {\rho} v\cdot  \nabla \omega )
=\Lambda^{\ell_2}_{\mm{h}}  (\mu\Delta \omega +\mathbf{L} -  {\mathbf{M}}-\mathbf{N} ).
\end{cases}
\end{equation}
Taking the inner products of \eqref{1aaaB}$_1$ resp.  \eqref{1aaaB}$_2$ and $\varepsilon ^2\langle t\rangle^{\ell_1}\Lambda^{\ell_2}_{\mm{h}}\nabla\Delta\tilde{\varrho}$ resp.  $ -\langle t\rangle^{\ell_1}\Lambda^{\ell_2}_{\mm{h}}\Delta\omega $    in $L^2$,  and then making use of the boundary conditions of $v$ in  \eqref{omega}, the integration by parts and the mass equation  {\eqref{1a}$_1$}, we  obtain
\begin{align}
&\frac{1}{2}\bigg( \frac{\mm{d}}{\mm{d}t}\left( \langle t\rangle^{\ell_1}{(\|\varepsilon \underline{\bar{\rho}}\Lambda^{\ell_2}_{\mm{h}}\nabla \Delta \tilde{\varrho}\|^2_{0} +\|\sqrt{\bar{\rho}} \Lambda^{\ell_2}_{\mm{h}}\nabla\omega  \|^2_{0})} \right)\nonumber\\& -{\ell_1}\langle t\rangle^{{\ell_1}-1}(\|\varepsilon \underline{\bar{\rho}}\Lambda^{\ell_2}_{\mm{h}}\nabla \Delta \tilde{\varrho}\|^2_{0} +\|\sqrt{\bar{\rho}} \Lambda^{\ell_2}_{\mm{h}}\nabla\omega  \|^2_{0})\bigg)+\mu\langle t\rangle^{\ell_1}\|\Delta \omega\|^2_{{\ell_2},0}=\sum_{i=5}^7{\mathcal{L}}_{ i} ,\label{1ee2}
\end{align}		
where we have defined that
\begin{align*}
{\mathcal{L}}_5:=& \int\langle t\rangle^{\ell_1}\left( \Lambda^{\ell_2}_{\mm{h}} (\mathbf{M} -\mathbf{L})\cdot \Lambda^{\ell_2}_{\mm{h}}\Delta\omega -\bar{\rho}'\Lambda^{\ell_2}_{\mm{h}}\omega_t\Lambda^{\ell_2}_{\mm{h}}\partial_{3}\omega-\varepsilon ^2{\underline{\bar{\rho}}^2} \Lambda^{\ell_2}_{\mm{h}}\nabla\Delta v_3\cdot\Lambda^{\ell_2}_{\mm{h}}\nabla\Delta \tilde{\varrho}\right)   \mm{d}x,\\
{\mathcal{L}}_6:=	&  \int\langle t\rangle^{\ell_1} \Lambda^{\ell_2}_{\mm{h}}\left( \varrho  \omega_t+   {\rho} v\cdot  \nabla \omega +\mathbf{N} \right)\cdot  \Lambda^{\ell_2}_{\mm{h}}\Delta\omega    \mm{d}x,	 \\
{\mathcal{L}}_7:=& \varepsilon ^2\int\langle t\rangle^{\ell_1}  {\underline{\bar{\rho}}^2} \Lambda^{\ell_2}_{\mm{h}}\nabla \Delta\left(  \underline{\varrho}v_3 - v\cdot\nabla \tilde{\varrho}  \right)\cdot\Lambda^{\ell_2}_{\mm{h}}\nabla\Delta\tilde{\varrho}  \mm{d}x. \end{align*}

Applying ${\underline{\bar{\rho}}^2}\Lambda^{\ell_2}_{\mm{h}}\Delta$ to  \eqref{1a}$_{1}$, we get
\begin{align}\label{t1}
{\underline{\bar{\rho}}^2}\Lambda^{\ell_2}_{\mm{h}}\Delta (\tilde{\varrho}_t+v_3-\underline{\varrho}v_3+v\cdot\nabla \tilde{\varrho})=0.
\end{align}
Taking inner products of \eqref{1aaaB}$_2$
resp. 	\eqref{t1}  and $ \langle t\rangle^{\ell_1}\Lambda^{\ell_2}_{\mm{h}}\omega _t$ resp. $\varepsilon ^2\langle t\rangle^{\ell_1}\Lambda^{\ell_2}_{\mm{h}}\Delta\tilde{\varrho}_t$ in $L^2$,  and then using the boundary conditions of $\omega$  in \eqref{omega}, the integration by parts and the mass equation  \eqref{1a}$_1$, we  deduce that
\begin{align}	
& \frac{\mm{d}}{\mm{d}t}\left( \langle t\rangle^{\ell_1} \left(\frac{\mu}{2}\|\Lambda^{\ell_2}_{\mm{h}}\nabla \omega\|^2_{0}+\varepsilon ^2\int  {\underline{\bar{\rho}}^2} \Lambda^{\ell_2}_{\mm{h}}\Delta \tilde{\varrho}\Lambda^{\ell_2}_{\mm{h}}\Delta v_3 \mm{d}x\right)\right)  +\langle t\rangle^{\ell_1}\|\sqrt{\bar{\rho}} \Lambda^{\ell_2}_{\mm{h}}\omega_t\|^2_{0}\nonumber
\\&	=\langle t\rangle^{{\ell_1}}\left({\ell_1}\langle t\rangle^{ -1} \left(\frac{\mu}{2}\|\Lambda^{\ell_2}_{\mm{h}}\nabla \omega\|^2_{0}+\varepsilon ^2\int  {\underline{\bar{\rho}}^2}\Lambda^{\ell_2}_{\mm{h}}\Delta \tilde{\varrho} \Lambda^{\ell_2}_{\mm{h}}\Delta v_3 \mm{d}x\right) + \|\varepsilon \underline{\bar{\rho}}\Lambda^{\ell_2}_{\mm{h}}\Delta \tilde{\varrho}_t\|^2_{0} \right) \nonumber \\
&\quad +\sum_{i=8}^{10}{\mathcal{L}}_{i},\label{1ee4}
\end{align}	
where we have defined that
\begin{align*}
&{\mathcal{L}}_8:= \int\langle t\rangle^{\ell_1}\left(  \Lambda^{\ell_2}_{\mm{h}}(\mathbf{L} -\mathbf{M})\cdot \Lambda^{\ell_2}_{\mm{h}}\omega_t\right)   \mm{d}x,\\
&{\mathcal{L}}_9:=	-  \int \langle t\rangle^{\ell_1}\Lambda^{\ell_2}_{\mm{h}}\left( \varrho  \omega_t+   {\rho} v\cdot  \nabla \omega +\mathbf{N} \right)\cdot  \Lambda^{\ell_2}_{\mm{h}}\omega_t    \mm{d}x,	 \\
&{\mathcal{L}}_{10}:= \varepsilon ^2\int \langle t\rangle^{\ell_1} {\underline{\bar{\rho}}^2}\Lambda^{\ell_2}_{\mm{h}} \Delta\left(\underline{\varrho}v_3 -v\cdot\nabla \tilde{\varrho} \right)\Lambda^{\ell_2}_{\mm{h}}\Delta\tilde{\varrho}_t \mm{d}x. \end{align*}
	
Integrating the above four identities \eqref{1ee1}, \eqref{1ee3}, \eqref{1ee2} and \eqref{1ee4} over $(0,T)$, resp., and then making use of the
\emph{a priori} assumption \eqref{imvarrho}, the bounds of density in \eqref{im1}, the stabilizing estimate \eqref{F}, the boundary condition of $v_3$, the incompressible condition, the integration by parts, the interpolation inequality \eqref{201807291850} and H\"older's inequality, we arrive at
\begin{align}
& \sup_{0\leqslant  t\leqslant  T}\left\lbrace  \langle t\rangle^{\ell_1}(\|\varrho\|^2_{{\ell_2},1}+\|v\|^2_{{\ell_2},0})\right\rbrace +\int_0^{T}\langle t\rangle^{\ell_1}\left( \| \nabla v\|^2_{{\ell_2},0}
-c\langle t\rangle^{-1}(\| \varrho\|^{{2}}_{\ell_2,1}+\| v\|^2_{{\ell_2},0})
\right) \mm{d}t \nonumber \\
& \lesssim \|\varrho^0\|^2_{{\ell_2},1}+\|v^0\|^2_{{\ell_2},0} +\int_0^{T}(|{\mathcal{L}}_1|+|{\mathcal{L}}_2|)\mm{d}t,
\label{dctv1q} \\
&\sup_{0\leqslant  t\leqslant  T} \big{\{}\langle t\rangle^{\ell_1}\left(  \|\nabla v\|^2_{{\ell_2},0}-c\|\varrho\|_{\ell_2,0}\|v\|_{\ell_2,1}  \right)\big{\}} +\int_{0}^{T}\langle t\rangle^{\ell_1}\big(\|v_t\|^2_{{\ell_2},0}- c\langle t\rangle^{-1} \|\nabla v\|^2_{{\ell_2},0}\big)\mm{d}t \nonumber \\& \lesssim
 \|( {\varrho}^0,v^0  )\|^2_{{\ell_2},1}
+\int_0^{T}\Big(|{\mathcal{L}}_{ 3}|+|{\mathcal{L}}_{4}|+\langle t\rangle^{\ell_1}( \langle t\rangle^{-1}\| \varrho\|_{{\ell_2},2}\| v\|_{{\ell_2},0}+ \|\varrho_{t}\|^2_{{\ell_2},1} )\Big) \mm{d}t{\color{red},} \label{dctv3q}
\\
& \sup_{0\leqslant  t\leqslant  T}\left\lbrace  \langle t\rangle^{\ell_1}\|(\nabla \Delta \tilde{\varrho}, \nabla\omega  )\|^2_{{\ell_2},0}\right\rbrace + \int_0^{T}\langle t\rangle^{\ell_1}\Big(\mu\|\Delta \omega\|^2_{{\ell_2},0}-{\ell_1}\langle t\rangle^{-1}(\|\varepsilon \underline{\bar{\rho}}\Lambda^{\ell_2}_{\mm{h}}\nabla \Delta \tilde{\varrho} \|^2_{0} \nonumber
\\&+\|\sqrt{\bar{\rho}} \Lambda^{\ell_2}_{\mm{h}}\nabla\omega  \|^2_{0})\Big)\mm{d}t\lesssim  \|(\nabla {\varrho}^0,v^0  )\|^2_{{\ell_2},2}+\sum_{i=5}^7\int_0^{T}|{\mathcal{L}}_{ i} |\mm{d}t,\label{dctv2q}
\end{align}
and
\begin{align}
& \sup_{0\leqslant  t\leqslant  T}\big{\{} \langle t\rangle^{\ell_1} \left( \| \nabla \omega\|^2_{\ell_2,0}- c\| \Delta v_3\|_{{\ell_2},0}\|\Delta \tilde{\varrho}\|_{{\ell_2},0}\right)\big{\}}  + \int_0^{T}\langle t\rangle^{\ell_1}\| \omega_t\|^2_{{\ell_2},0}\mm{d}t\nonumber\\&\lesssim \|(  {\varrho}^0,v^0  )\|^2_{{\ell_2},2}+\sum_{i=8}^{10}\int_0^{T}|{\mathcal{L}}_{ i} |\mm{d}t+\int_0^{T} \langle t\rangle^{\ell_1}\big{(}\|\varrho_t\|^2_{{\ell_2},2}\nonumber\\&\quad  + { \ell_1}\langle t\rangle^{ -1}(\|\varrho\|_{1+{\ell_2},2}\| \nabla v\|_{{\ell_2},0}+\|\nabla v\|^2_{{\ell_2},1} )\big{)}\mm{d}t.\label{dctv4q}
\end{align}	
To obtain the desired estimates in Lemma \ref{lem2dxfa}, next it suffices to estimate for ${\mathcal{L}}_{ 1} $--${\mathcal{L}}_{10}$.

(1) {Let us first consider the case that ${\ell_1}=0$ and ${\ell_2}=-s$}.
Recalling the definition of $\mathfrak{U}^{i}(\cdot)$ in \eqref{ugg} and the boundary condition of $v_3$, we obtain
	\begin{align}
\mathfrak{U}^{0}(v_3) &\lesssim\left(\|v_3\|_{1,1}\|v_3\|^{s}_{-s,1}\right) ^{\frac{1}{1+s}}\lesssim
\left( \|\Delta v_3\|^{s}_{-s,0}\sqrt{\underline{\mathcal{D}}}\right)^{\frac{1}{1+s}} \lesssim
\left( \|\nabla v\|^{s}_{1-s,0}\sqrt{\underline{\mathcal{D}}}\right)^{\frac{1}{1+s}} \nonumber\\
&\lesssim
\left( \left(  \|\nabla v\|^{s}_{-s,0}\|\nabla v\|_{1,0}\right)^\frac{s}{1+s} \sqrt{\underline{\mathcal{D}}}\right)^{\frac{1}{1+s}}\lesssim\left(\underline{\mathcal{D}}^{1+2s}\underline{\mathfrak{D}}^{s^2}\right) ^{\frac{1}{2(1+s)^2}},\label{dctv4qaa}
\end{align}
where we have used \eqref{xfsdsaf41252} in the second inequality, the incompressible condition in the  third inequality and \eqref{v11} in the fourth inequality.

Making use of  \eqref{im3a}, \eqref{dctv4qaa}, the interpolation inequality \eqref{201807291850}, the  anisotropic  inequality \eqref{4v313}, the definition of $\mathfrak{U}^{i}(\cdot)$ in \eqref{ugg},   and  H\"older's inequality, we  obtain
\begin{align}
&\int_0^{T}(|{\mathcal{L}}_1|+|{\mathcal{L}}_3|	+|{\mathcal{L}}_7|	+|{\mathcal{L}}_{10} |)\mm{d}t	
\nonumber\\ &\lesssim \int_0^{T}\left(\|\varrho\|_{-s,3}+\|\varrho_t\|_{-s,2} \right) \|(v\cdot\nabla\varrho,\nabla\Delta(\underline{\varrho} v_3 -v\cdot\nabla \tilde{\varrho}))\|_{-s,0}  \mm{d}t	
\nonumber\\ &\lesssim \int_0^{T}  \|(v\cdot\nabla\varrho,\nabla\Delta(\left( \bar \rho'\right) ^{-1}v\cdot\nabla\varrho ) \|_{-s,0}\sqrt{\underline{\mathfrak{E}}} \mm{d}t	
\nonumber\\ &\lesssim \int_0^{T}\Big( \|\varrho\|_{1,2}\mathfrak{U}^{1}( v_{\mm{h}})+\|\varrho\|_{3}\mathfrak{U}^{1}(v_3)+\|\varrho\|_{1,3}\mathfrak{U}^0( v_{\mm{h}})+\|\varrho\|_{4}\mathfrak{U}^{0}(v_3)
\nonumber\\ &\quad+\|v_3\|_{2}
\mathfrak{U}^{1}(\partial_{3}\varrho)+\|\nabla v\|_{2}\mathfrak{U}^{0}(\partial_{\mm{h}}\varrho)+
\|v_3\|_{3}\mathfrak{U}^{0}(\partial_{3}
\varrho)\Big)\sqrt{\underline{\mathfrak{E}}} \mm{d}t
\nonumber\\ &\lesssim \int_0^{T}\bigg(\|\nabla v\|_{2}\mathfrak{U}^{0}(\partial_{\mm{h}}\varrho)+\|\varrho\|_{1,3}\mathfrak{U}^{0}( v_{\mm{h}})+\|\varrho\|_{1,2}\mathfrak{U}^{0}( v_{\mm{h}})+\|v_3\|_{2}\mathfrak{U}^{1}(\partial_{3}\varrho)\nonumber\\ &\quad+\|v_3\|_{3}\mathfrak{U}^{0}(\partial_{3}\varrho)+\|\varrho\|_{3}\mathfrak{U}^{1}(v_3)+
\left((\|\varrho\| ^{s}_{1,3}\|\varrho\|_{-s,3})^{\frac{1}{1+s}}\|\varrho\|_{5}\right)^{1/2}\mathfrak{U}^{0}(v_3)\bigg) \sqrt{\underline{\mathfrak{E}}}\mm{d}t
\nonumber
\\	 &\lesssim\int_0^{T}  \sqrt{\underline{\mathfrak{E}}}\bigg(({\underline{\mathfrak{D}}}^{s}\mathcal{D}^{{2+s}})^{\frac{1}{2(1+s)}}
+{\underline{\mathfrak{E}}}^{\frac{2s}{4(1+s)}}
\bigg(\mathfrak{D}^{\frac{4+s}{6(1+s)}}\mathcal{D}^{\frac{1}{3}}+\|\nabla v\|_{1,0}{\mathcal{E}}^{\frac{2}{4(1+s)}}\nonumber\\ &\quad+
(\underline{\mathcal{E}}^{ {2}}+ \underline{\mathcal{E}}{\mathfrak{E}})^{\frac{1}{4(1+s)}}\sqrt{\mathfrak{D}}\bigg)+
(\mathfrak{E}({\underline{\mathfrak{D}}}^{ {s}}{\mathfrak{D}})^{\frac{1}{ (1+s)}})^{\frac{1}{2}}
+(\mathcal{E}({\underline{\mathfrak{E}}}\mathcal{D}^s)^{\frac{1}{(1+s)}})^{\frac{1}{4}}
(\underline{\mathcal{D}}^{ {1+2s}}\underline{\mathfrak{D}}^{ {s^2} })^{\frac{1}{2(1+s)^2}}\bigg) \mm{d}t  \label{a1sq311}
\end{align}	
and
\begin{align}
&\int_0^{T} (|	{\mathcal{L}}_2|+|	{\mathcal{L}} _4|+ |	{\mathcal{L}}_6|+|	{\mathcal{L}} _9|)\mm{d}t
\nonumber\\ \lesssim&\int_0^{T}\Big(\|\varrho\|_{1,2}\big(  \|v_{\mm{h}}\|_{-s,0}+\|\partial_{t}v_{\mm{h}}\|_{-s,0} +\|\varrho\|_{3}(\|v_{3}\|_{-s,0}+\|\partial_{t}v_{3}\|_{-s,0})\big) \mathfrak{U}^{1}(\varrho)\nonumber\\ &+\|\varrho\|_{2}(\|v_{\mm{h}}\|_{-s,0}+\|\partial_{t}v_{\mm{h}}\|_{-s,0})\mathfrak{U}^{0}(\partial_{\mm{h}}\varrho)
+\|v\|_{-s,0}\left(\|\nabla v\|_{0} \mathfrak{U}^{0}(v)+\|v_{t}\|_{0} \mathfrak{U}^{0}(\varrho)\right)
\nonumber\\ &+(  \|\nabla v\|_{-s,2}+\|v_t\|_{-s,1})\big((\|\varrho\|_{1,3}+  \|\nabla v\|_{0}+  \|v_{t}\|_{1})\mathfrak{U}^{2}(\varrho)+\|\nabla v\|_{1} \mathfrak{U}^{0}(v)\big)\Big)\mm{d}t	\nonumber\\ \lesssim&\int_0^{T}\Big(\big(  (\underline{\mathfrak{E}} +\underline{\mathfrak{D}} ){\mathfrak{D}} + {\mathfrak{E}}\underline{\mathfrak{D}}\big) ({\underline{\mathfrak{E}}} ^{ {s} }{\mathfrak{E}})^{\frac{1}{1+s}}+ {\mathfrak{E}}(\underline{\mathfrak{E}} +\underline{\mathfrak{D}} ) ({\underline{\mathfrak{D}}}^{{s}}{\mathfrak{D}})^{\frac{1}{1+s}}  +({\underline{\mathfrak{E}}}^{ {s}}({\mathfrak{E}}+ {\mathcal{E}}))^{\frac{1}{1+s}}\underline{\mathfrak{D}}{\mathcal{D}}\Big)^{\frac{1}{2}}\mm{d}t. \label{a1sq322}
\end{align}	

It is easy to get that	 	
\begin{align}		
&\int_0^{T}  \sqrt{\underline{\mathfrak{E}}}(\mathcal{D}^{2+s}{\underline{\mathfrak{D}}}^{s})^{\frac{1}{2(1+s)}} \mm{d}t \lesssim B^{1/2}(T)\left(\left( \int_0^{T}\mathcal{D}\mm{d}t\right)^{ {2+s} }\left( \int_0^{T}{\underline{\mathfrak{D}}}\mm{d}t\right)^s\right)^{\frac{1}{2(1+s)}}\lesssim	 B^{3/2}(T), \\
&\int_0^{T}  {\underline{\mathfrak{E}}}^{\frac{1+2s}{2(1+s)}}\mathfrak{D}^{\frac{4+s}{6(1+s)}}\mathcal{D}^{\frac{1}{3}}\mm{d}t \nonumber \\
&\lesssim\left(\mathcal{E}^{1+2s}_{\mm{total}} \left( \int_0^{T}\left( (\langle t\rangle^a\mathfrak{D})^{4+s}\mathcal{D}^{2(1+s)}\right)^{\frac{1}{3(2+s)}} \mm{d}t\right)^{2+s} \left( \int_0^{T}\langle t\rangle^{-\frac{a(4+s)}{3s}}\mm{d}t\right)^s\right) ^{\frac{1}{2(1+s)}} \lesssim	 B^{3/2}(T) ,\label{10030932}\\
&\int_0^{T}  {\underline{\mathfrak{E}}}^{\frac{1+2s}{2(1+s)}}{\mathcal{E}}^{\frac{1}{2(1+s)}}\|\nabla v\|_{1,0}\mm{d}t \lesssim B (T)\left( \int_0^{T}\langle t\rangle^{\frac{3a}{2}} {\underline{\mathcal{D}}}\mm{d}t  \int_0^{T}\langle t\rangle^{-\frac{3a}{2}}\mm{d}t\right)^{\frac{1}{2}} \lesssim	 B^{3/2}(T) ,\label{10030933}\\
&\int_0^{T}  {\underline{\mathfrak{E}}}^{\frac{1+2s}{2(1+s)}}\sqrt{\mathfrak{D}}(\underline{\mathcal{E}}^{2}+\underline{\mathcal{E}}{\mathfrak{E}})^{\frac{1}{4(1+s)}}\mm{d}t \lesssim B (T)\left( \int_0^{T}\langle t\rangle^{a} \mathfrak{D}\mm{d}t \int_0^{T}\langle t\rangle^{-\frac{a(9+4s)}{4(1+s)}}\mm{d}t\right)^{\frac{1}{2}}\lesssim	 B^{3/2}(T) ,\label{10030934}\\
&\int_0^{T}   \left({\underline{\mathfrak{E}}}\mathfrak{E}\left(  {\underline{\mathfrak{D}}}^s{\mathfrak{D}}\right)^{\frac{1}{1+s}} \right)^{\frac{1}{2}} \mm{d}t \lesssim B^{3/2}(T)\left( \int_0^{T}\langle t\rangle^{-\frac{a(2+s)}{1+s}}\mm{d}t\right)^{\frac{1}{2}}\lesssim	 B^{3/2}(T), \label{10030935}\\
&	\int_0^{T}  \sqrt{\underline{\mathfrak{E}}}\left( \mathcal{E}\left( {\underline{\mathfrak{E}}}\mathcal{D}^{s}\right)^{\frac{1}{1+s}} \right)^{\frac{1}{4}}	(\underline{\mathcal{D}}^{ {1+2s}}\underline{\mathfrak{D}}^{ {s^2} })^{\frac{1}{2(1+s)^2}}\mm{d}t\nonumber\\
&\lesssim\left( \mathcal{E}^{4+3s}_{\mm{total}}\left( \int_0^{T}\left( \mathcal{D}^{s}\left( {\underline{\mathfrak{D}}}^{2s^2}\left(\langle t\rangle^{\frac{3a}{2}} {\underline{\mathcal{D}}}\right) ^{2+2s}\right)^{\frac{1}{1+s}} \right)^{\frac{1}{2+3s}} \mm{d}t\right)^{2+3s} \left( \int_0^{T}\langle t\rangle^{-\frac{3a(1+2s)}{(1+s)(2+s)}}\mm{d}t\right)^{2 + s}\right) ^{\frac{1}{4(1 + s)}}\nonumber\\
&\lesssim B^{3/2}(T),\label{10030936}
\end{align}
  where we have used \eqref{a1sqa5} in \eqref{10030932}, \eqref{a1sqa1} in \eqref{10030933}, \eqref{a1sqa2} in \eqref{10030934},  \eqref{a1sqa4} in \eqref{10030935}  and \eqref{a1sqa3} in \eqref{10030936}. 

	Putting the above six estimates into \eqref{a1sq311} yields
\begin{align}
&\int_0^{T}(|{\mathcal{L}}_1|+|{\mathcal{L}}_3|	+|{\mathcal{L}}_7|	+|{\mathcal{L}}_{10} |)\mm{d}t\lesssim B^{3/2}(T).	   \label{a1sq31}
\end{align}		
By  \eqref{a1sqa4} and  H\"older's inequality, we get
	\begin{align}
			&	\int_0^{T} \left( \underline{\mathfrak{E}}\mathfrak{D}\left( {\underline{\mathfrak{E}}}^{s}{\mathfrak{E}}\right)^{\frac{1}{(1+s)}} \right)^{\frac{1}{2}}\mm{d}t\lesssim B(T)\left( \int_0^{T}\langle t\rangle^{-\frac{a(2+s)}{1+s}} \mm{d}t \int_0^{T}\langle t\rangle^{a}{\mathfrak{D}} \mm{d}t\right)^{\frac{1}{2}}\lesssim B^{3/2} (T),\nonumber\\
		&	\int_0^{T} \left( \underline{\mathfrak{D}}\mathfrak{D}\left( {\underline{\mathfrak{E}}}^{s}{\mathfrak{E}}\right)^{\frac{1}{(1+s)}} \right)^{\frac{1}{2}}\mm{d}t\lesssim \left(B (T) \int_0^{T}\left(\underline{\mathfrak{D}}\mathfrak{D}\right)^{\frac{1}{2}} \mm{d}t\right)\lesssim B^{3/2} (T),\nonumber\\
			&	\int_0^{T} \left( {\mathfrak{E}}\underline{\mathfrak{D}}\left( {\underline{\mathfrak{E}}}^{s}{\mathfrak{E}}\right)^{\frac{1}{(1+s)}} \right)^{\frac{1}{2}}\mm{d}t\lesssim B(T)\left( \int_0^{T}\langle t\rangle^{-\frac{a(2+s)}{1+s}} \mm{d}t \int_0^{T}\underline{\mathfrak{D}} \mm{d}t\right)^{\frac{1}{2}}\lesssim B^{3/2} (T),\nonumber\\
		&	\int_0^{T}  \left({\mathfrak{E}}\underline{\mathfrak{E}} \left( {\underline{\mathfrak{D}}}^{s}{\mathfrak{D}}\right) ^{\frac{1}{(1+s)}}  \right)^{\frac{1}{2}}\mm{d}t \lesssim B(T)\left( \int_0^{T}\langle t\rangle^{-\frac{a(2+s)}{1+s}} \mm{d}t \int_0^{T}\left( {\underline{\mathfrak{D}}}^{s}\langle t\rangle^{a}{\mathfrak{D}}\right) ^{\frac{1}{(1+s)}} \mm{d}t\right)^{\frac{1}{2}}\lesssim B^{3/2} (T),\nonumber\\
		&	\int_0^{T}  \left({\mathfrak{E}}\underline{\mathfrak{D}} \left( {\underline{\mathfrak{D}}}^{s}{\mathfrak{D}}\right) ^{\frac{1}{(1+s)}}  \right)^{\frac{1}{2}}\mm{d}t \lesssim \left(B (T) \int_0^{T}\underline{\mathfrak{D}} \mm{d}t  \int_0^{T}\left( {\underline{\mathfrak{D}}}^{s}{\mathfrak{D}}\right) ^{\frac{1}{(1+s)}} \mm{d}t\right)^{\frac{1}{2}} \lesssim B^{3/2} (T), \nonumber\\
		&	\int_0^{T}  \left( \underline{\mathfrak{D}}{\mathcal{D}}{\underline{\mathfrak{E}}}^{s} {\mathcal{E}}  ^{\frac{1}{1+s}} \right)^{\frac{1}{2}}\mm{d}t \lesssim B^{1/2}(T)  \int_0^{T}\left( \underline{\mathfrak{D}}{\mathcal{D}}\right)^{\frac{1}{2}} \mm{d}t\lesssim B^{3/2} (T). \nonumber
	\end{align}		
Inserting the above six estimates into \eqref{a1sq322} yields
		\begin{align}
&\int_0^{T} (|{\mathcal{L}}_2|+|{\mathcal{L}} _4|+ |	{\mathcal{L}} _6|+|{\mathcal{L}} _9|)\mm{d}t \lesssim B^{3/2} (T) 	. \label{a1sq32}
\end{align}

In addition, utilizing the boundary  conditions  of $(\bar{\rho}'',\omega_{\mm{h}},\partial_3^{2}\omega_{\mm{h}},\partial_3^{2}\varrho)$ in \eqref{0102n1} and \eqref{omega} and the integration by parts,   we get
\begin{align}
{\mathcal{L}}_{j}\lesssim&	\begin{cases}
(\|\varrho\|_{1-s,1}+\|\nabla v_t\|_{-s,0})\|\nabla v\|_{-s,2}&\mbox{for } j=5;
\\
(\|\varrho\|_{1-s,2}+\|v_t\|_{-s,0})\|\omega_t\|_{-s,0}&\mbox{for } j=8.
\end{cases}	\label{20223402201565ssa}
\end{align}
Exploiting \eqref{a1sq31}--\eqref{20223402201565ssa} and Young's inequality, we infer \eqref{dctv1}--\eqref{dctv4}  from \eqref{dctv1q}--\eqref{dctv4q} for the case that ${\ell_1}=0$ and ${\ell_2}=-s$.

(2)  {Now we turn to consider the case  ${\ell_1}=a\,\,\hbox{or}\,\,0$ and ${\ell_2}=0$}. Let
$$\tilde{\mathbf{Q}}^{N}_3:=  \partial_1\big(  {\partial_1 \varrho   }\left(\partial_3 \varrho - { \bar \rho' \varrho}/{\bar \rho } \right)/{ \varrho + \bar{\rho} }\big)+  \partial_2\big(  {\partial_2 \varrho   }\left(\partial_3 \varrho - { \bar \rho' \varrho}/{\bar \rho } \right)/({ \varrho + \bar{\rho} })\big).$$
Obviously it holds that
\begin{align*}
\mathbf{Q}^{\mm N}_3 =&  \partial_{3}\left( \left(  {( \bar{\rho}'\varrho/ \bar{\rho} )^2} +    {(\partial_3 \varrho)^2   } -{ 2\bar \rho' \varrho \partial_{3}\varrho}/{\bar \rho }\right)/({\varrho + \bar{\rho}})\right)
+{\tilde{\mathbf{Q}}^{\mm N}_3}.
\end{align*}

Exploiting the  boundary condition of $v_3$ in \eqref{n1},  the above identity,  the incompressible condition and the integration by parts, we rewrite ${\mathcal{L}}_4$ as follows
\begin{align*}
{\mathcal{L}}_4:=&-\int\langle t\rangle^{\ell_1}\bigg(\left( \varrho v_t+{\rho} v\cdot  \nabla v\right)  \cdot v_t+ {\varepsilon ^2}\left( \mathbf{Q}^{\mm N}_{\mm{h}}  \cdot\partial_{t} v_{\mm{h}}+  {\mathbf{Q}}^{\mm N}_3 \partial_{t} v_{3}\right)  \bigg)  \mm{d}x\nonumber\\=&-\int\langle t\rangle^{\ell_1}\bigg(\left( \varrho v_t+{\rho} v\cdot  \nabla v\right)  \cdot v_t+ {\varepsilon ^2}\left({( \mathbf{Q}^{\mm N}_{\mm{h}}+\tilde{\mathbf{Q}}^{\mm N}_{\mm{h}} ) } \cdot\partial_{t} v_{\mm{h}}+ \tilde{\mathbf{Q}}^{\mm N}_3 \partial_{t} v_{3}\right)  \bigg)\mm{d}x
{:=\tilde{{\mathcal{L}}}_4},
\end{align*}
where	we have defined that
\begin{align*}
\tilde{\mathbf{Q}}^{\mm N}_{\mm{h}}:=& \nabla_{\mm{h}}\left(  \left(   2 \bar \rho' \varrho\partial_{3}\varrho / \bar \rho -( \bar{\rho}'\varrho /  \bar{\rho})^2 - ( \partial_3 \varrho)^2  \right)/ (\varrho  + \bar \rho )\right) .	
\end{align*}		Consequently, by virtue of \eqref{imvarrho},
the product estimates in \eqref{fgestims} and  the bounds of  density  in \eqref{im3a}, we arrive at		
\begin{align}
&	\int_0^{T}(| {\mathcal{L}}_3|+|
{\tilde{{\mathcal{L}}}_{4}}|+|{\mathcal{L}}_9|+|{\mathcal{L}}_{10} |)\mm{d}t\nonumber	\\&\lesssim 	\int_0^{T}\langle t\rangle^{\ell_1}\Big(\left( \|\varrho\|_{2}\|v\|_{2}\|\nabla v\|_{0}+\|\varrho\|_{1,2}\|\varrho\|_{3}\right) \|v_t\|_{0}+(\|\varrho\|_{1,2}\|v_{\mm{h}}\|_{2}\nonumber\\&\quad+\|\varrho\|_{3}\|v_{3}\|_{3})\|\varrho_t\|_{2}+\left( \|\varrho\|_{3}(\|\varrho\|_{1,2}+\|v_t\|_{1})+\|v\|_{2}\|\nabla v\|_{2}\right) \|v_t\|_{1}\Big)\mm{d}t	\nonumber \\ &\lesssim\int_0^{T}	\langle t\rangle^{\ell_1}\sqrt{\mathfrak{E}}\mathfrak{D}\mm{d}t\lesssim	 B^{3/2}(T). \label{i610}
\end{align}
In addition, we can follow the argument of \eqref{20223402201565ssa} to obtain
\begin{align}
{\mathcal{L}}_{8}\lesssim&\langle t\rangle^{\ell_1}(\|\varrho\|_{1,2}+\|v_t\|_{0})\|\omega_t\|_{0}. \label{i6101}
\end{align}	
		
Using \eqref{i610}, \eqref{i6101}  and  Young's inequality, we infer \eqref{dctv3} and \eqref{dctv4}	from \eqref{dctv3q} and \eqref{dctv4q}   for the case that ${\ell_1}=a$ or $0$,  and ${\ell_2}=0$.
		
(3) Finally, arguing analogously to  \eqref{i610}, we  get
\begin{align*}
&\int_0^{T}( |{\mathcal{L}}_1|+ |{\mathcal{L}}_2|)\mm{d}t	
\nonumber\\ &\lesssim\int_0^{T}\langle t\rangle^{3a/2}\big( 	\|\varrho\|_{1,2}\left( (\|\varrho\|_{3}+\|v_t\|_{0})\|\nabla v\|_2 +\|\varrho\|_{1,2}\|v\|_{2}\right) +\|\varrho\|_2\|\nabla v\|_0\| v_t\|_1\\ &\quad  +\left(1+ \|\varrho\|_{3}\right)\|v\|_2\|\nabla v\|^2_2\big)\mm{d}t  \lesssim\int_0^{T}	\langle t\rangle^{3a/2}\sqrt{\mathfrak{E}}\mathfrak{D}\mm{d}t\lesssim	 B^{3/2}(T).
\end{align*}			
Putting the above estimate into \eqref{dctv1q} yields \eqref{dctv1} with $({\ell_1},{\ell_2})=(3a/2,1)$. The proof is completed. \hfill $\Box$
\end{pf}	
\subsection{Dissipative estimates of $\varrho$}
This section is devoted to the derivation of  the dissipative estimates of the tangential derivatives to $\varrho$.
\begin{lem}
\label{lem2dxfa1}
For any given $t\in\overline{I_T}$, we have
\begin{align}
& \|\varrho\|_{1+\ell,1}\lesssim 	\begin{cases}
\|\nabla v\|_{-s,1}+\| v_t\|_{-s,0}+\left(({\underline{\mathfrak{E}}}^{{s}}{{\mathfrak{E}}})^{\frac{1}{(1+s)}}\mathfrak{D}\right)^{1/2}&\mbox{for } \ell=-s;
\\
\|\nabla v\|^2_{1}+\|v_t\|^2_0 + \left(\sqrt{{\mathfrak{E}}}\mathcal{D}\right)^{1/2}&\mbox{for } \ell=0,
\end{cases}	 \label{Pfsda2}
\end{align}	
and \begin{align}
& \|\varrho\|_{1+\ell,j}\lesssim 	\begin{cases}
\|\varrho\|_{ 1-s,1}+\|\nabla v\|_{-s,2}+\| v_t\|_{-s,1} +
({\underline{\mathfrak{E}}}^{{s}}{{\mathfrak{E}}})^{\frac{1}{2(1+s)}}
 \underline{\mathfrak{D}} ^{\frac{1}{2}}&\mbox{for } \ell=-s,\ j=2	;\\
\|\varrho\|_{ 1,1}+\|\nabla v\|_{2}+\|v_t\|_{1}+ (\sqrt{{\mathfrak{E}}}\mathfrak{D})^{\frac{1}{2}}&\mbox{for } \ell=0,\ j=2	;\\
\|\varrho\|_{ 1,1}+\|\nabla v\|_{5}+\|v_t\|_4 + (\sqrt{{\mathcal{E}}}\mathcal{D})^{\frac{1}{2}}&\mbox{for } \ell=0,\ j=5.\label{r4}
\end{cases}
\end{align}		
\end{lem}	
\begin{pf}	
(1)	Taking the inner product of  \eqref{1aaaB}$_2$ with $\ell$ in place of $\ell_2$ and $\Lambda^{\ell }_{\mm{h}}\nabla_{\mm{h}}^{\perp}\varrho/{\bar{\rho}'}$ in $L^2$, and then using  the boundary conditions of $\varrho$ in \eqref{varrho} and the integration by parts, we get
\begin{align}
&E( \Lambda^{\ell }_{\mm{h}}\partial_1\varrho)+E( \Lambda^{\ell }_{\mm{h}}\partial_2\varrho)
={\mathcal{V}},\label{LEM2P1}
\end{align}
where we have defined that
\begin{align}
& {\mathcal{V}}:={\mathcal{V}}^{\mm{L}}+ {\mathcal{V}}^{\mm{N}},\nonumber \\
&{\mathcal{V}}^{\mm{L}}:= \int \frac{1}{\bar{\rho}'}\Bigg( \Lambda^{\ell }_{\mm{h}}\Bigg(\frac{\mu\bar{\rho}''} {\bar{\rho}'}\partial_{3}\omega_{\mm{h}} 	-{ \bar{\rho}}\partial_{t}\omega_{\mm{h}}-\mathbf{M} _{\mm{h}}  \Bigg)
\cdot \Lambda^{\ell }_{\mm{h}}\nabla_{\mm{h}}^{\perp}\varrho
-\mu {\Lambda^{\ell }_{\mm{h}}\partial_i \omega_{\mm{h}} }\cdot \partial_i \Lambda^{\ell }_{\mm{h}}
\nabla_{\mm{h}}^{\perp}\varrho
\Bigg)\mm{d}x,\label{2025081313071}\\
& {\mathcal{V}}^{\mm{N}}:=-\int \frac{1}{\bar{\rho}'}\Lambda^{\ell }_{\mm{h}}\Bigg(  	{ {\varrho}}\partial_{t}\omega_{\mm{h}}
+ \mathbf{N}_{\mm{h}}  +{\rho}  v\cdot  \nabla \omega_{\mm{h}}\Bigg)
\cdot \Lambda^{\ell }_{\mm{h}}\nabla_{\mm{h}}^{\perp}\varrho
 \mm{d}x.\label{202508131307}
\end{align}		
Exploiting \eqref{imvarrho}, \eqref{im3a}, the product estimates in  \eqref{fgestims}, and the anisotropic  inequality \eqref{4v313}, the boundary condition of $\varrho$ and the integration by parts, we  obtain
\begin{align}
{\mathcal{V}}\lesssim&	\begin{cases}
\|\varrho\|_{1-s,1} \Big(\|\nabla v\|_{-s,1}+\| v_t\|_{-s,0}&\\+ (\|\varrho\|_{1,2}+\|\nabla v\|_1+\|v_t\|_0 ) \left( \mathfrak{U}^{2}(\varrho)+\mathfrak{U}^{0}(v)\right)\Big)&\mbox{for } {\ell }=-s;
\\
\|\varrho\|_{1,1} (\|\nabla v\|_{1}+\|v_t\|_0 +\|\varrho\|_{1,2}\|\varrho\|_3 + \|v\|_2\|\nabla v\|_2 )&\mbox{for } {\ell }=0,
\end{cases}\label{LEM2P1a}
\end{align}
where the definitions of $\mathfrak{U}^{i}(\cdot)$ can be found in \eqref{ugg} for $i=0$ and $1$.

In addition, it is easy to deduce that
\begin{align}
\mathfrak{U}^{2}(\varrho)+\mathfrak{U}^1(v)\lesssim
({\underline{\mathfrak{E}}}^{{s}}{{\mathfrak{E}}})^{\frac{1}{2(1+s)}}. \label{u2p}
\end{align}
Consequently, using \eqref{LEM2P1a}, \eqref{u2p},   Young's inequality and the stabilizing estimate \eqref{F}, we obtain \eqref{Pfsda2} from \eqref{LEM2P1}.	
		
(2)	Let $\eta = \Lambda^{\ell }_{\mm{h}}\nabla^{\perp}_{\mm{h}}\varrho$. In view of the equation \eqref{1aaaB} with $\ell$
in place of $\ell_2$ and the boundary conditions of $\varrho$ in \eqref{varrho},
$\eta$  satisfies the following boundary value problem
\begin{equation}
\label{Stokesequson}
\begin{cases}
\Delta\eta= {\bar \rho }   \Lambda^{\ell }_{\mm{h}}({\rho}\partial _t \omega_{\mm{h}} -
g\eta_{\mm{h}}-\mu\Delta \omega_{\mm{h}}-\tilde{\mathbf{L}}_{\mm{h}} +{\mathbf{M}}_{\mm{h}}+  {\rho} v\cdot  \nabla \omega_{\mm{h}} + \mathbf{N}_{\mm{h}})/{\varepsilon ^2\bar \rho '} ,  \\
\eta  |_{\partial\Omega} =0,
\end{cases}
\end{equation}
where we have defined that $\tilde{\mathbf{L}}_{\mm{h}}:=\mathbf{L}_{\mm{h}}-\epsilon^2\bar{\rho}'\bar{\rho}^{-1} \nabla^{\perp}_{\mm{h}}\Delta\varrho $.
Making use of	\eqref{imvarrho}, \eqref{im3a}, the product estimates in \eqref{fgestims}, the  elliptic estimate \eqref{xfsdsaf41252}, and the anisotropic  inequality \eqref{4v3},   we get
\begin{align*}
\|\varrho\|_{{1+{\ell }},j}=\|\eta\|_{j}
\lesssim& 	\begin{cases}
\|\varrho\|_{ 1-s,1}+\|\nabla v\|_{-s,2}+\| v_t\|_{-s,1}&\\+ (\|\varrho\|_{1,2} +\|\nabla v\|_1+\|v_t\|_1)\left( \mathfrak{U}^{2}(\varrho)+\mathfrak{U}^{1}(v)\right)&\mbox{for } {\ell }=-s,\ j=2	;\\
\|\varrho\|_{ 1,1}+\|\nabla v\|_{2}+\|v_t\|_1 &\\+\left(\|\varrho\|_{ 3}+\|v\|_{2}\right) (
\|\varrho\|_{1,2}+\|\nabla v\|_1+\|v_t\|_1 )&\mbox{for } {\ell }=0,\ j=2	;\\
\|\varrho\|_{ 1,4}+\|\nabla v\|_{5}+\|v_t\|_4 + \|\varrho\|_{1,5}\|\varrho\|_{ 6} &\mbox{for } {\ell }=0,\ j=5,
\end{cases}
\end{align*}
which, together with  \eqref{u2p}, yields \eqref{r4}. The proof is completed.
\end{pf}

\subsection{Estimates of temporal derivatives}
We next derive some energy/dissipative estimates for the temporal derivatives of $(\varrho,v)$.
\begin{lem}
\label{lem2dxfasfa1}
For any given $t\in\overline{I_T}$, we have
\begin{align}
&\|\varrho_t\|_{\ell,j} \lesssim
\begin{cases}
\|v_3\|_{-s,j}+	\left(({\underline{\mathfrak{E}}}^s{{\mathcal{E}}})^{\frac{1}{1+s}}	{{\mathfrak{D}}} +{{\mathfrak{E}}}({\underline{\mathfrak{D}}}^{{s}}{{\mathfrak{D}}})^{\frac{1}{1+s}}\right)^{\frac{1}{2}}	&\mbox{for }\ell=-s,\ 1\leqslant j\leqslant2;\\	 \|v_3\|_{j} +\|\varrho\|_{1,j}\|v\|_{2} &\mbox{for }\ell=0,\ 0\leqslant j\leqslant1;\\
 \|v_3\|_{j} +\|\varrho\|_{1,j}\|v\|_{j} &\mbox{for }\ell=0,\ 2\leqslant j\leqslant5,
\end{cases}\label{202221401177321} \\
&\| v_{t}\|_{\ell,0} \lesssim\begin{cases}
\|(\varrho,v)\|_{-s,2} +(\|\varrho\|_3+\|v_t\|_{0}+\|v\|_2)({\underline{\mathfrak{E}}}^{{s}}{{\mathcal{E}}})^{\frac{1}{2(1+s)}}
&\mbox{for }\ell=-s;
\\\|(\varrho,v)\|_2 +\|\varrho\|_3(\|\varrho\|_3+\|v_t\|_{0})  &\mbox{for }\ell=0 ,\end{cases}\label{vt0} \\
&	\| v_t\|_3 \lesssim  \|(\varrho,v)\|_5+\|\varrho\|_3(\|\varrho\|_3+\|v_t\|_{0}) . \label{20222140117732}
\end{align}		
\end{lem}	
\begin{pf}
(1) Using the bounds of density in \eqref{im3a}, \eqref{im3}, the product estimates in \eqref{fgestims} and the product estimate \eqref{4v3}, we derive from \eqref{h1f}$_1$ that
\begin{align}
\|\varrho_t\|_{\ell,j} &\lesssim
\begin{cases}
\|v_3\|_{-s,j}+	\|\varrho\|_{3}\mathfrak{U}^1(v_3) +	 \|v_3\|_{2}\mathfrak{U}^{1}(\varrho) &\\+\|\varrho\|_{1,2}\mathfrak{U}^{1}(v_\mm{h})+\|v_\mm{h}\|_{2}\mathfrak{U}^0(\nabla_\mm{h}\varrho)&
\mbox{for }\ell=-s,\ 1\leqslant j\leqslant2;\\
(1+\|\varrho\|_{3})\|v_3\|_{j} +\|\varrho\|_{1,j}\|v\|_{2} &\mbox{for }\ell=0,\ 0\leqslant j\leqslant1;\\
(1+\|\varrho\|_{j+1})\|v_3\|_{j} +\|\varrho\|_{1,j}\|v\|_{j} &\mbox{for }\ell=0,\ 2\leqslant j\leqslant5,
\end{cases}
\label{2508181314}
\end{align}
which, together with \eqref{imvarrho} and the definition of $\mathfrak{U}^{i}(\cdot)$ in \eqref{ugg}, yields \eqref{202221401177321}.		
		
(2)	Taking the inner product of the momentum equations \eqref{h1f}$_2$  and $\Lambda^{\ell_2}_{\mm{h}}v_t$ in $L^2$, and then using  \eqref{imvarrho}, \eqref{im3a}, the product estimates in \eqref{fgestims}, and the anisotropic inequality \eqref{4v3},  we obtain
\begin{align*} \|v_{t}\|_{{\ell},0}^2&\lesssim   \|\sqrt{\bar\rho}\Lambda^{\ell}_{\mm{h}}v_{t}\|_0^2\\&=
\int \Lambda^{\ell_2}_{\mm{h}} (\mu\Delta v-g{\varrho}\mathbf{e}^3- {\varepsilon ^2}\mathbf{Q}-{\varrho} v_{t}-   \rho v\cdot  \nabla v)\cdot \Lambda^{\ell}_{\mm{h}}v_t\mm{d}x
\\&	\lesssim\begin{cases}
\|(\varrho,v)\|_{-s,2} +(\|\varrho\|_2+\|v_t\|_{0})\mathfrak{U}^{1}(\varrho)+
\|v\|_1\mathfrak{U}^0(v)\|v_t\|_{-s,0} &\mbox{for }{\ell}=-s;
\\
\big(\|(\varrho,v)\|_2 +\|\varrho\|_2(\|\varrho\|_3+\|v_t\|_{0})+\|v\|_2 \big)\|v_t\|_0 &\mbox{for }{\ell}=0 ,
\end{cases}
\end{align*}
which, together with  the definition of $\mathfrak{U}^{i}(\cdot)$ in \eqref{ugg},  implies  \eqref{vt0}.
		
(3) Applying $\|\cdot\|_2$ to the vortex equation  \eqref{l41a}, and then using \eqref{imvarrho},   the lower-bound of density in \eqref{im1},  \eqref{im3a} and the product estimates in \eqref{fgestims},  we  obtain
\begin{align}
\| \omega_t\|_2=\|\left( \mu\Delta \omega+ {\mathbf{L}}- {\mathbf{M}} -   {\rho} v\cdot  \nabla \omega - \mathbf{N}\right)/ {\rho}\|_2\lesssim \|(\varrho,v)\|_5+\|v_t\|_2. 	\label{202240131012024a}
\end{align}
Thus,  exploiting \eqref{im3}, \eqref{vt0} with $\ell=0$, the Hodge-type elliptic estimate \eqref{202005021302} and the incompressible condition \eqref{1a}$_3$,
we deduce from the above estimate that
\begin{align}
\label{202240131012024}
\|v_t\|_3\lesssim \| v_t\|_{0}+\| \nabla v_t\|_2 \lesssim  \|(\varrho,v)\|_5+\|\varrho\|_3(\|\varrho\|_3+\|v_t\|_{0}),
\end{align}
which  yields \eqref{20222140117732}. This completes the proof. \hfill $\Box$
\end{pf}

\subsection{Energy inequalities}

With Lemmas \ref{varrho1}--\ref{lem2dxfa1} in hand, we are in the position to the derivation of energy inequalities for the QRT problem.
\begin{pro}\label{lem3s}
For any given $t\in [0,T]$, we have
\begin{enumerate}
\item[(1)]  the high-order energy inequality
\begin{align}
\label{omessetsimQ}
&\sup_{0\leqslant  t\leqslant  T}\{\mathcal{E} (t)\}+\int_0^{T} \mathcal{D} (t) \mm{d}t \lesssim \mathcal{E}(0)+\sup_{0\leqslant  t\leqslant  T}\{\mathcal{E}^{3/2} (t)\}+\int_0^{T}\Big(\sqrt{\mathcal{E}} \mathcal{D}+ \mathfrak{N}(t)\Big)   \mm{d}t,
\end{align}
where the definition of $\mathfrak{N}$ can be found in \eqref{msaft}.
 \item[(2)] the energy inequality of the  tangential  derivatives
\begin{align}
&\sup_{0\leqslant  t\leqslant  T}\left\lbrace \langle t\rangle^{3a/2}\underline{\mathcal{E}}(t)\right\rbrace +\int_0^{T}\langle t\rangle^{3a/2}\underline{\mathcal{D}}(t) \mm{d}t
\nonumber \\
&	\lesssim \|\varrho^0\|_{{1},1}^2+  \|v^0\|_{{1},0}^2+ B^{3/2}(T)+\int_0^{T}\langle t\rangle^{a/2}\left( \|\varrho\|^2_{1,1}+\|v\|^2_{1,0}\right)  \mm{d}t.\label{qec1}
\end{align}	
\item[(3)]the low-order energy inequality
 \begin{align}
&\sup_{0\leqslant  t\leqslant  T}\left\lbrace \langle t\rangle^{a}{\mathfrak{E}}(t)\right\rbrace +\int_0^{T} \langle t\rangle^{a}{\mathfrak{D}}(t)  \mm{d}t\lesssim\mathfrak{E}(0)+B^{3/2}(T)+\int_0^{T} \langle t\rangle^{a-1}\left(\|\varrho\|^2_{3}+\|v\|^2_{2}\right)  \mm{d}t,	\label{qecc1}
\end{align}	
where the definition of $B(T)$ has been provided by \eqref{202508172246}.
 \item[(4)]the energy inequality  of the negative derivatives
\begin{align}
&\sup_{0\leqslant  t\leqslant  T}\{\underline{\mathfrak{E}}(t)\}+\int_0^{T} \underline{\mathfrak{D}}(t) \mm{d}t\lesssim \underline{\mathfrak{E}}(0)+	 B^{3/2}(T). \label{lem2dx1}
\end{align}		
\end{enumerate}	
\end{pro}
\begin{pf}
(1) Integrating \eqref{lemiq11x}, \eqref{202221401177321a} with $(k,j)=(0,2)$ and \eqref{2022sfa2140117732} over $(0,T)$, and then using the bounds of density in \eqref{im1}, the stabilizing estimate \eqref{F} and the product estimate in \eqref{fgestims},  we get
\begin{align}
&\sup_{0\leqslant  t\leqslant  T}\left\lbrace \|\varrho\|^2_{1} +\|  v\|^2_{0} \right\rbrace  + \int_0^{T}\|\nabla v \|^2_{ 0}\mm{d}t\lesssim \|\varrho^0\|^2_{1} +\|  v^0\|^2_{0}+\int_0^{T}\sqrt{\mathfrak{E}}\mathfrak{D}\mm{d}t,
\label{lemiq11ax}\\
& \sup_{0\leqslant  t\leqslant  T}\left\lbrace \| (\Delta^{3}\tilde{\varrho},\Delta^{2}\omega)\|^2_{0}
\right\rbrace + \int_0^{T}\|\nabla \Delta^{2}\omega\|^2_{0}\mm{d}t\nonumber \\
& \lesssim \| (\varrho^0,v^0)\|^2_3+\sup_{0\leqslant t\leqslant T}\{\|\varrho\|_3\|\varrho\|_6^2\}+\int_0^{T}\left( \sqrt{\mathcal{E}}\mathcal{D}+\mathfrak{N}(t)+(\|\varrho\|_{1,4}+ \| v_t\|_4)\|\nabla v\|_5\right) \mm{d}t
\label{202221401177321aa}
\end{align}	
and
\begin{align}
&\sup_{0\leqslant  t\leqslant  T}\left\lbrace
2\varepsilon ^2\int\Delta^2\left(  \frac{\bar{\rho}'}{\bar{\rho}}\Delta{\varrho}-\left( \frac{\bar{\rho}'}{\bar{\rho}}\right)'\partial_{3}{\varrho}\right) \Delta^2  v_3\mm{d}x +\mu\|\nabla \Delta^2 v\|_0^2\right\rbrace +c\int_0^{T}\| \Delta^2  v_t\|_{0}^2\mm{d}t\nonumber\\&\lesssim
\| (\nabla \varrho^0,v^0)\|^2_5
+ \int_0^{T}(\|\varrho\|_{1,4}^2+\|  v_3\|_5^2+\|  v_t\|_3^2+ \sqrt{\mathcal{E}} \mathcal{D})\mm{d}t.
\label{202221401177321aas}
\end{align}	
	It follows from   \eqref{1ee3} with $({\ell_1},{\ell_2})=(0,0)$,  \eqref{1ee4} with $({\ell_1},{\ell_2})=(0,0)$, \eqref{i610} and \eqref{i6101} that
\begin{align}
& \sup_{0\leqslant  t\leqslant  T} \left\lbrace \frac{\mu}{2}\|\nabla v\|^2_{0}-\int\left( g{\varrho}v_3+{\varepsilon ^2}\mathbf{Q}^{\mm L}\cdot v\right)  \mm{d}x\right\rbrace +c\int_{0}^{T}\|v_t\|^2_{0}\mm{d}t \nonumber\\&\lesssim \|(\varrho^0,v^0)\|_1^2+\int_0^{T}( \|\varrho_{t}\|^2_{1}+\sqrt{\mathfrak{E}}\mathfrak{D}) \mm{d}t.\label{dctv3aaaa}
\end{align}	

Now, making use of the boundary conditions of $v_3$ in \eqref{vbian}, \eqref{im3},  the  stabilizing estimate \eqref{F}, the dissipative estimate of $\nabla_{\mm{h}}\varrho$ in \eqref{Pfsda2} with $\ell=0$,
the estimate of  $\varrho_t$ in \eqref{202221401177321} with $(\ell,j)=(0,1)$,  the integration by parts and Young's inequality, we derive from \eqref{lemiq11ax} and \eqref{dctv3aaaa} that
\begin{align}
&\sup_{0\leqslant  t\leqslant  T}\big\{\tilde{\mathcal{E}}_1(t)\big\}+ \int_0^{T}\tilde{\mathcal{D}}_1(t)\mm{d}t  \lesssim  \|(\varrho^0,v^0)\|_1^2	+\int_0^{T}\sqrt{\mathcal{E}} \mathcal{D}\mm{d}t, \label{lemisds}
\end{align}	
where we have defined that
\begin{align*}
\tilde{\mathcal{E}}_1(t):=\|(\varrho, v)\|^2_{1}\mbox{ and }\tilde{\mathcal{D}}_1(t):=\|\varrho\|^2_{1,1}+\|(\nabla v,v_t)\|^2_{0}.
\end{align*}	
Utilizing \eqref{nav1}--\eqref{nav1var}, \eqref{im3}, \eqref{Pfsda2} with $\ell=0$,  \eqref{r4} with $(\ell,j)=(0,5)$,   \eqref{202221401177321}, the interpolation inequality \eqref{201807291850}  and  Young's inequality, we deduce from  \eqref{202221401177321aa},  \eqref{202221401177321aas} and \eqref{lemisds} that
\begin{align}
&\sup_{0\leqslant  t\leqslant  T}\big\{\tilde{\mathcal{E}}_2(t)\big\}+ \int_0^{T}\tilde{\mathcal{D}}_2(t)\mm{d}t\nonumber  \\
& \lesssim \|(\varrho^0,v^0)\|_5^2+ \sup_{0\leqslant  t\leqslant  T}\left\lbrace  {\|\varrho\|_1\|\varrho\|_6^2}+\| \varrho\|^2_{5} + \|  v\|^2_{4}\right\rbrace  	\nonumber\\&\quad +\int_0^{T}\left(\|\varrho\|^2_{1,4}+\|  \nabla v\|^2_{4}+\|  v_t\|^2_{3}+ \sqrt{\mathcal{E}} \mathcal{D}+\mathfrak{N}(t)\right) \mm{d}t, \label{lemisdf}
\end{align}	
where we have defined that
\begin{align*}
\tilde{\mathcal{E}}_2(t):= \|\varrho\|^2_{6} +\|  v\|^2_{5}\mbox{ and }\tilde{\mathcal{D}}_2(t):=\|\varrho\|^2_{1,5}+\|\nabla v\|^2_{5}+\|v_t\|^2_{4}.
\end{align*}
		
In addition, using  \eqref{imvarrho}, the estimate of $\varrho_t$ in  \eqref{202221401177321} with $(\ell,j)=(0,5)$, the estimate of $v_t$ in  \eqref{20222140117732},
we easily obtain that, for  sufficiently large  constant $\chi$,
\begin{align}
	&	{\mathcal{E}}(t)\lesssim\tilde{\mathcal{E}}_2(t)+\|\varrho_t\|^2_{5} +\|  v_t\|^2_{3} \lesssim \tilde{\mathcal{E}}_2(t)+\mathcal{E}^{3/2} (t)  \label{202402081929}
\end{align}
and
\begin{align}
{\mathcal{D}}(t)\lesssim\tilde{\mathcal{D}}_2(t)+\|\varrho_t\|^2_{5} \lesssim\tilde{\mathcal{D}}_2(t) +\sqrt{\mathcal{E}} \mathcal{D}. \label{202402081929a}
\end{align}
Finally, making use of \eqref{lemisds}, \eqref{202402081929}, \eqref{202402081929a}, the  interpolation inequality \eqref{201807291850} and Young's inequality, we deduce \eqref{omessetsimQ} from \eqref{lemisdf}.

(2) The estimate \eqref{qec1} immediately follows
from the tangential  estimates of $\varrho$ and $v$ in \eqref{dctv1} with $({\ell_1},{\ell_2})=(3a/2,1)$.

(3) Taking the inner products of \eqref{lemiq11x} resp. \eqref{202221401177321a} with $(k,j)=(1,0)$ and
$\langle t\rangle^{a}$  resp. $\langle t\rangle^{a}$ in $L^2(0,T)$, and then exploiting the bounds of density in \eqref{im1}, the stabilizing estimate \eqref{F} and the product estimates in \eqref{fgestims},  we have
\begin{align}
&\sup_{0\leqslant  t\leqslant  T}\left\lbrace  \langle t\rangle^{a}\left(\|\varrho\|^2_{1} +\|  v\|^2_{0} \right)\right\rbrace   + \int_0^{T}\langle t\rangle^{a}\|\nabla v \|^2_{ 0}\mm{d}t\nonumber \\
&\lesssim \|\varrho^0\|^2_{1} +\|  v^0\|^2_{0}+ B^{3/2}(T)+\int_0^{T}\langle t\rangle^{a-1}\left(\|\varrho\|^2_{1} +\|  v\|^2_{0} \right)\mm{d}t,
\label{lemiq11axc}
\end{align}
and
\begin{align}
& \sup_{0\leqslant  t\leqslant  T}\left\lbrace \langle t\rangle^{a}\left(\|\nabla \Delta\tilde{\varrho}\|^2_{0} +\|\nabla\omega\|^2_{0}
\right)\right\rbrace + \int_0^{T}\langle t\rangle^{a}\|\Delta\omega\|^2_{0}\mm{d}t\nonumber \\
&\lesssim  \|(\nabla {\varrho}^0,v^0)\|^2_{0}+ B^{3/2}(T)	+\int_0^{T}  \langle t\rangle^{a}\big((\|\varrho\|_{1,1}+\|v_t\|_{1})\|\nabla v\|_2\nonumber \\
&\quad+\mathfrak{T}(t)+\langle t\rangle^{-1}(\|\varrho\|^2_{3}+\|\nabla v\|^2_{1})\big)\mm{d}t.
\label{202221401177321aac}
\end{align}	

Recalling the definition of $\mathfrak{T}$ in \eqref{202508091522}, and then using   \eqref{imvarrho}, we have	
\begin{align}
\int_0^{T}\langle t\rangle^{a}\mathfrak{T}(t)\mm{d}t	&\lesssim B^{\frac{3}{2}}(T)\bigg(  \left( \int_0^{T}\langle t\rangle^{-\frac{3a}{2}}\mm{d}t\right)^{\frac{1}{2}}+\left( \int_0^{T}\langle t\rangle^{-2a}\mm{d}t\right)^{\frac{1}{2}}\bigg) \nonumber\\&\quad+ \int_0^{T} \langle t\rangle^{-\frac{5a}{8}} \left(\left( \langle t\rangle^{a}\mathfrak{E}\right) ^{{3}} \langle t\rangle^{\frac{3a}{2}} \underline{\mathcal{D}}\right) ^{\frac{1}{4}}\left(\langle t\rangle^{a} \mathfrak{D}\right) ^{\frac{1}{2}}\mm{d}t	\nonumber\\&
\lesssim B^{\frac{3}{4}}(T)\left(B^{{3}}(T)+ \int_0^{T}\langle t\rangle^{-\frac{5a}{2}}\mm{d}t\left( \int_0^{T}\left( \left(\langle t\rangle^{a} \mathfrak{D}\right) ^{2} \langle t\rangle^{\frac{3a}{2}} \underline{\mathcal{D}}\right) ^{\frac{1}{3}}\mm{d}t\right)^{ {3}}\right)^{\frac{1}{4}} \nonumber \\
&\lesssim B^{\frac{3}{2}}(T).
\label{lemisdfAa}
\end{align}

Exploiting the boundary condition of $v_3$ in \eqref{n1},  the estimate of $v_t$ in \eqref{dctv3}, the dissipative estimates of $\nabla_{\mm{h}}\varrho$ in \eqref{Pfsda2}, Young's inequality and the integration by parts, we deduce from \eqref{dctv1}, \eqref{dctv3} and \eqref{lemiq11axc} that for any sufficiently large  constant $\chi$,
\begin{align}
&\sup_{0\leqslant  t\leqslant  T} \Big\{\langle t\rangle^{\ell_1}\tilde{\underline{ {\mathfrak{E}}}}^{{\ell_1},{\ell_2}}_1\Big\} + \int_0^{T} \langle t\rangle^{\ell_1}\tilde{\underline{ {\mathfrak{D}}}}^{{\ell_1},{\ell_2}}_1  \mm{d}t  \nonumber\\&\lesssim \int_0^{T}  \ell_1\langle t\rangle^{{\ell_1}-1}\|(\varrho, v)\|^2_{{\ell_2},1} \mm{d}t  + \|(\varrho^0,v^0)\|^2_{{\ell_2},1}+	 B^{3/2}(T)\nonumber \\
&\qquad+\chi^{-1}\int_0^T \langle t\rangle^{{\ell_1}} \|\nabla v\|^2_{{\ell_2},1}\mm{d}t \hbox{ for }\begin{cases}
{\ell_1}=0,\  {\ell_2}=-s;\\ {\ell_1}=a,\ {\ell_2}=0,
\end{cases} \label{lemsafiqsaf11Ac}
\end{align}		
where we have defined that		
$$\tilde{\underline{ {\mathfrak{E}}}}^{{\ell_1},{\ell_2}}_{1}:=\|(\varrho,v)\|^2_{{\ell_2},1} $$ and
$$\tilde{\underline{ {\mathfrak{D}}}}^{{\ell_1},{\ell_2}}_1:= \chi^{-1}\|\varrho\|^2_{1+{\ell_2},1}+\|(\nabla v,v_t)\|^2_{{\ell_2},0}  .$$

Noting that
$$ \int_0^T\langle t\rangle^{{\ell_1}}({\underline{\mathfrak{E}}}^{{s}}{{\mathfrak{E}}})^{\frac{1}{2(1+s)}}
 \underline{\mathfrak{D}} ^{\frac{1}{2}}\mm{d}t\lesssim B(T),$$	
 and then making use of \eqref{nav1}--\eqref{nav1var},
   \eqref{im3}, the dissipative estimate of $\nabla_{\mm{h}}\varrho$ in  \eqref{r4} with $j=2$, the estimate of $\varrho_t$ in  \eqref{202221401177321}, \eqref{lemisdfAa},  the Hodge-type elliptic estimate \eqref{202005021302}, the interpolation inequality \eqref{201807291850} and  Young's inequality, we  infer from \eqref{dctv2}, \eqref{dctv4}, \eqref{202221401177321aac} and \eqref{lemsafiqsaf11Ac} that
\begin{align}
&\sup_{0\leqslant  t\leqslant  T}\Big\{\langle t\rangle^{\ell_1}\tilde{\underline{ {\mathfrak{E}}}}^{{\ell_1},{\ell_2}}_2\Big\}+ \int_0^{T} \langle t\rangle^{\ell_1}\tilde{\underline{ {\mathfrak{D}}}}^{{\ell_1},{\ell_2}}_2 \mm{d}t \nonumber\\& \lesssim
\int_0^{T}  \langle t\rangle^{{\ell_1}}\big( \ell_1\langle t\rangle^{-1}(\|\varrho\|^2_{{\ell_2},3} +\|  v\|^2_{{\ell_2},0}+\| \nabla v\|^2_{{\ell_2},1})  +\|  (\nabla v,v_t)\|^2_{{\ell_2},1} +\|\varrho\|^2_{1+{\ell_2},1}\big)\mm{d}t \nonumber \\
&\quad\ +(\|\varrho^0\|^2_{{\ell_2},3} +\|  v^0\|^2_{{\ell_2},2})+ \sup_{0\leqslant  t\leqslant  T}\big\{\langle t\rangle^{\ell_1}\left( \|  v\|^2_{{\ell_2},1}+\| \varrho\|^2_{{\ell_2},2} \right)\big\}\nonumber \\
&\quad\ +	 B^{3/2}(T) \hbox{ for }\,\begin{cases}
{\ell_1}=0,\, {\ell_2}=-s;\\  	{\ell_1}=a,\,{\ell_2}=0,
\end{cases} \label{lemisdfA}
\end{align}	
where we have defined that
$$\tilde{\underline{ {\mathfrak{E}}}}^{{\ell_1},{\ell_2}}_2:= \|\varrho\|^2_{{\ell_2},3} +\|  v\|^2_{{\ell_2},2}
$$
and
$$\tilde{\underline{ {\mathfrak{D}}}}^{{\ell_1},{\ell_2}}_2:= \|\varrho\|^2_{1+{\ell_2},2}+\|\nabla v\|^2_{{\ell_2},2} +\|v_t\|^2_{{\ell_2},1} .
$$

In addition, using  \eqref{imvarrho}, \eqref{im3},  the estimates of $\varrho_t$ in  \eqref{202221401177321} with $(\ell,j)=(-s,2)$ and $(\ell,j)=(0,2)$, the estimates of $v_t$ in  \eqref{vt0} and the relation
$$a=5/3-s,$$
we obtain that,  for sufficiently large  constant $\chi$,
\begin{align}
\langle t\rangle^{\ell_1}\underline{ {\mathfrak{E}}}&\lesssim\langle t\rangle^{\ell_1}\left( \|\varrho_t\|^2_{-s,2} +\|  v_t\|^2_{-s,0}+\chi \tilde{\underline{ {\mathfrak{E}}}}^{{\ell_1},{\ell_2}}_1 + \tilde{\underline{ {\mathfrak{E}}}}^{{\ell_1},{\ell_2}}_2\right)\nonumber\\&\lesssim\langle t\rangle^{\ell_1}
\Big(\chi\tilde{\underline{ {\mathfrak{E}}}}^{{\ell_1},{\ell_2}}_1+  \tilde{\underline{ {\mathfrak{E}}}}^{{\ell_1},{\ell_2}}_2\Big)+ B^{3/2}(T)\hbox{ for }{\ell_1}=0,\ {\ell_2}=-s, \label{eea1}\\
\langle t\rangle^{\ell_1}	{ {\mathfrak{E}}}&\lesssim	\langle t\rangle^{\ell_1}\left( \|\varrho_t\|^2_{2} + \|  v_t\|^2_{0} + \tilde{\underline{ {\mathfrak{E}}}}^{{\ell_1},{\ell_2}}_2+\chi \tilde{\underline{ {\mathfrak{E}}}}^{{\ell_1},{\ell_2}}_1\right)\nonumber
\\&\lesssim \langle t\rangle^{\ell_1}\Big(\chi\tilde{\underline{ {\mathfrak{E}}}}^{{\ell_1},{\ell_2}}_1+\tilde{{ {\mathfrak{E}}}}^{{\ell_1},{\ell_2}}_2\Big)+ B^{3/2}(T)\hbox{ for }
{\ell_1}=a,\ {\ell_2}=0, \label{eea2}
\\
\int_0^T\langle t\rangle^{\ell_1}	\underline{ {\mathfrak{D}}}\mm{d}t&\lesssim \int_0^T
\langle t\rangle^{\ell_1}\Big(\tilde{\underline{ {\mathfrak{D}}}}^{{\ell_1},{\ell_2}}_2+ \|\varrho_t\|^2_{-s,2} +\chi\tilde{\underline{ {\mathfrak{D}}}}^{{\ell_1},{\ell_2}}_1\Big)\mm{d}t\nonumber\\&
\lesssim  B^{3/2}(T)+
\int_0^T\langle t\rangle^{\ell_1} \Big(
\chi\tilde{\underline{ {\mathfrak{D}}}}^{{\ell_1},{\ell_2}}_1+\tilde{\underline{ {\mathfrak{D}}}}^{{\ell_1},{\ell_2}}_2\Big)\mm{d}t \hbox{ for } {\ell_1}=0,\ {\ell_2}=-s,
\end{align}
and
\begin{align}
\int_0^T\langle t\rangle^{\ell_1} 	{ {\mathfrak{D}}}\mm{d}t&\lesssim 	
\int_0^T\langle t\rangle^{\ell_1}\Big(\|\varrho_t\|^2_{2}+\chi\tilde{\underline{ {\mathfrak{D}}}}^{{\ell_1},{\ell_2}}_1 +
\tilde{\underline{ {\mathfrak{D}}}}^{{\ell_1},{\ell_2}}_2\Big)\mm{d}t\nonumber\\&\lesssim B^{3/2}(T)+\int_0^T\langle t\rangle^{\ell_1}\Big(\chi\tilde{\underline{ {\mathfrak{D}}}}^{{\ell_1},{\ell_2}}_1+ \tilde{{ {\mathfrak{D}}}}^{{\ell_1},{\ell_2}}_2\Big)\mm{d}t \hbox{ for }
{\ell_1}=a,\ {\ell_2}=0. \label{eea4}
\end{align}	

  Consequently, exploiting  \eqref{eea2}, \eqref{eea4}, the interpolation inequality \eqref{201807291850} and Young's inequality, we deduce \eqref{qecc1} from \eqref{lemsafiqsaf11Ac} and \eqref{lemisdfA}  with ${\ell_1}=a$ and ${\ell_2}=0$.

(4)  Similarly, using \eqref{eea1}, the interpolation inequality \eqref{201807291850} and Young's inequality, we infer \eqref{lem2dx1}  from  \eqref{lemsafiqsaf11Ac} and \eqref{lemisdfA} with the case that ${\ell_1}=0$, and ${\ell_2}=-s$.
\hfill $\Box$
\end{pf}

\begin{pro}
\label{2508191327}
Under the \emph{a priori} assumption \eqref{imvarrho} with sufficiently small $\delta^*$, the smallness of which depends on $g$,  $\mu$, $\varepsilon$ and $\bar{\rho}$, it holds that
\begin{align}
B(T)\lesssim \|(\nabla\varrho^0,v^0)\|_5^2 + \|(\nabla\varrho^0,v^0)\|_{-s,2}^2.
 \label{2025081321621}
 \end{align}
\end{pro}
 \begin{pf}
Utilizing the definition of $a$  in \eqref{disfa}, \eqref{201807291850},  and Young's inequality, we have, for sufficiently large  constant $\chi$,
\begin{align}
\langle t\rangle^{a-1}\|\varrho\|^2_{3}&\lesssim\langle t\rangle^{a-1}\left( \|\varrho\| _{-s,3}\|\varrho\|^{s}_{1,3}\right)^{\frac{2}{1+s}} \nonumber\\
&\lesssim\langle t\rangle^{a-1}\left( \|\varrho\|_{-s,3}\left( \|\varrho\|^{2}_{1,2}\|\varrho\|_{1,5}\right)^{\frac{s}{3}} \right) ^{\frac{2}{1+s}} \nonumber\\
&\lesssim\left( \chi\langle t\rangle^{(a-1)}\right) ^{1+s} \langle t\rangle^{-\frac{2as}{3}}\|\varrho\|^{2}_{-s,3} +\chi^{-\frac{1+s}{s}}\left( \langle t\rangle^{a}\|\varrho\|^{2}_{1,2}\|\varrho\|_{1,5}\right)^{\frac{2}{3}}.\label{2508311447a}
\end{align}	
Similarly, we can also get
\begin{align}
\langle t\rangle^{a-1}\|v\|^2_{2}&\lesssim
\chi^{1+s}\langle t\rangle^{(a-1)(1+s)-\frac{2as}{3}}\|v\|^{2}_{-s,0} +\chi^{-\frac{1+s}{s}}\langle t\rangle^{a}\|\nabla v\|^{2}_{2},\nonumber
\end{align}	
which, together with \eqref{2508311447a} and \eqref{ab1}, implies that
\begin{align}\label{2508311447}
	&\int_0^T\langle t \rangle^{a-1}(\|\varrho\|^2_{3} + \|v\|^2_{2}) \mathrm{d}t
	\nonumber\\&\lesssim \int_0^T \left( \chi^{1+s} \langle t \rangle^{(a-1)(1+s)-\frac{2as}{3}} \underline{\mathfrak{E}}(t) + \chi^{-\frac{1+s}{s}} \left( \langle t \rangle^{a} \mathfrak{D}(t) + \|\varrho\|^{2}_{1,5} \right)\right)\mathrm{d}t \nonumber\\&\lesssim  \chi^{1+s} \sup_{0 \leqslant t \leqslant T} \left\lbrace\underline{\mathfrak{E}}(t) \right\rbrace + \int_0^T \left(\chi^{-\frac{1+s}{s}} \left( \langle t \rangle^{a} \mathfrak{D}(t) + \|\varrho\|^{2}_{1,5} \right)\right)\mathrm{d}t.
\end{align}
In view of the above estimate, we deduce from \eqref{qecc1} and \eqref{qec1}  that
\begin{align}
&\sup_{0\leqslant  t\leqslant  T}\left\lbrace \underline{\mathfrak{E}}(t)+\langle t\rangle^{3a/2}\underline{\mathcal{E}}(t)+\langle t\rangle^{a}{\mathfrak{E}}(t)\right\rbrace +\int_0^{T}\big(\langle t\rangle^{3a/2}\underline{\mathcal{D}}(t)+ \underline{\mathfrak{D}}(t)+ \langle t\rangle^{a}{\mathfrak{D}}(t)\big)  \mm{d}t\nonumber\\&\lesssim\mathcal{E}_{\mm{total}}(0)+ B^{3/2}(T)
+\int_0^{T}\|\varrho\|^{2}_{1,5}\mm{d}t.\label{lem2dx3}
\end{align}	

In addition, using H\"older's inequality and the facts of $\theta\in(0,({19-\sqrt{349}})/{6})$ and $a=\theta+{2}/{3}$, we have	
\begin{align}
\int_0^{T}\mathfrak{N}(t)\mm{d}t\lesssim &\int_0^{T}
\sqrt{\mathcal{E}}
\bigg(\left( {\mathcal{E}}\langle t\rangle^{\frac{3a}{2}} \underline{\mathcal{D}}\right)^{\frac{1}{2}}\langle t\rangle^{-\frac{3a}{4}}+\left(\left( {\mathcal{E}}\langle t\rangle^{a} {\mathfrak{D}}\right)^{\frac{1}{2}}\right.\nonumber \\
&  \left.+\left( \langle t\rangle^{a} {\mathfrak{E}}{\mathcal{D}}\right)^{\frac{1}{2}}\right) \left( \langle t\rangle^{a} {\mathfrak{E}}\right)^{\frac{1}{2}}\langle t\rangle^{-a} + \left( \langle t\rangle^{\frac{5a}{2}}   {\mathfrak{E}}\mathcal{E}\mathcal{D} \underline{\mathcal{D}}\right)^{\frac{1}{4}}\langle t\rangle^{-\frac{5a}{8}}\bigg)\mm{d}t\nonumber\\ 
\lesssim& \left(B^{{3}}(T) \int_0^{T}\langle t\rangle^{-\frac{5a}{2}}\mm{d}t \left( \int_0^{T}\left(\left(\langle t\rangle^{a} \mathfrak{D}\right) ^{ {2} } \langle t\rangle^{\frac{3a}{2}} \underline{\mathcal{D}}\right) ^{\frac{1}{3}}{\mm{d}t}\right)^{3}  \right)^{\frac{1}{4}} +B^{\frac{3}{2}}(T) \nonumber \\
\lesssim& B^{\frac{3}{2}}(T). \label{2510031549}
 \end{align}
Thanks to the above estimate, we can derive from \eqref{omessetsimQ} and \eqref{lem2dx3} that
\begin{align} B(T) \lesssim \mathcal{E}_{\mathrm{total}}(0).
 \label{20250813216211}
 \end{align}

Finally, using \eqref{vt0}--\eqref{2508181314}, we have
$$ \mathcal{E}_{\mathrm{total}}(0)\lesssim \|(\nabla\varrho^0,v^0)\|_5^2 + \|(\nabla\varrho^0,v^0)\|_{-s,2}^2.$$
Putting the above estimate into \eqref{20250813216211} yields \eqref{2025081321621}.
This completes the proof.  \hfill$\Box$
\end{pf}	

\section{Proof of Theorem \ref{thm2}}\label{subsec:08a}
This section is devoted to the proof of Theorem \ref{thm2}. To begin with, we shall introduce a local well-posedness result for the QRT problem.
\begin{pro}\label{202102182115}
Let $\mu$, $\varepsilon$ be positive constants, $s\in (0,1)$ and $0<\bar{\rho}\in {C^8}[0,h]$. For any  $( \varrho^0,v^0)\in  { {H}}^{\bar{\rho},6}_{-s,3}\times  {H}^{\diamond,\sigma,-s}_{5,2}$ satisfying
\begin{align*}
 \|(\nabla \varrho,v)\|_5 + \|(\varrho,v)\|_{-s,2}   \leqslant b, \end{align*}
necessary compatibility conditions and the positive lower-bound condition of initial density (i.e. $\inf_{x\in\Omega}\big\{{\rho}^{0}(x)\big\}>0$,	where $\rho^0:=\varrho^0+\bar{\rho}$), then there exists  $T_*>0$, depending on  $\mu$,  $\varepsilon$, $s$ and $b$, such that  the QRT problem \eqref{1a}--\eqref{n1}  admits a unique local(-in-time) classical solution $(\varrho,v)$ with an associated pressure $\beta$. Moreover $(\varrho,v,\nabla\beta)\in{\mathfrak{P}} _{T_*}\times \mathfrak{W}_{T_*} \times C^0(\overline{I_{{T_*}}},{H}^5)$, $\beta\in  C^0(\overline{I_{{T_*}}},L^6)$ and \begin{align*}
&0<\inf\limits_{x\in\Omega}\big\{{\rho}^{0}(x)\big\}\leqslant {\rho}(t,x)\leqslant \sup\limits_{x\in\Omega}\big\{{\rho}^{0}(x)\big\}
\mbox{ for any }(t,x)\in I_{T_*}\times\Omega .
\end{align*}
\end{pro}
\begin{pf}
Since Proposition \ref{202102182115} can be easily proved by the standard iteration method as in \cite{JFJSWWWN,JFJSZYYO,MR4246913}, we omit the details of the proof. \hfill $\Box$
\end{pf}
\begin{rem}
{Since $\varrho\in{\mathfrak{P}} _{T_*} $, $\varrho$ satisfies the boundary conditions in \eqref{2508171135}. Such boundary conditions satisfied by $\varrho$ are deduced from the mass equation and the boundary conditions of $(\rho^0,v^0)$ in Proposition \ref{202102182115}. Next we verify this fact for reader's convenience.}

In view of the mass equation \eqref{1a}$_1$ and the boundary condition of $v_3$ in \eqref{n1}, it holds that
\begin{equation}\label{n22a1}
\varrho_t+v_{\mm{h}}\cdot \nabla_{\mm{h}}\varrho=0\  \text{on}\ \partial \Omega.
\end{equation}
Taking the inner product of the above identity
and  $\varrho$ in $L^2(\partial\Omega)$, and  then using the integration by parts and  the embedding inequality of $H^2\hookrightarrow C^0(\Omega)$ in \eqref{esmmdforinfty}, we derive that
\begin{equation}\label{n22a2}
\frac{\mm{d}}{\mm{d}t}	\int_{\partial\Omega}|\varrho|^2\mm{d}x_{\mm{h}}=-\frac{1}{2}\int_{\partial\Omega} v_{\mm{h}}\nabla_{\mathrm{h}} |\varrho|^2 \mm{d}x_{\mm{h}}  \lesssim \|\mm{div}_{\mathrm{h}}v_{\mm{h}}\|_{2} \int_{\partial\Omega} |\varrho|^2\mm{d}x_{\mm{h}}  .
\end{equation}
Noting  $\int_0^T \|\mm{div}_{\mathrm{h}}v_{\mm{h}}\|_{2}\mm{d}\tau<\infty$, thus applying Gronwall's inequality to \eqref{n22a2} yields
\begin{equation}\label{n22asf}
\|\varrho\|^2_{L^2(\partial\Omega)}=0,
\end{equation}
which implies
\begin{equation}\label{n22a}
\varrho|_{\partial\Omega}=0.
\end{equation}

Applying $\partial_3$ to \eqref{1a}$_1$, and then using  the boundary condition of $v_3$ in \eqref{n1}, we can compute out that
\begin{equation}\label{n22ac}
\partial_{3}\rho_t+v _{\mm{h}}\cdot \nabla_{\mm{h}} \partial_{3} \rho+\partial_{3}v_3\partial_{3}\rho=0 \text{ on } \partial \Omega.
\end{equation}
Exploiting the incompressible condition \eqref{1a}$_3$ and the initial condition $ \partial_3\rho^0 |_{\partial\Omega}=0$, we easily follow the argument of  \eqref{n22a} to deduce from \eqref{n22ac} that
\begin{equation}\label{n22sa}
\partial_3\rho |_{\partial\Omega}=0.
\end{equation}
In view of \eqref{nfdsa22aa},  \eqref{n22a} and \eqref{n22sa}, the boundary condition \eqref{2508171135} holds for $i=0$.

Applying $\partial_3^{2}$ to \eqref{1a}$_1$, and then making use of the boundary conditions of $ v $ in
\eqref{n1},  $\bar{\rho}''$ in \eqref{0102n1} and $\varrho$ in \eqref{n22a}, and the incompressible condition \eqref{1a}$_3$, we have
\begin{equation} 	\partial^2_{3}\varrho_t+v_{\mm{h}}\cdot \nabla_{\mm{h}}\partial^2_{3} \varrho +2\partial_{3}v_3\partial^2_{3}\varrho=0 \text{ on } \partial \Omega,
\label{2508171331}
\end{equation}
which obviously implies that
\begin{equation}\label{nfdsa22a}
\partial_3^2\varrho |_{\partial\Omega}=0.
\end{equation}
Thus we further deduce from the boundary condition of $\bar{\rho}''$ in \eqref{0102n1}, \eqref{nfdsa22aa},  \eqref{n22a} and \eqref{nfdsa22a}  that \eqref{2508171135} holds for $i=1$.
	
Now  recalling the derivation of \eqref{vbian}  with $i=1$, we have
\begin{align}
(\partial^3_3v_{\mm{h}},\partial^{2}_3v_{3}) |_{\partial\Omega}=0.\label{vbianssfaadfafsds}
\end{align}
Similarly to \eqref{n22ac},  by further using the boundary conditions of $(\rho,\varrho)$ in  \eqref{varrho}, $v$ in \eqref{vbianssfaadfafsds}   and \eqref{vbian} with $i=0$, we can compute out that
	\begin{equation*}
		\partial^3_{3}\rho_t+v _{\mm{h}}\cdot \nabla_{\mm{h}} \partial^3_{3} \rho+3\partial_{3}v_3\partial^3_{3}\rho=0 \text{ on } \partial \Omega,
	\end{equation*}
which, together the initial condition $ \partial_3\rho^0 |_{\partial\Omega}=0$, implies that
	\begin{equation}\label{n22sasafx}
		\partial^3_3\rho |_{\partial\Omega}=0.
	\end{equation}	

Similarly to \eqref{2508171331}, exploiting the boundary conditions of \eqref{0102n1}, \eqref{varrho}--\eqref{omega} with $i=0$ and $1$ and \eqref{vbianssfaadfafsds}, we have
		\begin{equation*} 	\partial^4_{3}\varrho_t+v_{\mm{h}}\cdot \nabla_{\mm{h}}\partial^4_{3} \varrho +4\partial_{3}v_3\partial^4_{3}\varrho=0 \text{ on } \partial \Omega,
		\end{equation*}
		which obviously implies that
		\begin{equation}\label{nfdsa22aasdx}
			\partial_3^4\varrho |_{\partial\Omega}=0.
		\end{equation}
Finally we further deduce from the boundary condition of $\bar{\rho}''$ in \eqref{0102n1}, \eqref{nfdsa22aa},  \eqref{n22sasafx} and \eqref{nfdsa22aasdx}  that \eqref{2508171135} holds for $i=2$. This completes the verification of the boundary conditions in \eqref{2508171135}.
\end{rem}

Thanks to the \emph{a priori} stability estimates in Propositions \ref{2508191327}, and the local well-posedness result in Proposition \ref{202102182115}, we can easily establish Theorem \ref{thm2}.
Next, we briefly sketch the proof.

Let $(\varrho^0,v^0)$ satisfy the assumptions in Theorem \ref{thm2}.
By the embedding inequality \eqref{esmmdforinfty},  there exists a constant $\delta_1>0$ such that, if $\|\varrho^0\|_6\leqslant \delta_1$, it holds that
$\inf_{x\in\Omega}\big\{{\rho}^{0}(x)\big\}>0$.  From now on, we choose  $\delta$ in Theorem \ref{thm2} to be less than $\delta_1$.

In view of Proposition \ref{202102182115}, there exists a unique local classical solution $(\varrho,v, \beta)$ to the QRT problem \eqref{1a}--\eqref{n1}  with the maximal existence time $T^{\max}$, which satisfies
\begin{itemize}
\item for any $a\in I_{T^{\max}}$,
$(\varrho,v, \nabla\beta)\in {\mathfrak{P}} _{a}\times {\mathfrak{W}_{a} }\times C^0(\overline{I_a},{H}^5) $;
\item $\limsup_{t\to T^{\max} }\|(\nabla \varrho, v)( t)\|_5=\infty$  if  $T^{\max}<\infty$.
\end{itemize}
Moreover,
\begin{equation}
0<\inf\limits_{x\in\Omega}\big\{{\rho}^{0}(x)\big\}\leqslant {\rho}(t,x)\leqslant \sup\limits_{x\in\Omega}\big\{{\rho}^{0}(x)\big\}\mbox{ for any }(t,x)\in I_{T^{\max}}\times\Omega .
\label{202402081718}
\end{equation}

Recalling the regularity of  $(\varrho,v, \beta)$ above, we can verify that the solution $(\varrho,v,\nabla\beta)$ satisfies \eqref{2025081321621} under the assumption \eqref{imvarrho} with sufficiently small $\delta^*$. More precisely,
there exist positive constants $c_1 \geqslant 1$ and $\delta_2\leqslant \delta_1$ such that
\begin{align}
\label{omessestsimQ}
&B(T)\leqslant {c}_1\big(\|(\nabla\varrho^0,v^0)\|_5^2 + \|(\nabla\varrho^0,v^0)\|_{-s,2}^2\big),
\end{align}
if
$$B(T)\leqslant  \delta_2^2,$$ where the constants $c_1$ and $\delta_2$ depend on the domain $\Omega$, the other known physical parameters/functions, and the parameter $\theta$.

Let $\delta\leqslant \delta_2/\sqrt{2c_1}$ and
\begin{align}
T^{*}=&\sup \left\{\tau\in I_{T^{\max}}~ \left|~B(T)\leqslant 2 c_1\delta^2 \ \mbox{ for any }
t\leqslant\tau \right.\right\}.\nonumber
\end{align}
We easily see from \eqref{omessestsimQ} that the definition of $T^*$ makes sense due to the assumption that
\begin{align}
\label{2510051658}
 \|(\nabla\varrho^0,v^0)\|_5  + \|(\nabla\varrho^0,v^0)\|_{-s,2} \leqslant    \delta   .\end{align} Thus, to show the existence of a global(-in-time) solution, it suffices to verify $T^*=\infty$. We shall prove this fact by contradiction below.
	
Assume $T^*<\infty$, by Proposition \ref{202102182115} and \eqref{202402081718}, we have
\begin{align}T^*\in I_{T^{\max}} \label{20222022519850}
\end{align}
and
\begin{equation*}
B(T^*) \leqslant  2 c_1 \delta^2 \leqslant  {\delta^2_2},
\end{equation*}
which, together with the assertion in \eqref{omessestsimQ} and the initial condition \eqref{2510051658}, further yields
\begin{align}
B(T^*)\leqslant c_1
\left( \|(\nabla \varrho,v)\|_5 + \|(\varrho,v)\|_{-s,2}\right)  \leqslant   c_1 \delta^2.  \label{2020103261534}
\end{align}

By \eqref{20222022519850}, \eqref{2020103261534} and the strong continuity of  $( \varrho,v)\in C^0([0,T^{\max}),H^6_{-s,3}\times H^5_{-s,2})$, we deduce that
there is a  constant $\tilde{T}\in (T^*,T^{\max})$, such that
\begin{align}
B(\tilde{T})  \leqslant 2 c_1\delta^2  , \nonumber
\end{align}
which contradicts with the definition of $T^*$. Hence, $T^*=\infty$, which implies $T^{\max}=\infty$.
This completes the proof to the existence of a global solution. The uniqueness of the global solution is obvious due to the one of the  local solutions in Proposition \ref{202102182115}. This completes the proof of Theorem \ref{thm2}.
	
\section{Appendix: analysis tools}\label{sec:09}
\renewcommand\thesection{A}
This section  is devoted to providing some mathematical results, some of which have been frequently used in previous sections. We should point out that $\Omega$ and the simplified notations appearing in what follows are  the same as those defined in Section \ref{202402081729}. In addition, $a\lesssim b$ still denotes $a\leqslant cb$, where the positive constant $c$ depends on the parameters and the domain in the lemmas in which $c$ appears.
\begin{lem}
\label{201806171834}  Embedding inequalities:
\begin{enumerate}[(1)]
  \item Let $D\subset \mathbb{R}^3$ be a domain satisfying the cone condition, then (see \cite[Theorems 4.12]{ARAJJFF})
\begin{align}
&\label{esmmdforinfty}\|f\|_{C^0(\bar{D})}= \|f\|_{L^\infty(D)}\lesssim\| f\|_{H^2(D)}\mbox{ for any }f\in  H^2(D).
\end{align} It should be noted that the domain $\mathbb{R}^2\times (0,h)$ satisfies the cone condition.
  \item  Let $0< s< 2$, $1< p< q< \infty$ and  $1/ q+ s/2= 1/ p$.  If  $\Lambda^{-s}_{\mathrm{h}}f$ and $f$  belong to the function space $L^2(\mathbb{R}^2)$, then we have  (see \cite[Section 1.2]{MR290095}):
\begin{align}\label{v2}
\left\|\Lambda^{-s}_{\mm{h}}f\right\|_{L^q(\mathbb{R}^2)}\lesssim\left\|f\right\|_{L^p(\mathbb{R}^2)}.
\end{align}
\end{enumerate}
\end{lem}
\begin{lem}
A Poincar\'e-type inequality \cite[Lemma 10.6]{YnaBVG}:
It holds that
\begin{equation}
\|{f}\|_{0}\lesssim\|{f}\|_{L^2(\partial\Omega)}+\|{\partial_3f}\|_{0}
\mbox{ for all }f\in H^{1}.
\label{2202402040948}
\end{equation}
\end{lem}
\begin{lem}\label{pro4a}
A	Hodge-type elliptic estimate  (referring to \cite[Lemma A.4]{ZHAOYUI}):
Let $i\geqslant1$,
then
\begin{align}
&\label{202005021302}
\|\nabla w\|_{i-1}
\lesssim\|(\mm{curl}w,\mm{div}w)\|_{i-1}\mbox{ for any  }w\in{H}^i\mbox{ satisfying }w_3|_{\partial\Omega}=0.
\end{align}
\end{lem}
\begin{lem}	\label{201806171834a}\begin{enumerate}[(1)]
\item Interpolation inequalities in $L^p$ \cite[Lemma A.3]{jiang2025asymptotical}:
It holds that
\begin{align}\label{leml3}
&\| f\|_{L^4}^2\lesssim\|\nabla_{\mm{h}} f\|_{0}\| f\|_{1}\mbox{ for any }f\in H^1,\\
&\|f\|_{L^\infty}^2\lesssim\|\nabla_{\mm{h}} f\|_{1}\| f\|_{2}\mbox{ for any }f\in H^2\label{lemge4}.
\end{align}
\item
Interpolation inequality in $H^j$ \cite[Theorem 5.2]{ARAJJFF}: Let $D$ be a domain in $\mathbb{R}^n$ satisfying the cone condition, then, for any given $0\leqslant j< i$,
\begin{equation}
\|f\|_{H^j(D)}\lesssim\|f\|_{L^2(D)}^{1-\frac{j}{i}}
\|f\|_{H^i(D)}^{\frac{j}{i}}\lesssim\epsilon^{-j/(i-j)}\|f\|_{L^2(D)} +\epsilon \|f\|_{H^i(D)} \label{201807291850}
\end{equation}  for any $f\in H^i(D)$ and  for any $\epsilon>0$,
where the two constants $c $ in \eqref{201807291850} are independent of $\varepsilon$.
\item Interpolation inequalities for negative derivatives
(referring to \cite[Lemma A.4]{MR3005540}):
Let $s \geqslant 0$ and $\ell \geqslant 0$. It holds that
\begin{align}\label{v1}
\| \Lambda^{\ell}_{\mm{h}} f\|_{L^2(\mathbb{R}^2)}^ {\ell+s+1} \lesssim  \|\Lambda^{-s}_{\mm{h}} f\|_{L^2(\mathbb{R}^2)}\|\Lambda^{1+\ell}_{\mm{h}} f\|_{L^2(\mathbb{R}^2)}^{ {\ell+s} } 
\end{align}
and
\begin{align}\label{v11}
\|\Lambda^{\ell-s}_{\mm{h}} f\|_{L^2(\mathbb{R}^2)}^{\ell+s}\leqslant  \| \Lambda^{-s}_{\mm{h}} f\|_{L^2(\mathbb{R}^2)}^s\| \Lambda^{\ell}_{\mm{h}} f\|_{L^2(\mathbb{R}^2)}^{\ell},
\end{align}
if the norms/semi-norms of $f$ on the right hand side of the above inequalities are finite.
\end{enumerate}
\end{lem}

 \begin{lem}\label{xfsddfsf20180508a}
Real interpolation  inequality: Let $2\leqslant p < +\infty$, $ 0\leqslant \theta\leqslant 1 $ and $ \alpha   \geqslant 0$ satisfy
\begin{align}\label{v111a}
\alpha +2\left( \frac{1}{2} - \frac{1}{p} \right) =   \ell\theta,
\end{align}
then we have
\begin{align}\label{v111}
\|\Lambda^{\alpha}_{\mm{h}} f\|_{L^p(\mathbb{R}^2)} \lesssim \| f\|_{L^2(\mathbb{R}^2)}^{1-\theta} \|\Lambda^{\ell}_{\mm{h}} f\|_{L^2(\mathbb{R}^2)}^{\theta},
\end{align}
where  we have assumed that the norms on the right hand side of the above inequality are finite.
\end{lem}
\begin{pf}
Let $s= \ell\theta -\alpha$, then 
$$s=2\left( \frac{1}{2} - \frac{1}{p} \right)\in (0,1).$$
Exploiting \eqref{v2}, H\"older's inequality and Parseval theorem, we obtain
\begin{align}
\|\Lambda^{\alpha}_{\mm{h}} f\|_{L^p(\mathbb{R}^2)}=&\|\Lambda^{-s}_{\mm{h}}(\Lambda^{s+\alpha}_{\mm{h}} f)\|_{L^p(\mathbb{R}^2)} \nonumber \\
\lesssim & \|\Lambda^{ \ell\theta }_{\mm{h}} f\|_{L^2(\mathbb{R}^2)}\lesssim \| f\|_{L^2(\mathbb{R}^2)}^{1-\theta} \|\Lambda^{\ell}_{\mm{h}} f\|_{L^2(\mathbb{R}^2)}^{\theta}.\nonumber 
\end{align}
 The proof is complete. \hfill $\Box$
\end{pf} 

\begin{lem}\label{xfsddfsf20180508} Product estimates in $H^i$  (referring to \cite[Lemma A.3]{jiang2023rayleigh}):
Let $D\subset \mathbb{R}^3$ be a   domain  satisfying the cone condition, then
\begin{align}
\label{fgestims}
\|f\varphi\|_{H^i(D)}\lesssim \begin{cases}
\|f\|_{H^1(D)}\|\varphi\|_{H^1(D)} & \hbox{for }i=0; \\
\|f\|_{H^i(D)}\|\varphi\|_{H^2(D)} & \hbox{for }0\leqslant i\leqslant 2;  \\
\|f\|_{H^2(D)}\|\varphi\|_{H^i(D)}+\|f\|_{H^i(D)}\|\varphi\|_{H^2(D)}& \hbox{for }i\geqslant3 ,
\end{cases}
\end{align}
if the norms on the right hand side of the above inequalities are finite.		
\end{lem}
	
\begin{lem}\label{2024v3}
Product estimates in $\dot{H}^{-s}$: Let $0<s<1$ and $i\in \mathbb{N}$, it holds that
\begin{align}\label{4v3}
\left\|f\varphi\right\|_{-s,0}&\lesssim	\left\|f\right\|_{0} \left(\|\varphi\|^{s}_{-s,0}\|\varphi\|_{1,0} (\|\partial_3 \varphi\|^{s}_{-s,0}\|\partial_3 \varphi\|_0+\|\varphi\|^{s}_{-s,0}\|\varphi\| _{1,0})\right)^{\frac{1}{2(1+s)}}.
\end{align}
In particular,  if $\varphi|_{\partial{\Omega}}=0$, we further have
\begin{align}\label{4v312}
\left\|f\varphi\right\|_{-s,0}&\lesssim\left\|f\right\|_{0}\left(\|\partial_3 \varphi\|^{s}_{-s,0}\|\partial_3 \varphi\|_{1,0}\right)^{\frac{1}{1+s}}
\end{align}
and
\begin{align}\label{4v313}
\left\|f\varphi\right\|_{-s,i}&\lesssim\sum_{j=0}^{j=i}\left\|f\right\|_{j}\mathfrak{U}^{i-j}(\varphi),
\end{align}
where	we have defined that
\begin{align}\label{ugg}
\mathfrak{U}^{i}(\varphi):=\left(\|\varphi\|^{s}_{-s,i}\|\varphi\|_{1,i} \left(\|\partial_3 \varphi\|^{s}_{-s,i}\|\partial_3 \varphi\|_{i}+\|\varphi\|^{s}_{-s,i}\|\varphi\| _{1,i}\right)\right)^{\frac{1}{2(1+s)}}.
\end{align}	
\end{lem}
\begin{pf}
(1)	In what follows,  $\|f\|_{L_{x_{\mm{h}}}^p}:=\| f\|_{L^p(\mathbb{R}^2)}$ and $\|f\|_{L_{x_3}^p}:=\| f\|_{L^p(0,h)}$, where $1\leqslant p\leqslant \infty$.

By Gagliardo--Nirenberg--Sobolev's inequality (referring to \cite[Theorem  2]{MR109940} with the one dimensional case), one has
\begin{equation}\label{keya}
\|\phi\|_{L^\infty(0,h)}^2\lesssim\|\phi\|_{L^2(0,h)}\| \phi'\|_{L^2(0,h)}+\|\phi\|_{L^2(0,h)}^2\mbox{ for }\phi\in W^{1,2}(0,h).
\end{equation}
Then, making use of \eqref{v2}, \eqref{keya}, Fubini theorem, Minkowski's inequality and H\"older's  inequality,    we have
\begin{align}
\left\|  f\varphi \right\|_{-s,0}	&\lesssim\left\| \|f\varphi\|_{L_{x_{\mm{h}}}^{\frac{2}{s+1}}}\right\| _{L_{x_3}^2}\lesssim\left\|\|f\|_{L_{x_{\mm{h}}}^2}\|\varphi\|_{L_{x_{\mm{h}}}^{\frac{2}{s}}}
\right\|_{L_{x_{3}}^2}\lesssim\|f\|_{L^2}\left\| \|\varphi\|_{L_{x_{\mm{h}}}^{\frac{2}{s}}}\right\|_{L_{x_3}^\infty}\nonumber\\&\lesssim\|f\|_0\left\| \|\varphi\|_{L_{x_3}^\infty}\right\|_{L_{x_{\mm{h}}}^{\frac{2}{s}}}
\lesssim\|f\|_{0}\left\|(\|\varphi\|_{L^2_{x_{3}}}\|\partial_3 \varphi\|_{L^2_{x_{3}}})^{\frac{1}{2}}+\|\varphi\|_{L^2_{x_{3}}}
\right\|_{L_{x_{\mm{h}}}^{\frac{2}{s}}}.\label{keya1}
\end{align}

Thanks to \eqref{v1} and \eqref{v111}, we have
\begin{align}\|f\|_{L_{x_{\mm{h}}}^{\frac{2}{s}}}\lesssim &
\| f\|_{L^2(\mathbb{R}^2)}^{s} \|\Lambda_{\mm{h}} f\|_{L^2(\mathbb{R}^2)}^{1-s}
\lesssim   \left(\| \nabla_{\mm{h}}f\| _{L^{2}_{x_{\mm{h}}}}
\|\Lambda^{-s}_{\mm{h}}f\|^{s}_{L^{2}_{x_{\mm{h}}}}\right)^{\frac{1}{1+s}} , 
\label{2510031236}
\end{align} 
if the norm and semi-norm of $f$ on the right hand side of the above inequalities are finite.
Thus, exploiting \eqref{keya}, \eqref{2510031236} and Minkowski's inequality, we obtain
\begin{align}
&	\left\| \left(\|\varphi\| _{L^2_{x_{3}}}\|\partial_3 \varphi\|_{L^2_{x_{3}}}\right)^{\frac{1}{2}}+\|\varphi\|_{L^2_{x_{3}}}
\right\|_{L_{x_{\mm{h}}}^{\frac{2}{s}}}\nonumber\\&\lesssim	\left\| \|\varphi\|_{L_{x_{\mm{h}}}^{\frac{2}{s}}}\right\| ^{\frac{1}{2}}_{L^2_{x_{3}}}\left\| \|\partial_3\varphi\|_{L_{x_{\mm{h}}}^{\frac{2}{s}}}\right\| ^{\frac{1}{2}}_{L^2_{x_{3}}}+\left\| \|\varphi\|_{L_{x_{\mm{h}}}^{\frac{2}{s}}}\right\|_{L^2_{x_{3}}}\nonumber	\\
&\lesssim	\left\| \left(\|\nabla_{\mm{h}}\varphi\|_{L^{2}_{x_{\mm{h}}}}\|\Lambda^{-s}_{\mm{h}}\varphi\|^{ {s} }_{L^{2}_{x_{\mm{h}}}}\right)^{\frac{1}{1+s}}\right\| ^{\frac{1}{2}}_{L^2_{x_{3}}}\left\|
\left(\|\nabla_{\mm{h}}\partial_3 \varphi\| _{L^{2}_{x_{\mm{h}}}}\|\Lambda^{-s}_{\mm{h}}\partial_3 \varphi\|_{L^{2}_{x_{\mm{h}}}}^s\right)^{\frac{1}{1+s}}\right\| ^{\frac{1}{2}}_{L^2_{x_{3}}}\nonumber \\
&\quad\ +\left\|\left(\| \nabla_{\mm{h}}\varphi\| _{L^{2}_{x_{\mm{h}}}}
\|\Lambda^{-s}_{\mm{h}}\varphi\|^{s}_{L^{2}_{x_{\mm{h}}}}\right)^{\frac{1}{1+s}}\right\|_{L^2_{x_{3}}}\nonumber
\\&\lesssim\left(\|\nabla_{\mm{h}}\varphi\|_{0}\|\varphi
\|^{ {s} }_{-s,0}\left(\|\nabla_{\mm{h}}\partial_3 \varphi\| _{0}\| \partial_3 \varphi\|^{ {s} }_{-s,0}+\|\nabla_{\mm{h}}\varphi\| _0
\| \varphi\|^{ {s} }_{-s,0}\right)\right)^{\frac{1}{2(1+s)}}.\label{keya2}
\end{align}
Thus, putting \eqref{keya2} into \eqref{keya1},
we arrive at \eqref{4v3}.	

If $\varphi|_{\partial{\Omega}}=0$, the estimate \eqref{4v312} follows from \eqref{2202402040948} by the
 Newton--Leibniz formula. In addition,  we can easily obtain \eqref{4v313} from \eqref{4v3} by using the Leibniz formula for the high-order derivatives of the product of two functions. The proof is complete. \hfill $\Box$
\end{pf}

\begin{lem}\label{20lem}The equivalent forms of ${\varepsilon_{\mm{c}}}$: Let $\bar{\rho}'$, ${\underline{\bar{\rho}} } \in L^\infty(0,h)$ and ${\varepsilon_{\mm{c}}}$ be defined by \eqref{2saf01504} with a positive constant $g$,
\begin{align*}
&a_1:={\sup_{0\neq \varpi\in H^1_0 }\frac{\int\bar{\rho}' \varpi ^2 {\mm{d}x}}{ \int |{\underline{\bar{\rho}} }\nabla  \varpi|^2 {\mm{d}x}}},\
a_2:={\sup_{0\neq \varpi\in H^1_0}\frac{\int\bar{\rho}' \varpi^{ 2} {\mm{d}x}}{ \int |{\underline{\bar{\rho}}}\partial_3 \varpi|^2 {\mm{d}x}}}
\mbox{ and } a_3:=\sup_{0\neq\phi\in H_0^1(0,h)}\frac{\int^{h}_{0}\bar{\rho}' \phi^{ 2} {\mm{d}s}}{\|{\underline{\bar{\rho}} }\phi'\|^2_{L^2(0,h)}},
\end{align*}
then it holds that
\begin{align}
\label{2022309032042}
a_1=a_2=a_3={\varepsilon^2_{\mm{c}}}/g.
\end{align}
\end{lem}
\begin{pf} 
Recalling the definitions of ${\varepsilon_{\mm{c}}}$, $a_1$ and $a_2$, it is easy to see that
$$a_2\geqslant a_1\geqslant{\varepsilon^2_{\mm{c}}}/g.$$
 Hence, to obtain \eqref{2022309032042}, it suffices to prove that
\begin{align}
\label{lema92asfd}
a_2\leqslant a_3\leqslant {\varepsilon^2_{\mm{c}}}/g.
\end{align}

{We first verify that}
\begin{equation}\label{k2e}
a_2\leqslant a_3.
\end{equation}
Let $\phi(\xi_{\mm{h}},x_3):=\phi_1(\xi_{\mm{h}},x_3)+\mm{i}\phi_2(\xi_{\mm{h}},x_3):=\hat{\varpi}(\xi_{\mm{h}},x_3)$, where $\phi_1(\xi_{\mm{h}},x_3)$ and $\phi_2(\xi_{\mm{h}},x_3)$ are real functions. Obviously, $\phi_{i}(0)=\phi_{i}(h)=0$, where $i=1,2$ due to $\varpi|_{\partial\Omega}=0$.   Using both Fubini and  Parseval theorems, we have
\begin{align}
\int\bar{\rho}' \varpi^2 {\mm{d}x}&=\int_{\mathbb{R}^2} \int_{0}^h\bar{\rho}'|\phi(\xi_{\mm{h}},x_3)|^2\mm{d}x_{3}\mm{d}\xi_{\mm{h}}\nonumber\\
&=\int_{\mathbb{R}^2}\int_{0}^h\bar{\rho}'\big(|\phi_1(\xi_{\mm{h}},x_3)|^2+|\phi_2(\xi_{\mm{h}},x_3)|^2\big)\mm{d}x_{3}\mm{d}\xi_{\mm{h}}\label{lema95}
\end{align}
and
\begin{align}
\int |{\underline{\bar{\rho}} }\partial_{3}\varpi|^2 {\mm{d}x}&=\int_{\mathbb{R}^2}\int_{0}^h|{\underline{\bar{\rho}}} \partial_{3}\phi(\xi_{\mm{h}},x_3)|^2\mm{d}x_{3}\mm{d}\xi_{\mm{h}}\nonumber\\&=\int_{\mathbb{R}^2}
\int_{0}^h|{\underline{\bar{\rho}}}|^2\big(|\partial_{3}\phi_1(\xi_{\mm{h}},x_3)|^2
+|\partial_{3}\phi_2(\xi_{\mm{h}},x_3)|^2\big)\mm{d}x_{3}\mm{d}\xi_{\mm{h}}.\label{lema96}
\end{align}
Keeping in mind that $\phi_{i}\in H^1_0(0,h)$ if $\varpi\in H_0$, and then using \eqref{lema95} and  \eqref{lema96}, we  deduce that for any $\varpi\in H_0$,
\begin{align}
a_3&\geqslant\frac{\int_{\mathbb{R}^2} \int_{0}^h\bar{\rho}'|\phi(\xi_{\mm{h}},x_3)|^2\mm{d}x_{3}\mm{d}\xi_{\mm{h}}}
{\int_{\mathbb{R}^2}{\int_{0}^h|{\underline{\bar{\rho}}} \partial_{3}\phi(\xi_{\mm{h}},x_3)|^2\mm{d}x_{3}\mm{d}\xi_{\mm{h}}} } =\frac{\int\bar{\rho}' \varpi^2 {\mm{d}x}}{\int |{\underline{\bar{\rho}} }\partial_{3}\varpi|^2 {\mm{d}x}},\label{lema96w}
\end{align}
which implies \eqref{k2e}.
		
{Now we turn to verify that}	
\begin{equation}\label{k21}
a_3\leqslant {\varepsilon^2_{\mm{c}}}/g.
\end{equation}	
Recalling the definition of $a_3$, there is a  function $0\neq\phi_0\in H^1_0(0,h)$ such that	
\begin{equation}\label{k212}
\frac{\int^{h}_{0}\bar{\rho}' \phi_0^2 {\mm{d}x_3}}{\|{\underline{\bar{\rho}} }\phi_0'\|^2_{L^2(0,h)}}=a_3,
\end{equation}
please refer to \cite[Lemma 4.4]{JFJSARMA2019} for the proof.
Since $\phi_0\in H^1_0(0,h)$, there is a function sequence $\{\phi_{j}\}^{\infty}_{j=1}\subset C^{\infty}(0,h)$ satisfying $\phi_{j}\neq0$, $\phi_{j}(0)=\phi_{j}(h)=0$,
\begin{equation}\label{k213}
\int_{0}^h|{\underline{\bar{\rho}} } \phi'_{j}|^2\mm{d}s\rightarrow \int_{0}^h|\underline{\bar{\rho}} \phi'_{0}|^2\mm{d}s  \mbox{ and }  {\int_{0}^h \bar{\rho}' {\phi}^2_j\mm{d}s} \to {\int_{0}^h \bar{\rho}' \phi^2_0\mm{d}s} \mbox{ as } {j}\rightarrow \infty.
\end{equation}
We further define that
\begin{align}
& P_{k}(x_1,x_2):=k^{-1}e^{-\frac{x_1^2+x_2^2}{k}}\mbox{ and } H_{k}(x_1):= {k}^{-1/2}
\int_{0}^{x_1/\sqrt{k}}e^{-s^2}\mmd s , \label{k214a}
\\
& w^{k,j}(x)=\left( - \phi_j'(x_3)H_{k}(x_1)e^{-\frac{x^2_2}{k}} ,0,\phi_j(x_3)P_{k}(x_1,x_2) \right).
\end{align}  It is easy to  check that ${w}^{k,j} \in C^{\infty}(\Omega)\cap H_\sigma$ and
\begin{align}
&\int_{\mathbb{R}^2}\left| P_{k}(x_1,x_2)\right| ^2\mmd x_{\mm{h}}=\frac{\pi}{2k},\ \partial_{1}H_{k}(x_1) =k^{-1}e^{-\frac{x_1^2}{k}},\label{k214}\\&\int_{\mathbb{R}^2}\left|\partial_{1}P_k(x_1,x_2) \right| ^2\mmd x_{\mm{h}}=\int_{\mathbb{R}^2}\left|\partial_{2}P_{k}(x_1,x_2)\right| ^2 \mmd x_{\mm{h}}=\frac{2\pi}{k^3}.\label{k215}
\end{align}	
Thanks to \eqref{k214} and \eqref{k215}, we further conclude  that
\begin{align}
&\frac{2k}{\pi} \int\bar{\rho}' (w^{k,j}_3)^2 {\mm{d}x}
= \int_0^h\bar{\rho}'|\phi_{j}|^2\mm{d}x_3, \label{k2141} \\
&\mbox{ and }\frac{2k}{\pi} \int |{\underline{\bar{\rho}}}\nabla w^{k,j}_3|^2 {\mm{d}x}= \|{\underline{\bar{\rho}}}\phi'_{j}\|^2_{L^2(0,h)}+8k^{-2}\|{\underline{\bar{\rho}}}\phi_{j}\|^2_{L^2(0,h)}. \label{k2142}
\end{align}
Finally, exploiting \eqref{k212}, \eqref{k2141}  and \eqref{k2142},  we have
\begin{align}
\label{k2143}
a_3 =\lim_{j\to \infty} \lim_{k\to \infty}  \frac{\int_0^h\bar{\rho}' |\phi_j|^2\mm{d}x_3   }{ \|{\underline{\bar{\rho}}}\phi'_j\|_{L^2(0,h)}^2+8k^{-2}\|{\underline{\bar{\rho}}}\phi_{j}\|^2_{L^2(0,h)} }=\lim_{j\to \infty} \lim_{k\to \infty}   \frac{\int\bar{\rho}' (w^{k,j}_3)^2 {\mm{d}x}}{\int |{\underline{\bar{\rho}}}\nabla  w^{k,j}_3|^2 {\mm{d}x}}\leqslant {\varepsilon^2_{\mm{c}}}/g,
\end{align}
which implies \eqref{k21}. This completes the proof of \eqref{lema92asfd}.
\hfill $\Box$	\end{pf}	
 \begin{lem}  \label{xfsddfs2212}
An elliptic estimate for the Dirichlet boundary value condition (referring to  \cite[Lemma A.7]{ZHAOYUI}):
Let $i\geqslant  0$, $ {f}^1\in H^{i} $ and $ {f}^2 \in H^{i+3/2}(\partial \Omega)$ be given, then there exists a unique solution $u\in H^{i+2} $, which
solves the boundary value problem 
\begin{equation*}
\begin{cases}
\Delta u =  {f}^1&\mbox{in }   \Omega, \\[1mm]
u= {f}^2  &\mbox{on }\partial   \Omega 
\end{cases}
\end{equation*}
and satisfies
\begin{equation}
\label{xfsdsaf41252}
\|u\|_{{i+2} } \lesssim
\| {f}^1\|_{i }+ \| {f}^2 \|_{H^{i+3/2}(\partial\Omega)} .
\end{equation}
\end{lem}
\begin{lem}  \label{xfs05072212}
An elliptic estimate for the Neumann boundary value condition  (referring  to \cite[Lemma A.5]{JiangZhangZhang}):
Let  $i\geqslant  0$,  $ {f} \in H^{i} $ and $f_3|_{\partial\Omega}=0$,  then there exists a unique solution $u\in H^{i+1} $, which
solves the boundary value problem
\begin{equation*}
\begin{cases}
 \Delta u =  \mm{div}{f}  &\mbox{in }   \Omega, \\[1mm]
\partial_3 u= 0  &\mbox{on }\partial   \Omega
\end{cases}
\end{equation*} and satisfies 
\begin{equation}
\label{xfsddfsf20170saf5141252}
\|\nabla u\|_{ {i} } \lesssim\|\mm{div}{f}  \|_{i-1 }+\|{f}  \|_{i } .
\end{equation}
\end{lem}
\begin{lem}  \label{xfs05072212a}
Let $a=\theta+2/3$, $s=1-\theta$  and $\theta\in(0,({19-\sqrt{349}})/{6})$, we have
\begin{align}
&{a(4+s)}>{3s},\label{a1sqa5}\\
&{3a}>2,\label{a1sqa1}\\
&{a(9+4s)} >{4(1+s)},\label{a1sqa2}\\
& a(2+s)>{1+s},\label{a1sqa4}\\
&{3a(1+2s)}>{(2+s)(1+s)} .\label{a1sqa3}
\end{align}
\end{lem}
\begin{pf}  
(1) We can compute out that
\begin{align*}\frac{a(4+s)}{3s}-1= \frac{\left(\theta + {2}/{3}\right)(5 - \theta)}{3(1 - \theta)} -1=\frac{\theta(10 - 3\theta)}{3(2-\theta)},\end{align*}
which, together with the fact that $10 - 3\theta> 0$ and $2-\theta>0$, yields \eqref{a1sqa5}. 

(2) Obviously
$$ {3a}-{2} = {3} \left( \theta +  {2}/{3} \right) = {3}\theta>0.$$
Hence \eqref{a1sqa1} holds.

(3) It holds that
\begin{align*}\frac{a(9+4s)}{4(1+s)}-1=&\frac{4\theta^2-({26}+{31}\theta)/3 }{4\theta-8}-1
=\frac{2+ 43\theta-12\theta^2}{12(2-\theta)}. \end{align*}
In addition, it is easy to check that $(19 - \sqrt{349})/6\approx 0.0531$ and the numerator $2+ 43\theta-12\theta^2$ is positive over $(0,({19-\sqrt{349}})/{6})$. Thus, we immediately get \eqref{a1sqa2}. 

(4) Since $10 - 3\theta> 0$ and $2-\theta>0$, we can obtain that
\begin{align*}a+\frac{a}{1+s}-1&= \frac{2 - \theta^2 +  {7}\theta/3}{2 - \theta}-1=\frac{\theta(10 - 3\theta)}{3(2-\theta)}>0. \end{align*}
Therefore \eqref{a1sqa4} holds.

(5) Noting that $10 - 7\theta> 0$ and $(3-\theta)(2-\theta)>0$, thus
\begin{align*}\frac{3a(1+2s)}{(2+s)(1+s)}-1=\frac{6+5\theta-6\theta^2}{(3-\theta)(2-\theta)}-1=\frac{10\theta - 7\theta^2}{(3-\theta)(2-\theta)}>0,
\end{align*}
 which yields \eqref{a1sqa3}. 
The proof is complete. \hfill $\Box$\end{pf}

\vspace{4mm}
\noindent\textbf{Acknowledgements.}
The research of Fei Jiang was supported by the NSFC (Grant Nos. 12525109, 12371233 and 12231016), the NSF of Fujian Province of China (Grant  Nos. 2024J011011 and 2022J01105), Fujian Alliance of Mathematics (Grant No. 2025SXLMMS01) and the Central Guidance on Local Science and Technology Development Fund of Fujian Province (Grant No. 2023L3003), and the research of Zhipeng Zhang by NSFC (Grant Nos. 12101305, 12471215 and 12331007), and the research of Youyi Zhao was supported by NSFC (Grant No. 12401289) and the Natural Science Foundation of Fujian Province of China (Grant No. 2024J08029). The second author would like to take this opportunity to thank Prof. Song Jiang for his encouragement and support.

	
\renewcommand\refname{References}
\renewenvironment{thebibliography}[1]{%
\section*{\refname}
\list{{\arabic{enumi}}}{\def\makelabel##1{\hss{##1}}\topsep=0mm
\parsep=0mm
\partopsep=0mm\itemsep=0mm
\labelsep=1ex\itemindent=0mm
\settowidth\labelwidth{\small[#1]}
\leftmargin\labelwidth \advance\leftmargin\labelsep
\advance\leftmargin -\itemindent
\usecounter{enumi}}\small
\def\newblock{\ }
\sloppy\clubpenalty4000\widowpenalty4000
\sfcode`\.=1000\relax}{\endlist}
\bibliographystyle{model1b-num-names}
\bibliography{refs}
\end{CJK*}
\end{document}